\documentclass[a4paper]{amsart}

\usepackage[ps2pdf]{hyperref}
\usepackage{amsmath,amssymb}
\usepackage{bbold}
\usepackage{psfrag}
\usepackage{graphicx}
\usepackage{subfigure}
\usepackage{leftidx}
\usepackage[all]{xy}
\usepackage[latin1]{inputenc}
\usepackage{enumerate}

\newtheorem{thm}{Theorem}[section]
\newtheorem{cor}[thm]{Corollary}
\newtheorem{lem}[thm]{Lemma}
\newtheorem{prop}[thm]{Proposition}

\theoremstyle{definition}

\newtheorem{rem}[thm]{Remark}
\newtheorem{exa}[thm]{Example}

\DeclareMathOperator{\cdiv}{{\div \hspace*{-.2em}\raisebox{.15\height}{\scalebox{1}[.573]{$\mid$}}}}

\newcommand{\biMod}[2]{\leftidx{_{#1}}{\mathrm{Mod}}{_{#2}}}
\newcommand{\bimod}[2]{\leftidx{_{#1}}{\mathrm{mod}}{_{#2}}}

\newcommand{\bi}[3]{\leftidx{_{#1}}{#2}{_{#3}}}

\newcommand{\lComod}[1]{\leftidx{_{#1}}{\mathrm{Comod}}{}}
\newcommand{\lsComod}[2]{\leftidx{_{#2}}{\mathrm{Comod}}{^{#1}}}
\newcommand{\ct}[2]{\leftidx{}{\otimes}{^{#1}_{#2}}}

\newcommand{\lMod}[1]{\leftidx{_{#1}}{\mathrm{Mod}}{}}
\newcommand{\lsMod}[2]{\leftidx{_{#1}}{\mathrm{Mod}}{_{#2}}}
\newcommand{\lsmod}[2]{\leftidx{_{#1}}{\mathrm{mod}}{_{#2}}}

\newcommand{\BiAlg}{\mathrm{BiAlg}}
\newcommand{\BiMon}{\mathrm{BiMon}}

\newcommand{\HopfAlg}{\mathrm{HopfAlg}}
\newcommand{\HopfMon}{\mathrm{HopfMon}}

\newcommand{\biotimes}{\boxtimes}

\newcommand{\Cat}{\mathrm{Cat}}
\newcommand{\MonCat}{\mathrm{MonCat}}

\newcommand{\blabla}[1]{\quad\text{#1}\quad}
\newcommand{\co}{\colon}
\newcommand{\id}{\mathrm{id}}
\newcommand{\kk}{\Bbbk}
\newcommand{\kt}{$\Bbbk$\nobreakdash-\hspace{0pt}}

\newcommand{\opp}{\mathrm{op}}
\newcommand{\cop}{\mathrm{cop}}

\newcommand{\cc}{\mathcal{C}}
\newcommand{\bb}{\mathcal{B}}

\newcommand{\dd}{\mathcal{D}}
\newcommand{\aaa}{\mathcal{A}}
\newcommand{\ee}{\mathcal{E}}
\newcommand{\ff}{\mathcal{F}}

\newcommand{\zz}{\mathcal{Z}}

\newcommand{\rr}{\mathcal{R}}
\newcommand{\AT}{\mathbb{A}}

\newcommand{\lmo}[2]{\leftidx{_{#1}}{#2}{}}

\newcommand{\rhom}[1]{[#1]^r}
\newcommand{\lhom}[1]{[#1]^l}

\newcommand{\iso}{\stackrel{\sim}{\longrightarrow}}
\newcommand{\shortiso}{\stackrel{\scriptscriptstyle\sim\phantom{i}}{\rightarrow}}

\newcommand{\HM}{\mathcal{H}}
\newcommand{\mo}[2]{{#2}^{#1}}
\newcommand{\como}[2]{{#2}_{#1}}

\newcommand{\ti}{\text{-}\,}
\newcommand{\un}{\mathbb{1}}
\newcommand{\cp}{\rtimes}

\newcommand{\Ob}{\mathrm{Ob}}

\newcommand{\End}{\mathrm{End}}
\newcommand{\Hom}{\mathrm{Hom}}
\newcommand{\Nat}{{\textsc{Hom}}}
\newcommand{\EndNat}{{\textsc{End}}}

\newcommand{\Vect}{{\mathrm{Vect}}}
\newcommand{\vect}{\mathrm{vect}}

\newcommand{\lev}{\mathrm{ev}}
\newcommand{\rev}{\widetilde{\mathrm{ev}}}
\newcommand{\lcoev}{\mathrm{coev}}
\newcommand{\rcoev}{\widetilde{\mathrm{coev}}}

\newcommand{\leva}[2]{\mathrm{ev}^{#1}_{#2}}
\newcommand{\reva}[2]{\widetilde{\mathrm{ev}}^{#1}_{#2}}
\newcommand{\lcoeva}[2]{\mathrm{coev}^{#1}_{#2}}
\newcommand{\rcoeva}[2]{\widetilde{\mathrm{coev}}^{#1}_{#2}}

\newcommand{\ldual}[1]{\leftidx{^\vee}{\!#1}{}}
\newcommand{\rdual}[1]{{#1}^\vee}

\newcommand{\adjunct}[2]{\!\!\raisebox{.6ex}{\xymatrix{\ar@/^.4pc/[r]^{#1}  &  \ar@/^.4pc/[l]^{#2}}}\!\!}
\newcommand{\rsdraw}[3]{\raisebox{-#1\height}{\scalebox{#2}{\includegraphics{#3.eps}}}}
\providecommand{\bysame}{\leavevmode\hbox to3em{\hrulefill}\thinspace}

\newcommand{\uHopf}{pre-Hopf}

\newcommand{\frob}{Hopf}

\newcommand{\ufrob}{pre-Hopf}

\newcommand{\FO}{{\mathbb H}}
\newcommand{\FOli}{{\mathbb H}^{l\stackrel{-1}{\phantom{.}}}}
\newcommand{\FOri}{{\mathbb H}^{r\stackrel{-1}{\phantom{.}}}}
\newcommand{\Hli}{H^{l\stackrel{-1}{\phantom{.}}}}
\newcommand{\uei}{u^{e-1}}
\newcommand{\Hri}{H^{r\stackrel{-1}{\phantom{.}}}}

\newcommand{\HAr}{\mathbb{H}'}
\newcommand{\HAl}{\mathbb{H}}
\newcommand{\HAri}{{\mathbb{H}'}^{-1}}
\newcommand{\HAli}{\mathbb{H}^{-1}}

\newcommand{\zzl}{\zz^{\mathrm{lax}}}

\newcommand{\Mon}{\mathrm{Mon}}

\newcommand{\Zed}{\mathbb{Z}}
\newcommand{\NN}{\mathbb{N}}

\newcommand{\prettydef}[1]{\left\{\begin{array}{ccl}#1\end{array}\right.}

\newcommand{\dar}[2]{\ar@<2pt>[r]^-{#1}\ar@<-2pt>[r]_-{#2}}

\newcommand{\labela}{\renewcommand{\labelenumi}{{\rm (\alph{enumi})}}}
\newcommand{\labeli}{\renewcommand{\labelenumi}{{\rm (\roman{enumi})}}}

\begin{document}
\title{Hopf monads on monoidal categories}
\author[A. Brugui\`eres]{Alain Brugui\`eres}
\author[S. Lack]{Steve Lack}
\author[A. Virelizier]{Alexis Virelizier}
\email{bruguier@math.univ-montp2.fr, s.lack@uws.edu.au, virelizi@math.univ-montp2.fr}

\date{\today}

\begin{abstract}
We define Hopf monads on an arbitrary monoidal category, extending the definition given in~\cite{BV2} for monoidal categories with duals. A Hopf monad is a bimonad (or opmonoidal monad) whose fusion operators are invertible. This definition can be formulated in terms of Hopf adjunctions, which are comonoidal adjunctions with an invertibility condition.
On a monoidal category with internal Homs, a Hopf monad is a bimonad admitting a left and a right antipode.

Hopf monads generalize Hopf algebras to the non-braided setting. They also generalize Hopf algebroids (which are linear Hopf monads on a category of bimodules admitting a right adjoint). We show that  any finite tensor category is the category of finite-dimensional modules over a Hopf algebroid.

Any Hopf algebra in the center of a monoidal category $\cc$ gives rise to a Hopf monad on
$\cc$. The Hopf monads so obtained are exactly the augmented Hopf monads. More generally if a Hopf monad $T$ is a retract of a Hopf monad $P$, then $P$ is a cross product of $T$ by a Hopf algebra of the center of the category of $T$\ti modules (generalizing the Radford-Majid bosonization of Hopf algebras).

We show that the comonoidal comonad of a Hopf adjunction is canonically represented by a cocommutative central coalgebra. As a corollary, we obtain an extension of Sweedler's Hopf module decomposition theorem to Hopf monads (in fact to the weaker notion of \uHopf\ monad).
\end{abstract}
\maketitle

\setcounter{tocdepth}{1} \tableofcontents

\section*{Introduction}

Hopf monads on autonomous categories (that is, monoidal categories with duals) were introduced in \cite{BV2} as a tool for understanding and comparing quantum invariants of $3$ manifolds, namely the Reshetikhin-Turaev invariant associated with a modular category and the Turaev-Viro invariant associated with a spherical category (as revisited by Barrett-Westbury).

In this paper we extend the notion of Hopf monad to any monoidal category. 
Hopf monads generalize classical Hopf algebras, as well as Hopf algebras in a braided category.
Hopf algebras are bialgebras with an extra condition: the existence of an invertible antipode. 
Similarly,  one expects Hopf monads to be bimonads satisfying some extra condition.

The concept of bimonad (also called opmonoidal monad) was introduced by Moerdijk  in
\cite{Moer}\footnote{Bimonads were introduced in \cite{Moer} under the name `Hopf monads', which we prefer to
reserve for bimonads with antipodes by analogy with Hopf algebras.}.
Recall that if $T$ is a monad on a category $\cc$, then one defines a category $\mo{T}{\cc}$ of $T$\ti modules in $\cc$ (often called $T$\ti algebras). A \emph{bimonad} on a monoidal category $\cc$ is a monad on $\cc$ such that $\mo{T}{\cc}$ is monoidal and the forgetful functor $U_T\co \mo{T}{\cc}\to \cc$ is strict monoidal. This means that $T$ is a comonoidal monad: it comes with a coassociative natural transformation $T_2(X,Y) \co T(X\otimes Y) \to TX \otimes TY$ and a counit $T_0 \co T\un \to \un$.  
For example, a bialgebra $A$ in a braided category $\bb$ gives rise to bimonads  $A \otimes ?$ and $? \otimes A$ on $\bb$. More generally, bialgebroids in the sense of Takeuchi are also examples of bimonads. More generally still, any comonoidal adjunction defines a bimonad,
so that bimonads exist in many settings.

The `extra condition' a bimonad should satisfy in order to deserve the title of Hopf monad is not obvious, as there is no straightforward generalization of the notion of antipode to the monoidal setting. When $\cc$ is autonomous, according to Tannaka theory, one expects that a bimonad $T$ be Hopf if and only if $\mo{T}{\cc}$ is autonomous. This turns out to be equivalent to the existence of a left antipode and a right antipode, which are natural transformations $s^l_X \co T(\ldual{T(X)}) \to \ldual{X}$ and $s^r_X \co T(\rdual{T(X)}) \to \rdual{X}$. That was precisely the definition of a Hopf monad given in~\cite{BV2}. While it is satisfactory for applications to quantum topology, as the categories involved are autonomous, this definition has some drawbacks for other applications: for instance, it doesn't encompass infinite-dimensional Hopf algebras since the category of vector spaces of arbitrary dimension is not autonomous. Therefore one is prompted to ask several questions:
\begin{itemize}
\item What are Hopf monads on arbitrary monoidal categories?
\item What are Hopf monads on closed monoidal categories (with internal Homs)?
\item Is it possible to characterize Hopf monads obtained from Hopf algebras?
\item Can one extend classical results of the theory  of Hopf algebras to Hopf monads on monoidal categories?
\item When does a bialgebroid define a Hopf monad?
\end{itemize}
The aim of this paper is to answer these questions.

In Section~\ref{sect-HopfMon}, we define Hopf monads on an arbitrary monoidal category. Our definition is inspired by the fact that a bialgebra $A$ is a Hopf algebra if and only if its fusion morphisms $H^l, H^r \co A \otimes A \to A \otimes A$, defined by $H^l(x\otimes y)=x^{(1)} \otimes x^{(2)}y$ and $H^r(x \otimes y)=x^{(2)}y \otimes x^{(1)}$, are invertible. If $T$ is a bimonad, we introduce the \emph{fusion operators} $H^l$ and $H^r$, which are natural transformations
\begin{align*}
H^l_{X,Y}=(TX \otimes \mu_Y)T_2(X,TY)&\co T(X \otimes TY) \to TX \otimes TY,\\
H^r_{X,Y}=(\mu_X \otimes TY)T_2(TX,Y)&\co T(TX \otimes Y) \to TX \otimes TY,
\end{align*}
and decree that $T$ is a \emph{Hopf monad} if $H^l$ and $H^r$ are invertible.
We also introduce the related notion of \emph{\frob\ adjunction}. The monad of a \frob\ adjunction is a Hopf monad, and a bimonad is a Hopf monad if and only if its adjunction is a \frob\ adjunction. It turns out that certain classical results on Hopf algebras extend naturally to Hopf monads (or more generally to pre-Hopf monads), such as Maschke's semisimplicity criterion and Sweedler's theorem on the structure of Hopf modules (see Section~\ref{sect-Hmod}).

In Section~\ref{sect-HM-closed}, we study Hopf monads on closed monoidal categories. Hopf monads on such categories can be characterized, as in the autonomous case, in categorical terms and also in terms of antipodes. More precisely, let $T$ be a bimonad on a closed monoidal category $\cc$, that is, a monoidal category with internal Homs. We show that $T$  is a Hopf monad if and only if its category of modules  $\mo{T}{\cc}$ is closed
and the forgetful functor $U_T$ preserves internal Homs. Also $T$ is a Hopf monad if and only if it admits a \emph{left antipode} and \emph{right antipode}, that is, natural transformations in two variables:
$$
s^l_{X,Y} \co T\lhom{TX,Y} \to \lhom{X,TY} \quad \text{and} \quad
s^r_{X,Y} \co T\rhom{TX,Y} \to \rhom{X,TY}
$$
where $\lhom{-,-}$ and $\rhom{-,-}$ denote the left and right internal Homs, each of them satisfying two axioms as expected. The proof of these results relies on a classification of adjunction liftings.
In the special case where $\cc$ is autonomous, we show that
the definition of a Hopf monad given in this paper specializes to the one given in \cite{BV2}.  In the special case where $\cc$ is $*$-autonomous (and so monoidal closed), a Hopf monad in the sense of \cite{DS3} is a Hopf monad in our sense (but the converse is not true).


In Section~\ref{sect-rep}, we study the relations between Hopf algebras and Hopf monads.  Given a lax central bialgebra of a monoidal category $\cc$, that is, a bialgebra $A$ in the lax center $\zzl(\cc)$ of  $\cc$, with
lax\footnote{Here \emph{lax} means that $\sigma$ need not be invertible.}
half-braiding $\sigma \co A \otimes ? \to ? \otimes A$,   the endofunctor $A \otimes ?$ of $\cc$ is a bimonad, denoted by $A \otimes_\sigma ?$ on $\cc$. This bimonad is augmented, that is, endowed with a bimonad morphism $A \otimes_\sigma ? \to 1_\cc$. It is a Hopf monad if and only if $A$ is a Hopf algebra in the center $\zz(\cc)$ of $\cc$.
The main result of the section is that this construction defines an equivalence of categories between  central Hopf algebras of $\cc$  (that is,   Hopf algebras in the center $\zz(\cc)$) and augmented Hopf monads on~$\cc$.
More generally, given a Hopf monad $T$ on $\cc$ and a  central Hopf algebra $(\AT,\sigma)$ of
the category of $T$\ti modules,    we construct a Hopf monad $\AT \cp_\sigma T$ on $\cc$
of which $T$ is a retract. Conversely,
under suitable exactness conditions (involving reflexive coequalizers), any Hopf monad $P$ of which $T$ is a retract is of the form $\AT \cp_\sigma T$. The proof of this result is based on  two general constructions involving Hopf monads: the cross product and the cross quotient,  which are studied in Section~\ref{sect-crosstruc}.

In Section~\ref{sect-Hmod}, we show that the comonoidal comonad of a \uHopf\ adjunction is canonically represented by a cocommutative central coalgebra.
Combining this with a descent result for monads, we obtain a generalization of Sweedler's Hopf module decomposition theorem to Hopf monads (in fact to pre-Hopf monads).  We study the close relationships between Hopf adjunctions, Hopf monads, and cocommutative central coalgebras.

 Finally, in Section~\ref{sect-hopf-oid}, we study  bialgebroids which, according to Szlach\'anyi~\cite{Szl03}, are linear bimonads on categories of bimodules admitting a right adjoint. A bialgebroid corresponds with a Hopf monad if and only if it is a Hopf algebroid in the sense of Schauenburg~\cite{Sch}. We also use Hopf monads to prove that any finite tensor category is naturally equivalent (as a tensor category) to the category of finite-dimensional modules over some finite dimensional Hopf algebroid.


%
%
%
%
%
%

\section{Preliminaries and notations}\label{sect-prelims}

 Unless otherwise specified, categories are small, and monoidal categories
are strict. We denote by $\Cat$ the category of small categories (which is not small).

If $\cc$ is a category, we denote by $\Ob(\cc)$ the set of objects of $\cc$ and by $\Hom_\cc(X,Y)$ the set of morphisms in
$\cc$ from an object $X$ to an object $Y$. The identity functor of~$\cc$ is denoted by $1_\cc$.

If $\cc$ is a category and $c$ an object of $\cc$, the \emph{category of objects of $\cc$ over $c$} is the category  $\cc/c$  whose objects are pairs $(a,\phi)$, with $a \in \Ob(\cc)$ and $\phi \in \Hom_\cc(a,c)$. Morphisms from $(a,\phi)$ to  $(b,\psi)$ in $\cc/c$ are morphism $f\co a \to b$ in $\cc$  satisfying the condition $\psi f= \phi$. They are called \emph{morphisms over $c$}.

Similarly the \emph{category of objects of $\cc$ under $c$} is the category $c \backslash \cc$ whose objects are pairs $(a,\phi)$, with $a \in \Ob(\cc)$ and $\phi \in \Hom_\cc(c,a)$.

A pair of parallel morphisms $$X\xymatrix{\ar@<2.5pt>[r]^f\ar@<0pt>[r]_g &}Y$$ is \emph{reflexive} (resp.\@ \emph{coreflexive}) if $f$ and $g$ have a common section (resp.\@ a common retraction), that is, if there exists a morphism $h \co Y \to X$ such that $fh=gh=\id_Y$ (resp.\@ $hf=hg=\id_X$).
A \emph{reflexive coequalizer} is a coequalizer of a reflexive pair. Similarly a \emph{coreflexive equalizer} is an equalizer of a coreflexive pair.

\subsection{Monoidal categories and functors}\label{sect-monofunctor}
Given an object $X$ of a monoidal category $\cc$, we denote by $X \otimes ?$ the endofunctor of $\cc$ defined on objects by $Y \mapsto X \otimes Y$ and on morphisms by $f \mapsto X \otimes f=\id_X \otimes f$. Similarly one defines the endofunctor $? \otimes X$ of $\cc$.

Let $(\cc,\otimes,\un)$ and $(\dd, \otimes, \un)$ be two monoidal categories.
A \emph{monoidal functor} from $\cc$ to $\dd$ is a triple
$(F,F_2,F_0)$, where $F\co \cc \to \dd$ is a functor, $F_2\co
F\otimes F \to F \otimes$ is a natural transformation, and $F_0\co\un
\to F(\un)$ is a morphism in~$\dd$, such that:
\begin{align*}
& F_2(X,Y \otimes Z) (\id_{F(X)} \otimes F_2(Y,Z))= F_2(X \otimes Y, Z)(F_2(X,Y) \otimes \id_{F(Z)}) ;\\
& F_2(X,\un)(\id_{F(X)} \otimes F_0)=\id_{F(X)}=F_2(\un,X)(F_0
\otimes \id_{F(X)}) ;
\end{align*}
for all objects $X,Y,Z$ of $\cc$.

A monoidal functor $(F,F_2,F_0)$ is said to be \emph{strong} (resp.\@ \emph{strict}) if $F_2$ and $F_0$ are
isomorphisms (resp.\@ identities).

A natural transformation $\varphi\co F \to G$ between monoidal functors is \emph{monoidal}
if it satisfies:
\begin{equation*}
\varphi_{X \otimes Y} F_2(X,Y)= G_2(X,Y) (\varphi_X \otimes
\varphi_Y) \quad \text{and} \quad G_0=\varphi_\un F_0.
\end{equation*}

We denote by $\MonCat$ the category of small monoidal categories, morphisms being strong monoidal functors.

\subsection{Comonoidal functors}\label{sect-comonofunctor}
Let $(\cc,\otimes,\un)$ and $(\dd, \otimes, \un)$ be two monoidal categories.
A \emph{comonoidal functor} (also called \emph{opmonoidal functor}) from $\cc$ to
$\dd$ is a triple $(F,F_2,F_0)$, where $F\co \cc \to \dd$ is a functor, $F_2\co F \otimes \to F\otimes F$ is a natural
transformation, and $F_0\co F(\un) \to \un$ is a morphism in $\dd$, such that:
\begin{align*}
& \bigl(\id_{F(X)} \otimes F_2(Y,Z)\bigr) F_2(X,Y \otimes Z)= \bigl(F_2(X,Y) \otimes \id_{F(Z)}\bigr) F_2(X \otimes Y, Z) ;\\
& (\id_{F(X)} \otimes F_0) F_2(X,\un)=\id_{F(X)}=(F_0 \otimes \id_{F(X)}) F_2(\un,X) ;
\end{align*}
for all objects $X,Y,Z$ of $\cc$.

A comonoidal functor $(F,F_2,F_0)$ is said to be \emph{strong} (resp.\@ \emph{strict}) if $F_2$ and $F_0$ are
isomorphisms (resp.\@ identities). In that case, $(F,F^{-1}_2,F^{-1}_0)$ is a strong (resp. strict) monoidal functor.

A natural transformation $\varphi\co F \to G$ between monoidal functors is \emph{comonoidal} if it satisfies:
\begin{equation*}
G_2(X,Y) \varphi_{X \otimes Y}= (\varphi_X \otimes \varphi_Y) F_2(X,Y)\quad \text{and} \quad G_0 \varphi_\un= F_0.
\end{equation*}

Note that the notions of comonoidal functor and comonoidal natural transformation are dual to the notions of monoidal functor and monoidal natural transformation.

\section{Hopf monads}\label{sect-HopfMon}

In this section, we define Hopf monads on an arbitrary monoidal category: they are the bimonads whose fusion operators are invertible. We also introduce the related notion of \frob\ adjunction: the monad of a \frob\ adjunction is a Hopf monad, and a bimonad is a Hopf monad if and only if its adjunction is a \frob\ adjunction.

\subsection{Monads}\label{monad}
Let $\cc$ be a category. Recall that the category $\End(\cc)$ of endofunctors of $\cc$ is strict monoidal with
composition for monoidal product and identity functor $1_\cc$ for unit object. A \emph{monad} on $\cc$ is an algebra in
$\End(\cc)$, that is, a triple $(T,\mu,\eta)$, where $T\co \cc \to \cc$ is a functor, $\mu\co T^2 \to T$ and $\eta\co
1_\cc \to T$ are natural transformations, such that:
\begin{equation*}
\mu_X T(\mu_X)=\mu_X\mu_{TX} \quad \text{and} \quad
\mu_X\eta_{TX}=\id_{TX}=\mu_X T(\eta_X)
\end{equation*}
for any object $X$ of $\cc$.

Monads on $\cc$ form a category $\Mon(\cc)$, a morphism from a monad $(T,\mu,\eta)$ to a monad $(T',\mu',\eta')$ being a natural transformation $f \co T \to T'$
such that $f\eta=\eta'$ and $f \mu=\mu'T(f)f_T$. The identity functor $1_\cc$ is a monad (with the identity for product and unit) and it is an initial object in
$\Mon(\cc)$.


\subsection{Modules over a monad}
Let $(T,\mu,\eta)$ be a monad on a category $\cc$. An \emph{action} of $T$ on an object $M$ of $\cc$ is a morphism
$r\co T(M) \to M$ in $\cc$ such that:
\begin{equation*}
r T(r)= r \mu_M \quad \text{and} \quad r \eta_M= \id_M.
\end{equation*}
The pair $(M,r)$ is then called a \emph{$T$\ti module in $\cc$}, or just a \emph{$T$-module}\footnote{Pairs $(M,r)$ are
usually called $T$-algebras in the literature (see \cite{ML1}).}.

Given two $T$-modules $(M,r)$ and $(N,s)$ in $\cc$, a \emph{morphism of $T$\ti modules} from $(M,r)$ to $(N,r)$ is a
morphism $f\in \Hom_\cc(M,N)$ which is \emph{$T$-linear}, that is, such that $f r=s T(f)$. This gives rise to the
category of $T$-modules (in $\cc$), with composition inherited from $\cc$. We denote this category by $\mo{T}{\cc}$ (the notation $T \ti \cc$ is used in~\cite{BV2}) .

The \emph{forgetful functor  $U_T\co
\mo{T}{\cc} \to \cc$ of $T$} is defined by $U_T(M,r)=M$ for any $T$-module $(M,r)$ and $U_T(f)=f$
for any $T$-linear morphism~$f$. It has a left adjoint $F_T \co \cc \to \mo{T}{\cc}$, called the \emph{free module functor},
defined by $F_T(X)=(TX,\mu_X)$ for any object $X$ of $\cc$ and $F_T(f)=Tf$ for any  morphism $f$ of $\cc$.

\subsection{Monads, adjunctions, and monadicity}
Let $(F \co \cc \to \dd, U \co \dd \to \cc)$ be an adjunction, with unit $\eta \co 1_\cc \to UF$ and counit $\varepsilon \co FU \to 1_\dd$. Then $T=UF$ is a monad with product $\mu=U(\varepsilon_F)$ and unit $\eta$.
There exists a unique functor $K \co \dd \to \mo{T}{\cc}$ such that $U_TK=U$ and $KF=F_T$. This functor $K$,  called
the \emph{comparison functor} of the adjunction $(F,U)$, is defined by $K(d)=(Ud, U\varepsilon_d)$.

An adjunction $(F,U)$ is \emph{monadic} if its comparison functor $K$ is an equivalence of categories.
For example, if $T$ is a monad on $\cc$, the adjunction $(F_T,U_T)$ has monad $T$ and comparison functor $K=1_{\mo{T}{\cc}}$, and so is monadic.

A functor $U$ is \emph{monadic} if it admits a left adjoint $F$ and the adjunction $(F,U)$ is monadic.
If such is the case, the monad $T=UF$ of the adjunction $(F,U)$ is called \emph{the
monad of $U$}. It is well-defined up to unique isomorphism of monads (as the left adjoint $F$ is unique up to unique
natural isomorphism).

\begin{thm}[Beck]\label{thm-beck} An adjunction $(F \co \cc \to \dd, U \co \dd \to \cc)$ is monadic if and only if the functor~$U$ satisfies the following conditions:
\begin{enumerate}
\labela
\item The functor $U$ is conservative, that is, $U$ reflects isomorphisms;
\item Any reflexive pair of morphisms in $\dd$ whose image by $U$ has a split coequalizer has a coequalizer, which is preserved by $U$.
\end{enumerate}
Moreover, if $(F,U)$ is monadic, the comparison functor $K$ is an isomorphism if and only if the functor $U$ satisfies the transport of structure condition:
\begin{enumerate}
\labela
\setcounter{enumi}{2}
\item For any isomorphism $f \co U(d) \to c$ in $\cc$, where $c \in \Ob(\cc)$ and $d \in \Ob(\dd)$, there exist a unique $\widetilde{c}\in \Ob(\dd)$ and a unique isomorphism $\widetilde{f}\co d \to \widetilde{c}$ in $\dd$ such that  $U(\widetilde{f})=f$.
\end{enumerate}
\end{thm}

\subsection{Bimonads}
A \emph{bimonad} on a monoidal category~$\cc$ is a monad
$(T,\mu,\eta)$ on $\cc$ such that the functor $T\co \cc \to \cc$ is comonoidal and the natural transformations $\mu\co
T^2 \to T$ and $\eta\co 1_\cc \to T$ are comonoidal. In other words, $T$ is endowed with a natural
transformation $T_2\co T \otimes \to T\otimes T$ and a morphism $T_0\co T(\un) \to \un$ in $\cc$ such that:
\begin{align*}
& \bigl(TX \otimes T_2(Y,Z)\bigr) T_2(X,Y \otimes Z)= \bigl(T_2(X,Y) \otimes TZ\bigr) T_2(X \otimes Y, Z) ,\\
& (TX\otimes T_0) T_2(X,\un)=\id_{TX}=(T_0 \otimes TX) T_2(\un,X),\\
& T_2(X,Y) \mu_{X \otimes Y}= (\mu_X \otimes \mu_Y) T_2(TX,TY) T(T_2(X,Y)),\\
& T_0 \mu_\un= T_0 T(T_0) , \qquad T_2(X,Y) \eta_{X \otimes Y}= \eta_X \otimes \eta_Y , \qquad T_0 \eta_\un= \id_\un.
\end{align*}


\begin{rem}\label{remoppositemonad}
A bimonad $T$ on a monoidal category $\cc=(\cc,\otimes,\un)$ may be viewed as a bimonad $T^\cop$ on the
monoidal category $\cc^{\otimes \opp}=(\cc,\otimes^\opp,\un)$, with comonoidal structure $T_2^\cop(X,Y)=T_2(Y,X)$ and $T_0^\cop=T_0$. The bimonad $T^\cop$  is called  the \emph{coopposite} of the bimonad $T$. We have: $\mo{T^\cop}{(\cc^{\otimes\opp})}=(\mo{T}{\cc})^{\otimes \opp}$.
\end{rem}

\begin{rem}\label{rembicomonad}
The dual notion of a bimonad is that of a bicomonad, that is, a monoidal comonad. An endofunctor $T$ of a monoidal category $\cc=(\cc,\otimes,\un)$ is a bicomonad if and only if the opposite endofunctor $T^\opp$ is a bimonad on  $\cc^\opp=(\cc^\opp,\otimes,\un)$.
\end{rem}

Bimonads on $\cc$ form a category $\BiMon(\cc)$, morphisms of bimonads being como\-noidal morphisms of monads. The identity functor $1_\cc$ is a bimonad on $\cc$, which is an initial object of
$\BiMon(\cc)$.

\subsection{Bimonads and comonoidal adjunctions}\label{sect-bimonadadj}

A \emph{comonoidal adjunction} is an adjunction
$(F\co \cc \to \dd, U\co \dd \to \cc)$, where $\cc$ and $\dd$ are monoidal categories, $F$ and $U$ are comonoidal functors, and the adjunction unit $\eta \co 1_\cc \to UF$ and counit $\varepsilon \co FU\to 1_\dd$ are comonoidal natural transformations.

If $(F,U)$ is a comonoidal adjunction, then $U$ is in fact a strong comonoidal functor, which we may view as a strong monoidal functor.
Conversely, if a strong monoidal functor $U \co \dd \to \cc$ admits a left adjoint $F$, then $F$ is comonoidal,
with comonoidal structure given by:
\begin{equation*}
F_2(X,Y)=\varepsilon_{FX \otimes FY} FU_2(FX,FY)F(\eta_X \otimes \eta_Y) \quad \text{and} \quad
F_0= \varepsilon_\un F(U_0),
\end{equation*}
and $(F,U)$ is a comonoidal adjunction (viewing $U$ as a strong comonoidal functor),  see \cite{MCC}.  A comonoidal adjunction is an instance of a doctrinal adjunction in the sense of \cite{Kel}.

The monad $T=UF$ of a comonoidal adjunction $(U,F)$ is a bimonad, and the comparison functor $K \co \dd \to \mo{T}{\cc}$ is strong monoidal and satisfies $U_T K=U$ as monoidal functors and  $KF=F_T$ as comonoidal functors (see  for instance  \cite[Theorem~2.6]{BV2}).

The comonad $\hat{T}=FU$ of a comonoidal adjunction $(U,F)$ is a \emph{comonoidal comonad}, that is, a comonad whose underlying endofunctor is endowed with a comonoidal structure so that its coproduct and counit are comonoidal.

\begin{exa}
The adjunction $(F_U,U_T)$ of a bimonad $T$ is a comonoidal adjunction (because $U_T$ is strong monoidal).
\end{exa}
\begin{rem} Comonoidal adjunctions are somewhat misleadingly called monoidal adjunctions in~\cite{BV2}. \end{rem}

\subsection{Fusion operators}

Let $T$ be a bimonad on a monoidal category $\cc$. The \emph{left fusion operator of $T$} is the natural transformation $H^l \co T(1_\cc \otimes T) \to T \otimes T$
defined by:
$$H^l_{X,Y}=(TX \otimes \mu_Y)T_2(X,TY)\co T(X\otimes TY) \to TX \otimes TY.$$
The \emph{right fusion operator of $T$} is the natural transformation $H^r \co T(T \otimes 1_\cc) \to T \otimes T$ defined by:
$$H^r_{X,Y}=(\mu_X \otimes TY)T_2(TX,Y)\co T(TX\otimes Y) \to TX \otimes TY.$$

From the axioms of a bimonad, we easily deduce:
\begin{prop} \label{prop-fusion}
The left fusion operator $H^l$ of a bimonad $T$ satisfies:
\begin{align*}
& H^l_{X,Y}T(X \otimes \mu_Y)=(TX \otimes \mu_Y)H^l_{X,TY}, \\
& H^l_{X,Y}T(X \otimes \eta_Y)=T_2(X,Y), \qquad H^l_{X,Y} \eta_{X\otimes TY}=\eta_X \otimes TY,\\
& (T_2(X,Y) \otimes TZ)H^l_{X \otimes Y,Z}=(TX \otimes H^l_{Y,Z})T_2(X,Y \otimes TZ),\\
& (T_0 \otimes TX)H^l_{\un,X}=\mu_X,\qquad (TX\otimes T_0)H^l_{X,\un}=T(X\otimes T_0),
\end{align*}
and the \emph{left pentagon equation:}
$$ (TX \otimes H^l_{Y,Z}) H^l_{X, Y \otimes TZ} = (H^l_{X,Y} \otimes TZ)
H^l_{X \otimes TY, Z} T(X \otimes H^l_{Y,Z}).
$$
Similarly the right fusion operator $H^r$ of $T$ satisfies:
\begin{align*}
& H^r_{X,Y}T(\mu_X \otimes Y)=(\mu_X \otimes TY)H^r_{TX,Y}, \\
& H^r_{X,Y}T(\eta_X \otimes Y)=T_2(X,Y), \qquad H^r_{X,Y} \eta_{TX\otimes Y}= TX \otimes \eta_Y,\\
& (TX \otimes T_2(Y,Z))H^r_{X, Y \otimes Z}=(H^r_{X,Y} \otimes TZ)T_2(TX \otimes Y,Z),\\
& (TX \otimes T_0)H^r_{X,\un}=\mu_X,\qquad (T_0\otimes TX)H^r_{\un,X}=T(T_0\otimes X),
\end{align*}
and the \emph{right pentagon equation:}
$$ (H^r_{X,Y} \otimes TZ) H^r_{TX\otimes Y,Z}= (TX \otimes H^r_{Y,Z})
H^r_{X, TY \otimes Z} T(H^l_{X,Y} \otimes Z).
$$
\end{prop}

\begin{rem}
A bimonad can be recovered from its left (or right) fusion operator. More precisely, let $T$ be an endofunctor of a monoidal category $\cc$ endowed with a natural transformation
$H_{X,Y} \co T(X \otimes TY) \to TX \otimes TY$ satisfying the left pentagon equation:
$$(TX \otimes H_{Y,Z}) H_{X, Y \otimes TZ} = (H_{X,Y} \otimes TZ)H_{X \otimes TY, Z} T(X \otimes H_{Y,Z}),$$
and with a morphism $T_0 \co T\un \to \un$ and a natural transformation $\eta_X \co X \to TX$ satisfying:
\begin{align*}
&H_{X,Y} \eta_{X\otimes TY}=\eta_X \otimes TY, && T_0\eta_\un=\id_\un,\\ 
&(TX\otimes T_0)H_{X,\un}=T(X\otimes T_0), && 
(T_0 \otimes TX)H_{\un,X}T(\eta_X)=\id_{TX}.
\end{align*}
Then $T$ admits a unique bimonad structure $(T,\mu,\eta,T_2,T_0)$ having left fusion operator $H$. The product $\mu$ and comonoidal structural morphism
$T_2$ are given by: $$\mu_X=( T_0 \otimes TX)H_{\un,X} \quad \text{and} \quad T_2(X,Y)=H_{X,Y}T(X \otimes \eta_Y).$$
\end{rem}

\subsection{Hopf monads and \uHopf\ monads}\label{sect-hopfmon}
Let $\cc$ be a monoidal category. A \emph{left} (resp.\@ a \emph{right}) \emph{Hopf monad} on $\cc$ is a bimonad on $\cc$ whose left
fusion operator $H^l$ (resp.\@ right fusion operator $H^r$) is an isomorphism.

A \emph{Hopf monad on $\cc$} is a bimonad on $\cc$ such that both left and right fusion operators are isomorphisms.
 Hopf monads on $\cc$ form a full subcategory $\HopfMon(\cc)$ of the category $\BiMon(\cc)$ of bimonads. The identity functor $1_\cc$ is a Hopf monad on $\cc$, which is an initial object of
$\HopfMon(\cc)$.

It is convenient to consider a weaker notion: a \emph{left} (resp.\@ \emph{right}) \emph{\uHopf\ monad} on $\cc$ is a bimonad on $\cc$ such that, for any object $X$ of $\cc$, the morphism $H^l_{\un,X}$ (resp.\@ $H^r_{X,\un}$) is invertible.

A \emph{\uHopf\ monad} is a bimonad which is a left and a right \uHopf\ monad.
Clearly any
Hopf monad is a \uHopf\ monad, but the converse is false:

\begin{exa}\label{exa-pre}
We provide an example of a \uHopf\ monad on a monoidal (even autonomous) category which is not a Hopf monad.
Let $\Zed\ti \vect_\kk$ be the autonomous category of finite dimensional $\Zed$\ti graded vector spaces on a field $\kk$, and let $\NN\ti\vect_\kk$ be its full subcategory of graded vector spaces with support in $\NN$.
The inclusion functor $\iota \co \NN\ti\vect_\kk \to \Zed\ti\vect_\kk$ has a left adjoint $\pi$, which sends a $\Zed$\ti graded vector space to its non-negative part. The adjunction $(\pi,\iota)$ is monoidal. The bimonad $T=\iota\pi$ on $\Zed\ti\vect_\kk$ of this adjunction (see Section~\ref{sect-bimonadadj}) is a \uHopf\ monad  but not a Hopf monad.
\end{exa}

\begin{rem}
Certain general results on Hopf algebras extend naturally to \uHopf\ monads, such as Sweedler's theorem on the structure of Hopf modules (see Section~\ref{sect-Hmod}). Also, Maschke's semisimplicity theorem for Hopf monads on autonomous categories given in \cite[Theorem 6.5]{BV2} holds  word for word for \uHopf\ monads in arbitrary monoidal categories. Indeed the proof given in \cite{BV2},
which relied on the properties of a certain  natural transformation $\Gamma_X\co X \otimes T\un \to T^2X$, extends in a straightforward way, observing that $\Gamma_X=\Hri_{X,\un}(\eta_X \otimes T\un)$.
\end{rem}

\begin{exa} \label{exa-rep-braided}
Given a Hopf algebra $A$ in a braided category, we depict its product $m$, unit $u$,
coproduct $\Delta$, counit $\varepsilon$, and invertible antipode $S$ as follows:
\begin{center}
\psfrag{A}[Bc][Bc]{\scalebox{.8}{$A$}} $m=$\rsdraw{.45}{.9}{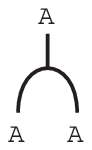}, \quad $u=$\rsdraw{.45}{.9}{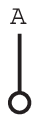}, \quad
$\Delta=$\rsdraw{.45}{.9}{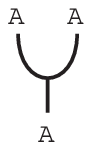}, \quad $\varepsilon=$\rsdraw{.45}{.9}{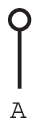}, \quad $S=\,\rsdraw{.45}{.9}{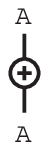}$\,,
\quad $S^{-1}=\,\rsdraw{.45}{.9}{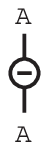}$\,.
\end{center}
Let $\bb$ be a braided category with braiding $\tau$, and $A$ a bialgebra in $\bb$. As shown in~\cite{BV2}, the endofunctor $A \otimes ?$ of $\bb$ is a bimonad on $\bb$, with structure maps:
\begin{center}
\psfrag{A}[Bc][Bc]{\scalebox{.8}{$A$}} \psfrag{X}[Bc][Bc]{\scalebox{.8}{$X$}} \psfrag{Y}[Bc][Bc]{\scalebox{.8}{$Y$}}
  $\mu_X=\rsdraw{.45}{.9}{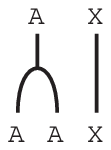}$, \quad $\eta_X=\rsdraw{.45}{.9}{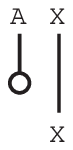}$, \quad $(A \otimes ?)_2(X,Y)=\rsdraw{.45}{.9}{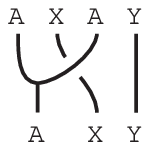}$, \quad
$(A \otimes ?)_0=\rsdraw{.45}{.9}{epsA}$\,.
\end{center}
Its fusion operators are:
\begin{center}
\psfrag{A}[Bc][Bc]{\scalebox{.8}{$A$}} \psfrag{X}[Bc][Bc]{\scalebox{.8}{$X$}} \psfrag{Y}[Bc][Bc]{\scalebox{.8}{$Y$}}
  $H^l_{X,Y}=\rsdraw{.45}{.9}{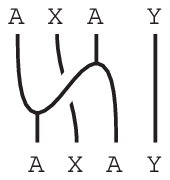}$\quad \text{and} \quad $H^r_{X,Y}=\rsdraw{.45}{.9}{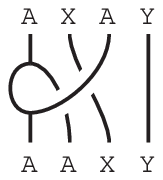}$.
\end{center}
If $A$ is a Hopf algebra with invertible antipode $S$, then $A \otimes ?$ is a Hopf monad, the inverses of the fusion operators being:
\begin{center}
\psfrag{A}[Bc][Bc]{\scalebox{.8}{$A$}} \psfrag{X}[Bc][Bc]{\scalebox{.8}{$X$}} \psfrag{Y}[Bc][Bc]{\scalebox{.8}{$Y$}}
  $\Hli_{X,Y}=\rsdraw{.45}{.9}{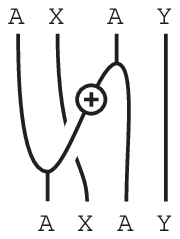}$\quad \text{and} \quad $\Hri_{X,Y}=\rsdraw{.45}{.9}{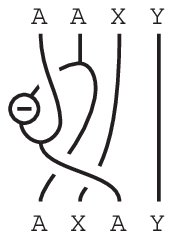}$\,.
\end{center}
Similarly, if $A$ is a Hopf algebra in $\bb$ with invertible antipode, then $? \otimes A$ is a Hopf monad on $\bb$.
Thus Hopf monads generalize Hopf algebras in braided categories. In particular, a Hopf algebra over a commutative ring $\kk$ defines a Hopf monad
on the category of \kt modules.
See Section~\ref{sect-rep} for a detailed discussion of Hopf monads associated with Hopf algebras.
\end{exa}

\begin{rem}\label{rem-copHopf}
Let $T$ be a bimonad on a monoidal category $\cc$. Then $T$ is a right (pre\ti)Hopf monad if and only if its opposite bimonad $T^\cop$ on $\cc^{\otimes \opp}$ (see Remark~\ref{remoppositemonad}) is a left  (pre\ti)Hopf monad.
\end{rem}

\subsection{Hopf monads and \frob\ adjunctions}
In view of the relation between bimonads and comonoidal adjunctions recalled in Section~\ref{sect-bimonadadj}, it is natural to look for a characterization of Hopf monads in terms of adjunctions. This leads to the notion of a Hopf adjunction.


Let $(F\co \cc \to \dd, U\co \dd \to \cc)$ be a comonoidal adjunction between monoidal categories (see Section~\ref{sect-bimonadadj}). The \emph{left \frob\ operator} and the \emph{right \frob\ operator} of $(F,U)$ are the natural transformations $$\FO^l \co F(1_\cc \otimes U) \to F \otimes 1_\dd \quad \text{and} \quad \FO^r \co F(U \otimes 1_\cc) \to 1_\dd \otimes F$$ defined by:
\begin{align*}
&\FO^l_{c,d}=(Fc\otimes \varepsilon_d) F_2(c,Ud)\co F(c \otimes Ud) \to Fc \otimes d, \\
& \FO^r_{d,c} =(\varepsilon_d \otimes  Fc ) F_2(Ud,c)\co F(Ud \otimes c) \to d \otimes Fc,
\end{align*}
for $c \in \Ob(\cc)$ and $d \in \Ob(\dd)$.

 \begin{rem} \frob\ adjunctions were initially introduced by Lawvere in the context of cartesian categories under the name of Frobenius adjunctions~\cite{Law}.
\end{rem}

\begin{rem} \label{rem-fus-from-frob}
Let $T=UF$ be the bimonad of the comonoidal adjunction $(F,U)$. The fusion operators $H^l$ and $H^r$ of $T$ are related to the \frob\ operators
$\FO^l$ and $\FO^r$ of $(F,U)$
as follows:
$$
H^l_{X,Y}=U_2(FX,FY)U(\FO^l_{X,FY}) \quad \text{and} \quad
H^r_{X,Y}=U_2(FX,FY)U(\FO^r_{FX,Y})
$$
for all  $X,Y\in \Ob(\cc)$.
\end{rem}

A \emph{left (resp.\@ right) \frob\ adjunction} is a comonoidal adjunction $(F,U)$ such that $\FO^l$ (resp.\@ $\FO^r$) is invertible.
A \emph{\frob\ adjunction} is a  comonoidal adjunction such that both $\FO^l$ and $\FO^r$ are invertible.

A \emph{left (resp.\@ right) \ufrob\ adjunction} is a comonoidal adjunction $(F,U)$ such that $\FO^l_{\un,-}$ (resp.\@ $\FO^r_{-,\un}$) is invertible.
A \emph{\ufrob\ adjunction} is a  comonoidal adjunction such that both $\FO^l_{\un,-}$ and $\FO^r_{-,\un}$ are invertible.

From Remark~\ref{rem-fus-from-frob}, we easily deduce:
\begin{prop}\label{thm-adj-frob-hopf}
\begin{enumerate}
\labela
\item The monad of a left (resp.\@ right) \frob\ adjunction is a left (resp.\@ right) Hopf monad. In particular the monad of a \frob\ adjunction is a Hopf monad.
\item  The monad of a left (resp.\@ right) \ufrob\ adjunction is a left (resp.\@ right) \uHopf\ monad. In particular the monad of \ufrob\ adjunction is a \uHopf\ monad.
\end{enumerate}
\end{prop}

On the other hand,  a bimonad is a Hopf monad if and only if its associated comonoidal adjunction is a Hopf adjunction:

\begin{thm}\label{thm-hopf-frob}
Let $T$ be a bimonad on a monoidal category $\cc$.
\begin{enumerate}
\labela
\item $T$ is a left (resp.\@ right) Hopf monad if and
only if the comonoidal adjunction $(F_T,U_T)$ is a left (resp.\@ right) \frob\ adjunction. In particular $T$ is a Hopf monad if and
only if $(F_T,U_T)$ is a \frob\ adjunction.
\item  $T$ is a left (resp.\@ right) \uHopf\ monad if and
only if the comonoidal adjunction $(F_T,U_T)$ is a left (resp.\@ right) \ufrob\ adjunction. In particular $T$ is a \uHopf\ monad if and
only if $(F_T,U_T)$ is a \ufrob\ adjunction.
\end{enumerate}
\end{thm}

We prove Theorem~\ref{thm-hopf-frob} in Section~\ref{sect-proofhopfadj}.

Hopf adjunctions are stable under composition:

\begin{prop}\label{prop-comp-frob}
The composite of two left (resp.\@ right) \frob\ adjunctions is a left (resp.\@ right) Hopf adjunction. In particular the composite of two \frob\ adjunctions is a Hopf adjunction.
\end{prop}
Proposition~\ref{prop-comp-frob} is a direct consequence of the following lemma:

\begin{lem}\label{lem-frob-comp}
Let $(F\co \cc \to \dd,U\co \dd \to \cc)$ and $(G\co \dd \to \ee, V\co \ee \to \dd)$ be two comonoidal adjunctions.
Denote by $\FO^l$ (resp. $\FO'^l$, resp. $\FO''^l$) and $\FO^r$ (resp. $\FO'^r$, resp. $\FO''^r$) the left and right \frob\ operators of $(F,U)$ (resp. $(G,V)$, resp. $(GF,UV)$). Then
$$\FO''^l_{c,e}= \FO'^l_{Fc,e}\, G(\FO^l_{c,Ve})\quad\text{and}\quad \FO''^r_{e,c}= \FO'^r_{e, Fc}\, G(\FO^r_{Ve, c})$$
for all $c \in \Ob(\cc)$ and $e \in \Ob(\ee)$.
\end{lem}

\subsection{Proof of Theorem~\ref{thm-hopf-frob}} \label{sect-proofhopfadj}
The `if' part of each assertion results immediately from Proposition~\ref{thm-adj-frob-hopf}, since $T$ is the bimonad of its comonoidal adjunction. The `only if'  part, less straightforward, results from  the following lemma:



\begin{lem}\label{prop-phi-inv}
Let $T$ be a bimonad on a monoidal category $\cc$. Denote by $H^l$ and $H^r$ its fusion operators and $\FO^l$, $\FO^r$ the \frob\ operators of the
adjunction $(F_T,U_T)$ of~$T$. Let $X$ be an object $\cc$. Then $H^l_{X,-}$ is invertible if and only if $\FO^l_{X,-}$ is invertible, and in that case their inverses are related by:
$$\Hli_{X,Y}=\FOli_{X,F_T Y}\quad \text{and} \quad \FOli_{X,(M,r)}=T(\id_X \otimes r)\Hli_{X,M} (\id_{TX}\otimes\eta_M ).$$
Similarly $H^r_{-,X}$ is invertible if and only if $\FO^r_{-,X}$ is invertible, and in that case:
$$\Hri_{Y,X}=\FOri_{F_T Y,X}\quad \text{and} \quad \FOri_{(M,r),X}=T(r \otimes \id_X)\Hri_{M,X} (\eta_M \otimes\id_{TX}).$$
\end{lem}

\begin{proof}  By Remark~\ref{rem-fus-from-frob}, the forgetful functor $U_T$ being strict monoidal, we have $H^l_{X,Y}=\FO^l_{X,F_T(Y)}$  and $H^r_{X,Y}=\FO^r_{F_T(X),Y}$.
Hence the `if' parts and the expressions given for inverses  of fusion operators.

Let us prove the `only if' part of the left-handed case (the right-handed case can be done similarly).
Assume $H^l_{X,-}$ is invertible. Set $A= T(X \otimes ?)$,
$B=TX \otimes ?$, and $\alpha=U_T\FO^l_{X, -}\co AU_T \to BU_T$. We have $\alpha_{F_T}=H^l_{X,-}$ and so $\alpha_{F_T}$ is invertible. Therefore $\alpha$ is invertible by Lemma~\ref{lem-inv-free} below. Thus $\FO^l$ is invertible ($U_T$ being conservative) and
$\FOli_{X,(M,r)}=T(\id_X \otimes r)\Hli_{X,M} (\id_{TX}\otimes\eta_M )$
for any $T$\ti module $(M,r)$.
\end{proof}

\begin{lem}\label{lem-inv-free}
Let $\alpha\co  A U_T \to B U_T$ be a natural transformation, where $T$ is a monad on a category $\cc$ and $A,B \co \cc \to \dd$ are two functors. If $\alpha_{F_T}$ is invertible, so is $\alpha$,
and $$\alpha^{-1}_{(M,r)}= A(r)\alpha^{-1}_{F_TM}B(\eta_M)$$ for any $T$\ti module.
\end{lem}

\begin{proof}
Let $(M,r)$ be a $T$\ti module. The fork $$\xymatrix{T^2M \ar@/^.4pc/[r]^{\mu_M} \ar@/^-.4pc/[r]_{Tr} &TM \ar[r]^r& M}$$ in $\cc$ is split by
$\xymatrix{T^2M &TM \ar[l]_-{\eta_{TM}}& M\ar[l]_-{\eta_M}}$.
As a result, in the diagram:
$$\xymatrix{
A T^2M \ar[d]_{\alpha_{F_TTM}} \dar{A\mu_M}{ATr}
&ATM \ar[d]^{\alpha_{F_TM}}\ar[r]^{Ar}& \ar[d]^{\alpha_{(M,r)}}AM\\
BT^2M \dar{B\mu_M}{BTr} &BTM \ar[r]^{Br}& BM
}$$
the two rows are split coequalizers and the first two columns are invertible by assumption. Therefore the third column is also invertible. Since $r\co F_TM \to (M,r)$ is $T$\ti linear, we obtain:
$\alpha_{(M,r)}^{-1} =\alpha_{(M,r)}^{-1}B(r \eta_M)=A(r) \alpha_{F_T M}^{-1}B(\eta_M)$.
\end{proof}

\section{Hopf monads on closed monoidal categories}\label{sect-HM-closed}

In this section we define binary left and right antipodes for a bimonad $T$ on a closed monoidal category $\cc$ and show that $T$ is a Hopf monad if and only if $T$ admits binary left and right antipodes, or equivalently, if the category of $T$\ti modules is closed monoidal and the forgetful functor $U_T$ preserves internal Homs. When $\cc$ is autonomous, Hopf monads as defined in the present paper coincide with Hopf monads defined in~\cite{BV2} in terms of unary antipodes.

The general results on Hopf monads on closed monoidal categories are stated in Section~\ref{sect-HMclosed-results} and the autonomous case is studied in Section~\ref{sect-auto}. The rest of the section is devoted to the proofs which are based on a classification of adjunction liftings (see Section~\ref{sect-lift-adj}).

\subsection{Closed monoidal categories}\label{subsect-closed}
 See \cite{EK} for a general reference.  Let $\cc$ be a monoidal category.
Let $X,Y$ be two objects of $\cc$. A \emph{left internal Hom from $X$ to $Y$} is an object $\lhom{X,Y}$ endowed with a morphism $\leva{X}{Y} \co \lhom{X,Y} \otimes X \to Y$
such that,
for each object $Z$ of $\cc$, the mapping
$$\left\{\begin{array}{ccc}\Hom_\cc(Z,\lhom{X,Y})&\to &\Hom_\cc(Z \otimes X,Y)\\
f &\mapsto &\leva{X}{Y}(f \otimes \id_X)
\end{array}\right.$$
is a bijection.
If a left internal Hom from $X$ to $Y$ exists, it is unique up to unique isomorphism.

A monoidal category $\cc$ is \emph{left closed} if left internal Homs exist in $\cc$.
This is equivalent to saying that, for every object
$X$ of $\cc$, the endofunctor $? \otimes X$ admits a right adjoint $\lhom{X,?}$, with adjunction unit and counit:
$$\leva{X}{Y} \co  \lhom{X,Y}\otimes  X\to Y \quad\text{and}\quad\lcoeva{X}{Y} \co Y \to \lhom{X, Y\otimes X},$$
called respectively the \emph{left evaluation} and the \emph{left coevaluation}.

Let $\cc$ be a left closed monoidal category. The left internal Homs of~$\cc$ give rise to a functor:
$$
\lhom{-,-}\co \cc^\opp \times \cc \to \cc
$$
where $\cc^\opp$ is the category opposite to $\cc$. Moreover, from the associativity and unitarity of the monoidal product of $\cc$, we deduce isomorphisms
$$\lhom{X \otimes Y, Z} \simeq \lhom{X, \lhom{Y,Z}} \quad \text{and} \quad \lhom{\un,X} \simeq X$$
which we will abstain from writing down in formulae. The composition $$c_{X,Y,Z} \co\lhom{Y,Z} \otimes \lhom{X,Y}  \to \lhom{X,Z}$$ of internal Homs
is the natural transformation defined by:
$$\leva{X}{Z}(c_{X,Y,Z} \otimes X)=\leva{Y}{Z}(\lhom{Y,Z} \otimes \leva{X}{Y}).$$

\begin{rem}
If $X$ is an object of a monoidal category $\cc$ admitting a left dual $(\ldual{X},\lev_X,\lcoev_X)$ then, for every object $Y$ of $\cc$,
$\lhom{X,Y}=Y\otimes \ldual{X}$ is a left internal Hom from $X$ to $Y$, with evaluation morphism  $\leva{X}{Y}=Y \otimes \lev_X$.
Therefore any left autonomous category is left closed monoidal.
\end{rem}
\begin{rem}\label{lem-carotte}
A left closed  monoidal category $\cc$  is left autonomous if and only if
$c_{X,\un,X}\co \lhom{\un,X} \otimes  \lhom{X,\un}\to \lhom{X,X}$ is an isomorphism for all object $X$ of $\cc$. In that case, $\ldual{X}=\lhom{X,\un}$ is a left dual of $X$, with evaluation $\lev_X=\leva{X}{\un}$ and coevaluation $\lcoev_X=(\leva{X}{\un}\otimes \id_{\ldual{X}} )c^{-1}_{X,\un,X}\ \lcoeva{X}{\un}$.
\end{rem}

One defines similarly \emph{right internal Homs} and \emph{right closed} monoidal categories. A monoidal category is right closed if and only if, for every object
$X$ of $\cc$, the endofunctor $X \otimes ?$ has a right adjoint $\rhom{X,?}$, with adjunction unit and counit:
$$\reva{X}{Y} \co  X \otimes \rhom{X,Y}\to Y \quad\text{and}\quad\rcoeva{X}{Y} \co Y \to \lhom{X, X\otimes Y},$$
called respectively the \emph{right evaluation} and the \emph{right coevaluation}.
The right internal Homs of a monoidal left right closed category~$\cc$ give rise to a functor:
$$
 \rhom{-,-}\co \cc^\opp \times \cc \to \cc .
$$

\begin{rem}\label{rem-leftorightclosed}
A right internal Hom in a monoidal category $\cc$ is a left internal Hom in $\cc^{\otimes \opp}$, and $\cc$ is right closed
 if and only if $\cc^{\otimes \opp}$ is left closed.
\end{rem}

A \emph{closed monoidal category} is a monoidal category which is both left and right closed.

\subsection{Functors preserving internal Homs}\label{subsect-closed-fun}
Let $X,Y$ be objects of a monoidal category $\dd$ which have a left internal Hom $\lhom{X,Y}$, with evaluation morphism
$\leva{X}{Y}\co \lhom{X,Y}\otimes X \to Y$.
A monoidal functor $U\co \dd \to \cc$  is said to \emph{preserve the left internal Hom from $X$ to $Y$} if $U\lhom{X,Y}$, endowed with the evaluation
$$U(\leva{X}{Y})U_2(\lhom{X,Y},X)\co U\lhom{X,Y} \otimes UX \to UY\,,$$
is a left internal Hom from $UX$ to $UY$.

A monoidal functor $U\co \dd \to \cc$ between left closed monoidal categories is \emph{left closed} if it preserves all left internal Homs.

Let $U\co \dd \to \cc$  be a monoidal functor between left closed monoidal categories.
The natural transformation
$U(\leva{X}{Y}) U_2(\lhom{X,Y},X)\co U\lhom{X,Y} \otimes UX \to UY$
induces by universal property of internal Homs a natural transformation:
$$U^l_{X,Y}\co U\lhom{X,Y} \to \lhom{UX,UY}.$$
The monoidal functor $U$ is left closed if and only if $U^l$ is an isomorphism.

Similarly one defines monoidal functors \emph{preserving right internal Homs} and \emph{right closed} monoidal functors.

\begin{lem}\label{lem-closed-aut} Let  $U \co \dd \to \cc$ be a strong monoidal functor between left closed monoidal categories.
If $U$ is conservative, left  closed, and $\cc$ is left autonomous, then $\dd$ is left autonomous.
\end{lem}

\begin{proof}
According to Remark~\ref{lem-carotte}, it is enough to show that,
for any object $X$ of $\dd$, the composition morphism $c_{X,\un,X}\co \lhom{\un,X} \otimes \lhom{X,\un} \to \lhom{X,X}$ is an isomorphism. Since $U$ is strong monoidal, $U_2$ and $U_0$ are isomorphisms.
Consider the following commutative diagram:
$$\xymatrix{
U(\lhom{\un,X} \otimes \lhom{X,\un}) \ar[r]^-{U(c_{X,\un,X})}& U\lhom{X,X} \ar[r]^-{U^l_{X,X}} & \lhom{UX,UX}\\
U\lhom{\un,X} \otimes U\lhom{X,\un}  \ar[u]_-{U_2(\lhom{\un,X},\lhom{X,\un})}  \ar[rd]|-{U^l_{\un,X} \otimes U^l_{X,\un}} &  & \lhom{\un,UX} \otimes \lhom{UX,\un} \ar[u]^-{c_{UX,\un,UX}}\\
&\lhom{U\un,UX} \otimes \lhom{UX,U\un}\ar[ru]|-{\lhom{U_0,UX} \otimes \lhom{UX,U_0^{-1}}} &
}$$
Since $U^l$ is an isomorphism ($U$ being left closed) and $c_{UX,\un,UX}$ is invertible (by Remark~\ref{lem-carotte}), we obtain that $U(c_{X,\un,X})$ is invertible. Now $U$ is conservative. Hence $c_{X,\un,X}$ is an isomorphism.
%
%
\end{proof}

\begin{prop}\label{thm-closed-adj-hopf}
Let  $(F\co \cc \to \dd,U\co \dd \to \cc)$ be a comonoidal adjunction
between monoidal left (resp.\@ right) closed  categories.
Then $(F,U)$ is a left (resp.\@ right) \frob\ adjunction if and only if $U$ is left (resp.\@ right) closed.
\end{prop}

\begin{proof}
We prove the left-handed version (the right one can be done similarly).
Let $(F \co \cc \to \dd, U \co \dd \to \cc)$ be a comonoidal adjunction between left closed monoidal categories.
For any $c \in \Ob(\cc)$ and $d,e \in \Ob(\dd)$, set
$$
h_{c,d}^e\co \left\{ \begin{array}{ccc} \Hom_\dd(F(c) \otimes d, e) &\to &\Hom_\dd(F(c \otimes Ud), e) \\
\alpha & \mapsto & \alpha \FO^l_{c,d} \end{array} \right.  $$
where $\FO^l$ is the left Hopf operator of $(F,U)$. Note that $h_{c,d}^e$ is natural in $c,d,e$ and
one verifies easily that it is
the composition:
\begin{align*}
\Hom&_\dd(Fc \otimes d,e) \iso \Hom_\dd(Fc,\lhom{d,e})\iso \Hom_\cc(c,U\lhom{d,e}) \\& \xymatrix{\ar[r]^-{u_{c,d,e}} & \Hom_\cc(c,\lhom{Ud,Ue}) }
\iso \Hom_\cc(c \otimes Ud,Ue) \iso \Hom_\dd(F(c\otimes Ud),e),
\end{align*}
where $u_c^{d,e}=\Hom_\cc(c,U^l_{d,e})$ and all other maps are adjunction bijections.

Assume that $U$ is left closed. Let $c \in \Ob(\cc)$ and $d \in \Ob(\dd)$. Since $U^l_{d,-}$ is an isomorphism, $u_c^{d,-}$ is an isomorphism, and so is $h_{c,d}^-$. Therefore $\FO^l_{c,d}$ is invertible by the Yoneda lemma. Hence $(F,U)$ is a left Hopf adjunction.

Conversely, suppose that $(F,U)$ is a left Hopf adjunction.  Let $d,e \in \Ob(\dd)$. Since $\FO^l_{-,d}$ is an isomorphism, $h_{-,d}^e$ is an isomorphism, and so is $u_-^{d,e}$. Therefore $U^l_{d,e}$ is invertible by the Yoneda lemma. Hence $U$ is left closed.
\end{proof}

\subsection{Hopf monads and antipodes in the closed monoidal setting}\label{sect-HMclosed-results}

Let $T$ be a bimonad on a monoidal category $\cc$.

If $\cc$ is left closed,  a \emph{binary left antipode} for $T$, or simply \emph{left antipode} for $T$, is a natural transformation
$$s^l=\{s^l_{X,Y} \co T\lhom{TX,Y} \to \lhom{X,TY}\}_{X,Y \in \Ob(\cc)}$$
satisfying the following two axioms:
\begin{subequations}
\begin{align}& T\bigl(\leva{X}{Y}(\lhom{\eta_X,Y}\otimes X)\bigr)=
\leva{TX}{TY}(s^l_{TX,Y} T\lhom{\mu_X,Y} \otimes TX)T_2(\lhom{TX,Y},X),\label{anti1}
\\
&\lhom{X,TY \otimes \eta_X}\lcoeva{X}{TY}=\lhom{X,(TY \otimes \mu_X)T_2(Y,TX)}s^l_{X,Y\otimes TX} T(\lcoeva{TX}{Y})
 ,\label{anti2}
\end{align}
\end{subequations}
for all objects $X,Y$ of $\cc$.

Similarly if $\cc$ is right closed, a \emph{binary right antipode} for $T$, or simply \emph{right antipode} for $T$, is a natural transformation
$$s^r=\{s^r_{X,Y} \co T\rhom{TX,Y} \to \rhom{X,TY}\}_{X,Y \in \Ob(\cc)}$$
satisfying:

\begin{subequations}
\begin{align}& T\bigl(\reva{X}{Y}(X\otimes \rhom{\eta_X,Y})\bigr)=
\reva{TX}{TY}(TX \otimes s^r_{TX,Y} T\rhom{\mu_X,Y} )T_2(X,\rhom{TX,Y}),
\\
& \rhom{X, \eta_X\otimes TY}\rcoeva{X}{TY}=\rhom{X,(\mu_X\otimes TY)T_2(TX,Y)}s^r_{X,TX \otimes Y} T(\rcoeva{TX}{Y})
 ,
\end{align}
\end{subequations}
for all objects $X,Y$ of $\cc$.

With this definition of (binary) antipodes, we have:

\begin{thm}\label{thm-hopfmon-closed} Let $T$ be a bimonad on a left (resp.\@ right) closed monoidal category~$\cc$.
The following assertions are equivalent:
\begin{enumerate}
\renewcommand{\labelenumi}{{\rm (\roman{enumi})}}
\item  The bimonad $T$ is a left (resp.\@ right) Hopf monad on $\cc$;
\item The monoidal category  $\mo{T}{\cc}$ is left (resp.\@ right) closed and the forgetful functor $U_T$ is left (resp.\@ right) closed;
\item The bimonad $T$ admits a left (resp.\@ right) binary antipode.
\end{enumerate}
\end{thm}
This theorem is proved in Section~\ref{sect-HMclosed-proof}.

\begin{rem}
If the equivalent conditions of Theorem~\ref{thm-hopfmon-closed} are satisfied, internal Homs in $\mo{T}{\cc}$ may be constructed in terms of the antipodes of $T$ as follows.
If $T$ is a left Hopf monad and $\cc$ is left closed monoidal, then a left internal Hom for any two $T$\ti modules $(M,r)$ and $(N,t)$ is given by:
\begin{align*}
& \lhom{(M,r),(N,t)}=\bigl(\lhom{M,N}, \lhom{M,t}s^l_{M,N} T\lhom{r,N}\bigr), \\
& \leva{(M,r)}{(N,t)}=\leva{M}{N}, \quad \text{and} \quad \lcoeva{(M,r)}{(N,t)}=\lcoeva{M}{N}.
\end{align*}
Similarly, if $T$ is a right Hopf monad and $\cc$ is right closed monoidal, then a right internal Hom for any two $T$\ti modules $(M,r)$ and $(N,t)$ is given by:
\begin{align*}
& \rhom{(M,r),(N,t)}=\bigl(\rhom{M,N}, \rhom{M,t}s^r_{M,N} T\rhom{r,N}\bigr), \\
& \reva{(M,r)}{(N,t)}=\reva{M}{N}, \quad \text{and} \quad \rcoeva{(M,r)}{(N,t)}=\rcoeva{M}{N}.
\end{align*}
\end{rem}

%

In addition to characterizing Hopf monads on closed monoidal categories, the left and right antipodes, when they exist, are unique and well-behaved with respect to the bimonad structure:



\begin{prop}\label{prop-anti-multi} Let $T$ be a bimonad on a monoidal category $\cc$.
\begin{enumerate}
\labela
\item If $\cc$ is left (resp.\@ right) closed and $T$ admits a left (resp.\@ right) antipode, then this antipode is unique.
\item Assume $\cc$ is left closed and $T$ is a left Hopf monad. Then the left antipode $s^l$ for $T$ satisfies:
\begin{align*}
& s^l_{X,Y}\mu_{\lhom{TX,Y}}=\lhom{X,\mu_Y} s^l_{X,TY}T(s^l_{TX,Y})T^2\lhom{\mu_X,Y},\\
& s^l_{X,Y} \eta_{\lhom{TX,Y}}=\lhom{\eta_X,\eta_Y}, \\
& s^l_{X \otimes Y,Z}T\lhom{T_2(X,Y),Z}=\lhom{X,s^l_{Y,Z}}s^l_{X,\lhom{T(Y),Z}}, \\
& s^l_{T(\un),X}\lhom{T_0,X}=\id_{TX},
\end{align*}
for all objects $X,Y,Z$ of $\cc$.
\item Assume $\cc$ is right closed and $T$ is a right Hopf monad. Then the  right antipode
$s^r$ for $T$ satisfies:
\begin{align*}
& s^r_{X,Y}\mu_{\rhom{TX,Y}}=\rhom{X,\mu_Y} s^r_{X,TY}T(s^r_{TX,Y})T^2\rhom{\mu_X,Y},\\
& s^r_{X,Y} \eta_{\rhom{TX,Y}}=\rhom{\eta_X,\eta_Y},\\
& s^r_{X \otimes Y,Z}T\rhom{T_2(X,Y),Z}=\rhom{X,s^r_{Y,Z}}s^r_{X,\rhom{T(Y),Z}}, \\
& s^r_{T(\un),X}\rhom{T_0,X}=\id_{TX},
\end{align*}
for all objects $X,Y,Z$ of $\cc$.
\end{enumerate}
\end{prop}
The proposition  is proved in Section~\ref{sect-HMclosed-proof}.

Lastly, the antipodes and the inverses of the fusion operators of a Hopf monad
can be expressed in terms of one another, as follows:

\begin{prop}\label{prop-Hs}
If $T$ is a left Hopf monad on a left closed monoidal category $\cc$,
then the inverse of the left fusion operator $H^l$ and the left antipode $s^l$ are related as follows:
\begin{align*}
\Hli_{X,Y}&=T(X \otimes \mu_Y)\leva{TY}{T(X \otimes T^2Y)}(s^l_{TY,X\otimes T^2Y}T(\lcoeva{T^2Y}{X})\otimes \id_{TY})\,,\\
s^l_{X,Y}&=\lhom{X,T\leva{TX}{Y}} \lhom{\eta_X,{\Hli_{X,\lhom{TX,Y}}}}\lcoeva{TX}{T\lhom{TX,Y}}\,.
\end{align*}
Similarly if $T$ is a right Hopf monad on a right closed monoidal category $\cc$,
then the inverse of the right fusion operator $H^r$ and the right antipode $s^l$ are related as follows:
\begin{align*}
\Hri_{X,Y}&=T(\mu_X \otimes Y)\reva{TX}{T(T^2 X \otimes Y)}(\id_{TX} \otimes s^r_{TX,T^2X\otimes Y}T(\rcoeva{T^2X}{Y})\,,\\
s^r_{X,Y}&=\rhom{X,T\reva{TX}{Y}} \rhom{\eta_X,{\Hri_{X,\rhom{TX,Y}}}}\rcoeva{TX}{T\rhom{TX,Y}}\,.
\end{align*}
\end{prop}
The proposition is proved in Section~\ref{sect-HMclosed-proof}.

\subsection{Hopf monads on autonomous categories}\label{sect-auto}

\newcommand{\s}{{\mathfrak s}}

The notion of Hopf monad introduced in this paper is a generalization of the notion of Hopf monad
on an autonomous category introduced in~\cite{BV2}.

If $T$ is a bimonad on a left autonomous category $\cc$, a \emph{unary left antipode} for $T$, or simply \emph{left antipode} for $T$, is a natural transformation
$$\s^l=\{\s^l_X\co {T(\ldual{TX})}\to \ldual{X}\}_{X\in \Ob(\cc)}$$
satisfying:
\begin{align*}
& T_0 T(\lev_X)T(\ldual{\eta_X} \otimes X)=\lev_{TX}\bigl(\s^l_{TX}T(\ldual{\mu}_X) \otimes
TX\bigr)T_2(\ldual{TX},X); \\
& (\eta_X \otimes \ldual{X})\lcoev_X T_0=(\mu_X \otimes \s^l_X) T_2(TX,\ldual{TX})T(\lcoev_{TX});
\end{align*}
for every object $X$ of $\cc$.

Similarly if $T$ is a bimonad on a right autonomous category $\cc$, a \emph{unary right antipode} for $T$, or simply \emph{left antipode} for $T$, is a natural transformation
$$\s^r=\{\s^r_X\co {T(\rdual{(TX)})}\to \rdual{X}\}_{X\in \Ob(\cc)}$$
satisfying:
\begin{align*}
& T_0 T(\rev_X)T(X \otimes \eta_X^\vee)=\rev_{TX}\bigl(TX
\otimes \s^r_{TX}T(\mu_X^\vee)\bigr)T_2(X,\rdual{(TX)});\\
& (X^\vee\otimes \eta_X )\rcoev_X T_0=(\s^r_X \otimes\mu_X) T_2(\rdual{(TX)},TX)T(\rcoev_{TX});
\end{align*}

In \cite{BV2}, a left (resp.\@ right) Hopf monad $T$ on a left (resp.\@ right) autonomous category $\cc$ is defined as a bimonad on $\cc$ which admits a left (resp.\@ right) unary antipode or, equivalently by \cite[Theorem 3.8]{BV2}, whose category of modules $\mo{T}{\cc}$ is left (resp.\@ right) autonomous. This definition, which makes sense only in the autonomous setting, agrees with that given in Section~\ref{sect-hopfmon}:

\begin{thm}\label{thm-anti-auto}
Let $\cc$ be a left (resp.\@ right) autonomous category and $T$ be a bimonad on $\cc$.
Then the following assertions are equivalent:
\begin{enumerate}
\labeli
\item The bimonad $T$ has a left (resp.\@ right) unary antipode;
\item The bimonad $T$ has a left (resp.\@ right) binary antipode;
\item The bimonad $T$ is a left (resp.\@ right) Hopf monad.
\end{enumerate}
\end{thm}
The theorem is proved in Section~\ref{sect-HMclosed-proof}

\begin{rem}
The binary left antipode $s^l_{X,Y}$
and unary left antipode $\s^l_X$ of a left Hopf monad $T$ on a left autonomous category are related as follows:
$$s^l_{X,Y}=(TX\otimes \s^l_Y)T_2(X,\ldual{TY}) \blabla{and} \s^l_X=(T_0 \otimes \ldual{X})s^l_{X,\un}\,.$$
Similarly the binary right antipode $s^r_{X,Y}$
and unary right antipode $\s^l_X$ of a right Hopf monad $T$ on a right autonomous category are related as follows:
$$s^r_{X,Y}=(\s^r_Y\otimes TX)T_2(\ldual{TY},X) \blabla{and} \s^r_X=(\ldual{X}\otimes T_0) s^r_{X,\un}\,.$$
\end{rem}

\subsection{Lifting adjunctions} \label{sect-lift-adj}
In this section, $(T,\mu,\eta)$ is a monad on a category $\cc$ and $(T',\mu',\eta')$ is a monad on
a category $\cc'$.

A \emph{lift of a functor $G\co\cc \to \cc'$ along $(T,T')$} is a functor $\widetilde{G} \co \mo{T}{\cc} \to \mo{T'}{\cc'}$ such that
$U_{T'}\widetilde{G}=G U_T$.
It is a well-known fact that such lifts $\widetilde{G}$ are in bijective correspondence with
natural transformations $\zeta \co T'G \to G T$ satisfying:
\begin{equation*}
\zeta \mu'_G=G(\mu)\zeta_{T} T'(\zeta) \quad \text{and} \quad
\zeta \eta'_G=G(\eta).
\end{equation*}
Such a natural transformation $\zeta$ is called a \emph{lifting datum for $G$ along $(T,T')$}.

The lift $\widetilde{G}^\zeta$ corresponding with a  lifting datum $\zeta$ is defined by
$$\widetilde{G}^\zeta(M,r)=(G(M), G(r)\zeta_M).$$
Conversely, the lifting datum associated with a lift $\widetilde{G}$
is
$$\zeta=U_{T'}(\varepsilon'_{\widetilde{G} F_T}) T'G(\eta),$$
where $\varepsilon'$ denotes the counit of the adjunction
$(F_{T'}, U_{T'})$.

Consider two functors $G, G'\co \cc \to \cc'$, a lifting datum $\zeta$ for $G$, and a lifting datum $\zeta'$ for $G'$. Then a natural transformation $\alpha\co G \to G'$ lifts to a natural transformation
$$\widetilde{\alpha}\co \widetilde{G}^\zeta \to \widetilde{G'}^{\zeta'}$$
(in the sense that  $U_{T'} (\widetilde{\alpha})= \alpha_{U_T}$)  if and only if it satisfies $\zeta' T'(\alpha)=\alpha_T \zeta$.

\begin{exa} \label{rem-frob-lift}
Let $T$ be a bimonad on a monoidal category $\cc$ and
 $(M,r)$ be a $T$\ti module. Then the endofunctors  $? \otimes M$ and $M \otimes ?$ of $\cc$ lift to endofunctors  $?\otimes (M,r)$ and $(M,r) \otimes ?$ of $\mo{T}{\cc}$.
The lifting data corresponding with these lifts are
the \frob\ operators
$\FO^l_{-,(M,r)}$ and $\FO^r_{(M,r),-}$ of the comonoidal adjunction $(F_T,U_T)$.
\end{exa}

Now let $(G\co \cc \to \cc',V\co \cc' \to \cc)$ be an adjunction, with unit $h\co 1_\cc \to VG$ and counit $e\co GV \to 1_{\cc'}$.

A \emph{lift of the adjunction $(G,V)$ along $(T,T'$)} is an adjunction $(\widetilde{G},\widetilde{V})$, where $\widetilde{G}\co \mo{T}{\cc} \to \mo{T'}{\cc'}$ is a lift of $G$ along $(T,T')$, $\widetilde{V}\co \mo{T'}{\cc'} \to \mo{T}{\cc}$ is a lift of $V$ along $(T',T)$, and the unit $\widetilde{h}$ and counit $\widetilde{e}$ of $(\widetilde{G},\widetilde{V})$ are lifts of
$h$ and $e$ respectively.

Lifts of the adjunction $(G,V)$ are in bijective correspondence with pairs $(\zeta,\xi)$,
where $\zeta\co T'G \to GT$ and $\xi\co TV \to VT'$
are natural transformations satisfying the following axioms:
\begin{subequations}\label{liftadj}
\begin{align}
\zeta \mu'_G&=G(\mu)\zeta_{T} T'(\zeta)\,,\label{liftzeta1}\\
\zeta \eta'_G&=G(\eta)\,,\label{liftzeta2}\\
\xi \mu_V&=V(\mu')\xi_{T'} T(\xi)\,,\label{liftxi1}\\
\xi \eta_V&=V(\eta')\,,\label{liftxi2}\\
T'(e)&=e_{T'} G(\xi) \zeta_V\,,\label{ehlinear1}\\
h_T&=V(\zeta) \xi_G T(h)\label{ehlinear2}\,.
\end{align}
\end{subequations}
Such a pair $(\zeta,\xi)$ is called a \emph{lifting datum for the adjunction $(G,V)$ along $(T,T')$}.


%
By adjunction, we have a bijection
$$\Phi\co\left\{ \begin{array}{ccl}
\Nat(TV,VT')&\to &\Nat(GT,T'G)\\
\xi &\mapsto &\Phi(\xi)= e_{T'G} G(\xi_G) GT(h)
\end{array}\right.$$
whose inverse is
given by $\Phi^{-1}(\alpha)=VT'(e) V(\alpha_V) h_{TV}$.

\begin{thm} \label{thm-lift-adj}
Let $\zeta\co T'G \to GT'$ be a lifting datum for $G$ along $(T,T')$.
Then the following assertions are equivalent:
\begin{enumerate}
\labeli
\item There exists a natural transformation $\xi\co TV \to VT'$ such that $(\zeta,\xi)$ is a lifting datum for the adjunction $(G,V)$ along $(T,T')$.
\item $\zeta$ is invertible.
\end{enumerate}
If such is the case,  $\xi$ is unique and $\xi=\Phi^{-1}(\zeta^{-1})$.
\end{thm}

The theorem, which may be interpreted in terms of doctrinal adjunctions, results immediately from the following lemma:

\begin{lem}\label{lem1}
Let $\zeta\co T'G \to GT$ and $\xi\co TV \to VT'$
be natural transformations.
\begin{enumerate}
\labela
\item Axiom~\eqref{ehlinear1} is equivalent to $\Phi(\xi)\zeta=\id_{T'G}$, and \eqref{ehlinear2} to $\zeta\Phi(\xi)=\id_{GT}$.
\item If \eqref{ehlinear1} and~\eqref{ehlinear2} hold, then \eqref{liftzeta1} is equivalent to \eqref{liftxi1},
and \eqref{liftzeta2}  to \eqref{liftxi2}.
\end{enumerate}
\end{lem}

\begin{proof}
The adjunction bijection $\Nat(T'GV,T') \shortiso \Nat(T'G,T'G)$, defined by $\beta \mapsto \beta_G T'G(h)$,
sends $T'(e)$ to
$\id_{T'G}$, and $e_{T'} G(\xi) \zeta_V$ to $e_{T'G} G(\xi_G) \zeta_{VG}T'G(h)=\Phi(\xi) \zeta$.
Similarly the adjunction bijection $\Nat(T,VGT) \shortiso \Nat(GT,GT)$ sends $h_T$ to
$\id_{GT}$ and $V(\zeta) \xi_G T(h)$ to $\zeta\Phi(\xi)$. Hence Part (a).

Now assume that Axioms~\eqref{ehlinear1} and~\eqref{ehlinear2} hold. In other words, $\zeta$ is invertible
and $\zeta^{-1}=\Phi(\xi)$.
Then Axiom~\eqref{liftzeta1} and Axiom~\eqref{liftzeta2} can be re-written  as $\Phi(\xi) G(\mu)=\mu'_GT'(\Phi(\xi))\Phi(\xi)_T $ and $\Phi(\xi) G(\eta)=\eta'_G$,
which translate respectively to Axiom~\eqref{liftxi1} and  Axiom~\eqref{liftxi2} via the adjunction bijections $\Nat(GT^2,T'G) \shortiso \Nat(T^2V,VT')$ and $\Nat(G,T'G) \shortiso \Nat(V,VT')$. Hence Part (b).
\end{proof}

\subsection{Bimonads and lifting adjunctions} \label{sect-HMclosed-proof}
Here, by applying the results of Section~\ref{sect-lift-adj} to bimonads in closed monoidal categories, we prove Theorems~\ref{thm-hopfmon-closed} and \ref{thm-anti-auto} and Propositions~\ref{prop-anti-multi} and \ref{prop-Hs}. We deal with the left closed case, from which the right closed case results  using the coopposite bimonad (see Remarks~\ref{rem-copHopf} and~\ref{rem-leftorightclosed}).

Let $\cc$ be a monoidal category and $T$ be a bimonad on $\cc$.
Note that $T \times 1_{\mo{T}{\cc}}$ is a bimonad on $\cc \times \mo{T}{\cc}$.
The monoidal tensor product $\otimes \co \mo{T}{\cc} \times \mo{T}{\cc} \to \mo{T}{\cc}$ of $\mo{T}{\cc}$ is a lift of the functor $\otimes(1_\cc \times U_T) \co \cc \times \mo{T}{\cc} \to \cc$ along $(T \times 1_{\mo{T}{\cc}},T)$:
$$\xymatrix{
\mo{T}{\cc} \times \mo{T}{\cc} \ar[d]_{U_T\times 1_{\mo{T}{\cc}}} \ar[rr]^-{\otimes}&& \mo{T}{\cc} \ar[d]^{U_T}\\
\cc \times \mo{T}{\cc} \ar[rr]_-{\otimes(1_\cc \times U_T)} && \cc
}$$
The corresponding lifting datum $\zeta \co T( 1_\cc \otimes U_T) \to T \otimes U_T$ is given by:
$$\zeta^{(M,r)}_X=(X \otimes r)T_2(X,M)\co T(X \otimes M)
\to T(X) \otimes M.$$ Note that $\zeta=U_T(\FO^l)$, where $\FO^l$ denotes  the left \frob\ operator  of the comonoidal
adjunction $(F_T,U_T)$, and so, $U_T$ being conservative, $\zeta$ is invertible if and only if $T$ is a left Hopf monad.

Assume now that $\cc$ is left closed, that is,  we have an adjunction
$(? \otimes X, \lhom{X,?})$ for each $X \in \Ob(\cc)$. In particular $(? \otimes M)_{(M,r) \in \mo{T}{\cc}}$ is a family of endofunctors of~$\cc$ admitting right adjoints indexed by $\mo{T}{\cc}$.


\begin{lem}\label{lem-lift-bimod}
The following assertions are equivalent:
\begin{enumerate}
\labeli
\item The category $\mo{T}{\cc}$ is left closed monoidal and $U_T$ is left closed;
\item For each $T$\ti module $(M,r)$, the adjunction $(? \otimes M,\lhom{M,?})$ lifts to an
adjunction $(? \otimes (M,r), \widetilde{V}_{(M,r)})$.
\end{enumerate}
\end{lem}

\begin{proof}
Let us prove that (i) implies (ii). Recall that since $U_T$ is left closed, 
we have a natural isomorphism $U_T^l \co U_T \lhom{\,\,\, ,\,\, } \to \lhom{U_T,U_T}$, see Section~\ref{subsect-closed-fun}.
Thus, by transport of structure, we may choose left internal Homs in $\mo{T}{\cc}$ so that $U_T \lhom{(M,r),(N,r)}$  is equal to $\lhom{M,N}$, $U^l_T$ being the identity. Then the adjunction $(? \otimes (M,r), \lhom{(M,r),?})$
is a lift of the adjunction $(? \otimes M, \lhom{M,?})$.

Conversely (ii) implies (i) since the existence of an adjunction $(? \otimes (M,r), \widetilde{V}_{(M,r)})$ lifting  the adjunction $(? \otimes M, \lhom{M,?})$ means firstly that $\mo{T}{\cc}$ is left closed monoidal, with $\lhom{(M,r),?}=\widetilde{V}_{(M,r)}$, and secondly that  $U^l_T$
is the identity (and so $U_T$ is left closed).
\end{proof}

Let us prove Theorem~\ref{thm-hopfmon-closed}, Propositions~\ref{prop-anti-multi} and \ref{prop-Hs}, and Theorem~\ref{thm-anti-auto}.

\begin{proof}[Proof of Theorem~\ref{thm-hopfmon-closed}]
According to Theorem~\ref{thm-lift-adj}, given a $T$\ti module $(M,r)$, the adjunction $(? \otimes M,\lhom{M,?})$ lifts to an
adjunction $(? \otimes (M,r), \widetilde{V}^{(M,r)})$ if and only if the lifting datum
$\zeta^{(M,r)}$ is invertible. Therefore, by Lemma~\ref{lem-lift-bimod}, $\mo{T}{\cc}$ is left closed monoidal and $U_T$ is left closed if and only if $\zeta$ is invertible, and so if and only if $T$ is a left Hopf monad. Hence the equivalence of
assertions (i) and (ii).

Assume (i) holds, so that $\zeta$ is invertible. By Theorem~\ref{thm-lift-adj}, for any $T$\ti module $(M,r)$, there exists a unique natural transformation $$\xi^{(M,r)}\co T\lhom{M,?} \to \lhom{M,T}$$ such that
$(\zeta^{(M,r)},\xi^{(M,r)})$ is a lifting datum for the adjunction $(? \otimes M,\lhom{M,?})$ along $(T,T)$, which is given by
$$\xi^{(M,r)}_X=\lhom{M,T(\leva{M}{X})} \lhom{M,\zeta^{(M,r)\stackrel{-1}{\phantom{.}}}_{\lhom{M,X}}} \lcoeva{M}{T\lhom{M,X}}.$$
Note that $\xi$ is natural in $(M,r)$. Axioms~\eqref{ehlinear1} and \eqref{ehlinear2} for this lifting datum are:
\begin{subequations}
\begin{align}T(\leva{M}{X})&=\leva{M}{TX} (\xi^{(M,r)}_X \otimes M) \zeta^{(M,r)}_{\lhom{M,X}},\label{newehlinear1}\\
\lcoeva{M}{TX}&=\lhom{M,\zeta^{(M,r)}_X} \xi^{(M,r)}_{X \otimes M} T(\lcoeva{M}{X})\label{newehlinear2}.
\end{align}
\end{subequations}
They translate to Axioms~\eqref{anti1} and~\eqref{anti2} of a left antipode under the adjunction bijection:
$$\Psi\co\left \{\begin{array}{ccl}
\Nat(T\lhom{U_T,1_\cc},\lhom{U_T,T})&\to &\Nat(T\lhom{T,1_\cc},\lhom{1_\cc,T}) \\
\xi &\mapsto & s^l=\{s^l_{X,Y}=\lhom{\eta_X,TY}\xi^{F_TX}_Y\}_{X,Y\in\Ob(\cc)}
\end{array} \right.$$
Hence assertion (iii).

Conversely assume (iii) holds. Denote by $s^l$ the left antipode of $T$. Set $\xi=\Psi^{-1}(s^l)$, that is, $\xi^{(M,r)}_Y=s^l_{M,Y} T\lhom{r,Y}$. Under $\Psi^{-1}$, Axioms~\eqref{anti1} and~\eqref{anti2} for $s^l$ translate to \eqref{newehlinear1} and~\eqref{newehlinear2}. In particular, for any $T$\ti module $(M,r)$, $\xi^{(M,r)}$ satisfies~\eqref{ehlinear1} and~\eqref{ehlinear2}. Furthermore, Axioms~\eqref{liftzeta1} and~\eqref{liftzeta2} hold for $\zeta^{(M,r)}$
as it is a lifting datum for $? \otimes M$. Thus, by Lemma~\ref{lem1}, Axioms~\eqref{liftxi1} and~\eqref{liftxi2} hold for $\xi^{(M,r)}$, that is:
\begin{subequations}
\begin{align}
\xi^{(M,r)}_X \mu_{\lhom{M,X}}&=\lhom{M,\mu_X}\xi^{(M,r)}_{TX} T(\xi^{(M,r)}_X),\label{newliftxi1}\\
\xi^{(M,r)}_X \eta_{\lhom{M,X}}&=\lhom{M,\eta_X}.\label{newliftxi2}
\end{align}
\end{subequations}
Therefore $(\zeta^{(M,r)},\xi^{(M,r)})$ is a lifting datum for the adjunction $(? \otimes M,\lhom{M,?})$ along $(T,T)$. Hence (ii) by Lemma~\ref{lem-lift-bimod}. This concludes the proof of Theorem~\ref{thm-hopfmon-closed}.
\end{proof}

\begin{proof}[Proof of Proposition~\ref{prop-anti-multi}]
Part (a) results from the fact that if a natural transformation $\xi$ satisfying \eqref{newehlinear1} and~\eqref{newehlinear2} exists,
it is unique by Theorem~\ref{thm-lift-adj}.

Let us prove Part (b). Assume that $T$ admits a left antipode $s^l$. When translated in terms of $s^l$, Axioms~\eqref{newehlinear1} and~\eqref{newehlinear2} yield the compatibility of $s^l$ with $\mu$ and $\eta$.
%
Given two $T$\ti modules $(M,r)$ and $(N,t)$, the $T$\ti action of the left internal Hom $\lhom{(M,r),(N,t)}$ obtained by lifting $\lhom{M,N}$
is $\lhom{M,t}s^l_{M,N} T\lhom{r,N}$.
Given a third $T$\ti module $(P,p)$, the $T$\ti linearity of the canonical isomorphism
$$\lhom{(M,r) \otimes (N,t),(P,p)} \simeq \lhom{(M,r),\lhom{(N,t),(P,p)}},$$
translated in terms of $s^l$, yields the compatibility of $s^l$ to $T_2$. Similarly the $T$\ti linearity of the canonical isomorphism
$$\lhom{(\un,T_0),(M,r)} \simeq (M,r)$$ yields the compatibility of $s^l$ to $T_0$. Hence Proposition~\ref{prop-anti-multi}.
\end{proof}

\begin{proof}[Proof of Proposition~\ref{prop-Hs}] Denote by $s^l$ the left antipode of $T$ and set $\xi=\Psi^{-1}(s^l)$. Recall that $\xi^{(M,r)}_Y=s^l_{M,Y} T\lhom{r,Y}$ and $s^l_{X,Y}=\lhom{\eta_X,TY}\xi^{F_TX}_Y$.
By Theorem~\ref{thm-lift-adj},
\begin{align*}&\xi^{(M,r)}_X =\lhom{M,T(\leva{M}{X})} \lhom{M,\zeta^{(M,r)\stackrel{-1}{\phantom{.}}}_{\lhom{M,X}}} \lcoeva{M}{T\lhom{M,X}}, \\
&\zeta^{(M,r)\stackrel{-1}{\phantom{.}}}_X= \lev^M_{X \otimes M}(\xi^{(M,r)}_{X\otimes M} \otimes M)T(\lcoev^M_X),
\end{align*}
where  $\zeta^{(M,r)}_X=\FO^l_{X,(M,r)}$. By Lemma~\ref{prop-phi-inv}, we have: $\Hli_{X,Y}={\zeta^{F_TY}_X}^{-1}$ and
$$\zeta^{(M,r)\stackrel{-1}{\phantom{.}}}_X=T(\id_X \otimes r)\Hli_{X,M} (\id_{TX}\otimes\eta_M ).$$
Hence the expression of $s^l$ in terms of $\Hli$, and conversely.
\end{proof}

\begin{proof}[Proof of Theorem~\ref{thm-anti-auto}]  We prove the left handed version. Assertions (ii) and (iii) are equivalent by Theorem~\ref{thm-hopfmon-closed}. Assertion (iii) is equivalent to $\mo{T}{\cc}$ and $U_T$ being left closed, and so to $\mo{T}{\cc}$ being left autonomous (using Lemma~\ref{lem-closed-aut} and the fact that a strong monoidal functor preserves left duals). Hence (ii) is equivalent to (i) by \cite[Theorem 3.8]{BV2}.
\end{proof}

\section{Cross products and related constructions}\label{sect-crosstruc}

In this section we study the cross product of Hopf monads (previously introduced in \cite{BV3} for Hopf monads on autonomous categories). In particular we introduce the inverse operation, called the cross quotient.

\subsection{Functoriality of categories of modules}\label{sect-functodogmod}

Let $\cc$ be a category. If $T$ is a monad on $\cc$, then $(\mo{T}{\cc},U_T)$ is a category over $\cc$, that is, an object of $\Cat/\cc$.
Any morphism $f \co T \to P$ of monads  on $\cc$ induces a functor
$$f^* \co\left\{\begin{array}{ccl}
 \mo{P}{\cc} &\to &\mo{T}{\cc}\\
(M,r) &\mapsto &(M,rf_M)
\end{array}\right. $$
 over $\cc$, that is, $U_T f^* = U_{P}$.  Moreover, any functor $F \co \mo{P}{\cc} \to \mo{T}{\cc}$ over $\cc$ is of this form.
This construction defines a fully faithful functor
$$\left\{\begin{array}{ccc}\Mon(\cc)^\opp &\to & \Cat/\cc\\
T& \mapsto& (\mo{T}{\cc}, U_T)
\end{array}\right. $$

If $f \co T \to P$ is a morphism of bimonads on a monoidal category $\cc$, then $f^* \co \mo{P}{\cc} \to \mo{T}{\cc}$ is a strict monoidal functor over $\cc$, and any strong monoidal functor $F \co \mo{P}{\cc} \to \mo{T}{\cc}$ over $\cc$ (that is, such that $U_T F=U_{P}$ as monoidal functors) is of this form (see \cite[Lemma 2.9]{BV2}). Hence a fully faithful functor $$\BiMon(\cc)^\opp \to \MonCat/\cc.$$

\subsection{Exactness properties of monads}
A Hopf monad $T$ on a monoidal category~$\cc$ admits a right adjoint if $\cc$ is autonomous (see \cite[Corollary 3.12]{BV2}), but not in general.
In many cases, the existence of a right adjoint can be replaced by the weaker condition of preservation of reflexive coequalizers (defined in Section~\ref{sect-prelims}).

\begin{lem}[\cite{Lin}]\label{lem-pres-creat}  Let $\cc$ be a category and $T$ be a monad on $\cc$ preserving reflexive coequalizers. Then:
\begin{enumerate}
\labela
\item  A reflexive pair of morphisms of $\mo{T}{\cc}$ whose image by $U_T$ has a coequalizer, has a coequalizer, and this coequalizer is preserved by $U_T$;
\item If reflexive coequalizers exist in $\cc$, they exist also in $\mo{T}{\cc}$ and $U_T$ preserves them.
\end{enumerate}
\end{lem}





\begin{lem}\label{lem-tens-pres-coeq} Let $\cc$ be a monoidal category admitting  reflexive coequalizers, which are preserved by monoidal product on the left (resp.\@ right). If $T$ is a bimonad  on $\cc$ preserving reflexive coequalizers, then $\mo{T}{\cc}$ has reflexive coequalizers which are preserved  by monoidal product on the left (resp.\@ right).
\end{lem}
%

\begin{proof} Let us prove the right-handed version.
According to Lemma~\ref{lem-pres-creat}, $\mo{T}{\cc}$ has reflexive coequalizers and $U_T$ preserves them.
Let $h$ be a coequalizer of a reflexive pair $(f,g)$ in $\mo{T}{\cc}$, and $d$ be an object of $\mo{T}{\cc}$. Denoting $k$ be a coequalizer of the reflexive pair $(f \otimes d,g \otimes d)$, the morphism  $h \otimes d$ factorizes uniquely as $\phi k$.
Both  $U_T(h \otimes d)$ and $U_T k$ are coequalizers of $(U_Tf,U_T g)$  (because $U_T$ and $U_T \otimes U_Td$ preserve reflexive coequalizers) so $U_T\phi$ is an isomorphism. Hence $\phi$ is an isomorphism, since $U_T$ is conservative. Thus $h \otimes d$ is a coequalizer of $(f \otimes d, g \otimes d)$.
%
\end{proof}

\subsection{Cross products}

Let $T$ be a monad on a category $\cc$. If $Q$ is a monad on the category $\mo{T}{\cc}$ of $T$\ti modules,
the monad of the composite adjunction  $$\mo{Q}{\bigl(\mo{T}{\cc}\bigr)}\adjunct{U_Q}{F_Q}\mo{T}{\cc}\adjunct{U_T}{F_T}\cc$$
is called the \emph{cross product of $T$ by $Q$} and denoted by $Q \cp T$ (see \cite[Section 3.7]{BV3}).  As an endofunctor of~$\cc$, $Q \cp T=U_TQF_T$. The product $p$ and unit $e$ of $Q \cp T$ are:
\begin{equation*}
p=q_{F_T} Q(\varepsilon_{Q F_T})  \quad \text{and}\quad e=v_{F_T}\eta,
\end{equation*}
where $q$ and $v$ are the product and the unit of $Q$, and $\eta$ and $\varepsilon$ are the unit and counit of the adjunction $(F_T,U_T)$.

Note that the composition of two monadic functors is not monadic in general. However:

\begin{prop}[\cite{BaWe}]\label{prop-cp-monadic}
If $T$ is a monad on a category $\cc$ and $Q$ is a monad on $\mo{T}{\cc}$ which preserves reflexive coequalizers, then the forgetful functor $U_T U_Q$ is monadic with monad $Q \cp T$. Moreover the comparison functor $$K \co (\mo{T}{\cc})^Q \to \mo{Q\cp T}{\cc}$$ is an isomorphism of categories.
\end{prop}


If $T$ is a bimonad on a monoidal category $\cc$ and $Q$ is a bimonad on $\mo{T}{\cc}$, then  $Q\cp T=U_TQF_T$ is a bimonad on $\cc$ (since a composition of comonoidal adjunctions is a comonoidal adjunction), with comonoidal structure given by:
\begin{align*}
& (Q\cp  T)_2(X,Y)=Q_2\bigl(F_T(X),F_T(Y)\bigr)\, Q\bigl((F_T)_2(X,Y)\bigr),\\
& (Q\cp  T)_0=Q_0\, Q\bigl((F_T)_0\bigl).
\end{align*}
In that case the comparison functor $K \co (\mo{T}{\cc})^Q \to \mo{Q\cp T}{\cc}$
is strict monoidal.

The cross product is functorial in $Q$: the assignment $Q \mapsto Q \cp T$ defines a functor $? \cp T \co \BiMon(\mo{T}{\cc}) \to \BiMon(\cc)$.

From Proposition~\ref{prop-comp-frob} and Proposition~\ref{thm-adj-frob-hopf}, we deduce:

\begin{prop}\label{prop-cp-hopf}
The cross product of two left (resp.\@ right) Hopf monads is a left (resp.\@ right) Hopf monad.
 In particular, the cross product of two Hopf monads is a Hopf monad.
\end{prop}

\begin{exa}
Let $H$ be a bialgebra over a field $\kk$ and $A$ be a $H$-module algebra, that is, an algebra in the monoidal category $\lMod{H}$ of left $H$\ti modules. In this situation, we may form the cross product $A \rtimes H$, which is a $\kk$-algebra (see \cite{Maj2}). Recall that $H \otimes ?$ is a monad on $\Vect_\kk$ and $A \otimes ?$ is a monad on $\lMod{H}$. Then:
\begin{equation*}
(A \otimes ?)\cp (H \otimes ?)=(A\rtimes H) \otimes ?
\end{equation*}
as monads. Moreover, if $H$ is a quasitriangular Hopf algebra and $A$ is a $H$-module Hopf algebra, that is, a Hopf algebra in the braided category $\lMod{H}$, then $A\rtimes H$ is a Hopf algebra over $\kk$, and
$(A \otimes ?)\cp (H \otimes ?)=(A\rtimes H) \otimes ?$ as Hopf monads.
\end{exa}

\begin{exa}\label{exa-double}
Let $T$ be a Hopf monad on an autonomous category $\cc$. Assume $T$ is \emph{centralizable}, that is, for all object $X$ of $\cc$, the coend $$Z_T(X)=\int^{Y \in \cc} \ldual{T(Y)}\otimes X \otimes Y$$
exists (see \cite{BV3}). In that case, the assignment  $X \mapsto Z_T(X)$ is a Hopf monad on $\cc$, called the centralizer of $T$ and denoted by $Z_T$.
The centralizer $Z_T$ of $T$ lifts canonically to a Hopf monad $\widetilde{Z_T}$ on $\mo{T}{\cc}$, which is the centralizer of $1_{\mo{T}{\cc}}$. Then, by \cite[Theorem 6.5]{BV3},  $D_T=\widetilde{Z_T} \cp T$ is a quasitriangular Hopf monad, called the \emph{double} of $T$, and satisfies $\zz(\mo{T}{\cc})\cong \mo{D_T}{\cc}$ as braided categories, where $\zz$ denotes the categorical center.
\end{exa}

 Distributive laws, introduced by Beck \cite{Beck1} to encode composition and lifting of monads, have been adapted to Hopf monads on autonomous categories in \cite{BV3}. The next corollary deals with the case of Hopf monads on arbitrary monoidal categories.

\begin{cor}\label{cor-crossprod-distlaw}
Let $P$ and $T$ be Hopf monads on a monoidal category $\cc$ and $\Omega \co TP \to PT$ be a comonoidal distributive law of $T$ over $P$.
\begin{enumerate}
\labela
\item If $P$ is a Hopf monad, then the lift $\Tilde{P}^\Omega$ of $P$ is a Hopf monad on~$\mo{T}{\cc}$.
\item If $P$ and $T$ are Hopf monads, then the composition $P\circ_\Omega T$ is a Hopf monad on $\cc$ and $\Tilde{P}^\Omega \cp  T=P\circ_\Omega T$ as Hopf monads.
\end{enumerate}
\end{cor}
\begin{proof}
Recall from \cite[Theorem 4.7]{BV3} that $\Omega$ defines a bimonad $\Tilde{P}^\Omega$ on $\mo{T}{\cc}$, which is a lift of the bimonad $P$, and a bimonad $P\circ_\Omega T$ on $\cc$, with underlying functor $PT$.
Moreover
$P\circ_\Omega T=\Tilde{P}^\Omega \cp  T$ as bimonads on $\cc$. The forgetful functor $U_T$ maps the fusion operators of $\Tilde{P}^\Omega$ to those of $P$. Therefore if $P$ is a Hopf monad, so is $\Tilde{P}^\Omega$ (as $U_T$ is conservative). If both $P$ and $T$ are Hopf monads, then $\Tilde{P}^\Omega \cp  T$ is a Hopf monad by Proposition~\ref{prop-cp-hopf}, and so is $P\circ_\Omega T$.
\end{proof}

\begin{lem}\label{lem-Q-Hopf}
Let $T$ be a bimonad on a monoidal category $\cc$ and $Q$ be a bimonad on $\mo{T}{\cc}$.
Assume that the monoidal products of $\mo{T}{\cc}$ and $(\mo{T}{\cc})^Q$ preserve reflexive coequalizers in the left (resp.\@ right) variable.
If the adjunction $(F_Q F_T, U_T U_Q)$ is a left (resp.\@ right) Hopf adjunction and $T$ is a left (resp.\@ right) Hopf monad, then
 $Q$ is a left (resp.\@ right) Hopf monad.
\end{lem}

\begin{proof} Let us prove the left handed version. Denote by $\FO^l$,  $\FO'^l$, and
$\FO''^l$ the left \frob\ operators of the adjunctions $(F_T,U_T)$,  $(F_Q,U_Q)$, and $(F_QF_T,U_TU_Q)$ respectively. Assume that $(F_QF_T,U_TU_Q)$ is a left Hopf adjunction, that is $\FO''^l$ is invertible. Assume also that $T$ is a left Hopf monad. By Theorem~\ref{thm-hopf-frob},  $\FO^l$ is invertible and
it is enough to show that $\FO'^l$ is also invertible. Let $e$ be an arbitrary object of $(\mo{T}{\cc})^Q$. The natural transformation
$
\FO'^l_{-,e} \co F_Q(? \otimes U_Q e) \to F_Q \otimes e
$
is invertible on the image of $F_T$, since $ \FO'^l_{F_T,e}=\FO''^l_{-,e}\,F_Q(\FO^l_{-,U_Q(e)})^{-1}$ by Lemma~\ref{lem-frob-comp}.
Now let $(M,r)$ be a $T$\ti module. The coequalizer
$$\xymatrix{F_TT(M)  \ar@<2pt>[r]^-{\mu_M}\ar@<-2pt>[r]_-{F_T(r)}  &F_T(M) \ar[r]^{r}& (M,r)}$$
is reflexive because $F_T(r) F_T(\eta_M)=\id_{F_T(M)}=\mu_M F_T(\eta_M)$. This reflexive coequalizer is preserved by the
functors $ F_Q(? \otimes U_Q (e))$ and $F_Q \otimes e$, because $F_Q$ is a left adjoint and $? \otimes U_Q(e)$ and $? \otimes e$ preserve reflexive coequalizers (by hypothesis). Hence $\FO'^l_{(M,r),e}$ is invertible.
\end{proof}

%
%
%
%
%
%

%
%
%

\subsection{Cross quotients}\label{sect-croco}
Let $f \co T \to P$ be a morphism of monads on a category $\cc$. We say that $f$ is \emph{cross quotientable} if the functor $f^* \co \mo{P}{\cc} \to \mo{T}{\cc}$ is monadic. In that case, the monad of $f^*$ (on $\mo{T}{\cc}$) is called the \emph{cross quotient} of $f$ and is denoted by $P \cdiv_f T$ or simply $P \cdiv T$. Note that the comparison functor
$$ \xymatrix@R=1.4em@C=.8em{\cc^P \ar[rr]^-K \ar[dr]_-{f^*}& & (\cc^T)^{P \cdiv T} \ar[dl]^{U_{P \cdiv T}} \\
&\cc^T&\\
}$$
is then an isomorphism of categories (by the last assertion of Theorem~\ref{thm-beck}).

%


\begin{lem}[\cite{Lin}]\label{lem-exist-quotient}
Let $f \co T \to P$ be a morphism of monads on a category $\cc$. The following assertions are equivalent:
\begin{enumerate}
\labeli
\item The morphism $f$ is cross quotientable;
\item The functor $f^*$ admits a left adjoint;
\item For any $T$\ti module $(M,r)$, the reflexive pair
$$\xymatrix@C=4.5em{F_P T M \ar@<2pt>[r]^-{P(r)}\ar@<-2pt>[r]_-{p_M P(f_M)} &F_{P} M}$$
admits a coequalizer $F_{P}M \to G(M,r)$ in $\mo{P}{\cc}$, where $p$ is the product of $P$.
\end{enumerate}
If these conditions hold, a left adjoint $f_!$ of $f^*$ is given by $f_!(M,r)=G(M,r)$.
\end{lem}

\begin{proof}
It results from Beck's theorem (see Theorem~\ref{thm-beck}) that if $U$ and $V$ are composable functors such that both $UV$ and $U$ are monadic, then $V$ is monadic if and only if it admits a left adjoint. Thus (i) is equivalent to (ii).

Now let $(M,r)$ be a $T$\ti module and $(N,\rho)$ be a $P$\ti module. The pair of Assertion~(iii) is reflexive (since $F_P(\eta_M)$ is a common retraction). Via the
adjunction bijection
$$\Hom_{\mo{P}{\cc}}(F_PM,(N,\rho))\simeq \Hom_{\cc}(M,U_P(N,\rho))=\Hom_{\cc}(M,N),$$
morphisms $F_P M\to (N,\rho)$ which coequalize that pair  correspond with $T$\ti linear morphisms $(M,r) \to f^*(N,\rho)$.
Therefore (ii) is equivalent to~(iii). We conclude using the last assertion of Theorem~\ref{thm-beck}.
\end{proof}

\begin{rem}\label{rem-exist-P/T}
A morphism $f \co T \to P$ of monads on $\cc$  is cross quotientable whenever $\cc$ admits coequalizers of reflexive pairs and $P$ preserve them.
\end{rem}

A cross quotient of bimonads is a bimonad: let $f \co T \to P$ be a cross quotientable morphism of bimonads on a monoidal category~$\cc$. Since $f^*$ is strong monoidal,  $P \cdiv_f T$ is a bimonad on $\cc^T$ and the comparison functor
$K \co \mo{P}{\cc} \to (\mo{T}{\cc})^{P \cdiv_f T}$ is an isomorphism of monoidal categories.

The cross quotient is inverse to the cross product in the following sense:

\begin{prop}\label{prop-prod-quot} Let $T$ be a (bi)monad on a (monoidal) category $\cc$.
\begin{enumerate} \labela
\item
If $T \to P$ is a cross quotientable morphism of (bi)monads on $\cc$, then $$(P \cdiv T) \cp T\simeq P$$ as (bi)monads.
\item Let $Q$ be a (bi)monad on $\mo{T}{\cc}$ such that $U_TU_Q$ is monadic. Then the unit  of $Q$ defines a cross quotientable morphism of (bi)monads $T \to Q\cp T$ and $$(Q\cp T)\cdiv T \simeq Q$$ as (bi)monads.
\end{enumerate}
\end{prop}

\begin{proof}
Let us prove Part (a). Since $\mo{P}{\cc} \simeq (\mo{T}{\cc})^{P \cdiv T}$, the functor $U_{P \cdiv T} U_T$ is monadic. Hence an isomorphism $\cc^P \simeq \cc^{(P \cdiv T) \cp T}$ of (monoidal) categories over~$\cc$, which induces an isomorphism $(P \cdiv T) \cp T \simeq P$ of (bi)monads on $\cc$.

Let us prove Part (b). Set $f=u \cp T \co T \to Q \cp T$, where $u$ is the unit of $Q$.  We have a commutative diagram of (monoidal) functors:
$$\xymatrix{
(\mo{T}{\cc})^Q \ar[r]^K \ar[d]_{U_Q} &\mo{Q\cp T}{\cc}\ar[ld]_{f^*}\ar[d]^{U_{Q\cp T}}\\
\mo{T}{\cc} \ar[r]_{U_T}&\cc
}$$
where $K$ is the comparison functor of the adjunction $(F_QF_T,U_T U_Q)$. Since $K$ is a equivalence, the functor $f^*$ is monadic, with (bi)monad $(Q\cp T)\cdiv T$. Hence $K$ induces an isomorphism of (bi)monads $(Q\cp T)\cdiv T \simeq Q$.
\end{proof}

\begin{rem} \label{rem-funct-crossquot}
Let $T$ be a bimonad on a monoidal category $\cc$. Let $\BiMon(\mo{T}{\cc})_m$ be the full subcategory of $\BiMon(\mo{T}{\cc})$ whose objects are monads $Q$ on $\mo{T}{\cc}$ such that $U_TU_Q$ is monadic. Let $T\backslash \BiMon(\cc)_q$ be the full subcategory of $T\backslash \BiMon(\cc)$ whose objects are quotientable morphisms of bimonads from $T$.
Then the functor
$$\prettydef{
\BiMon(\mo{T}{\cc}) &\to &T\backslash \BiMon(\cc)\\
Q &\mapsto& (Q, T \to Q \cp T)}$$
induces an equivalence of categories
$\BiMon(\mo{T}{\cc})_m \simeq T\backslash \BiMon(\cc)_q$,
with quasi-inverse given by $(T \to P) \mapsto (P\cdiv T)$.
\end{rem}

Under suitable exactness assumptions, if $P$ and $T$ are Hopf monads, so is $P \cdiv T$:

\begin{prop}\label{prop-P/T-Hopf} Let $\cc$ be a monoidal category admitting reflexive coequalizers, and whose monoidal product preserves reflexive coequalizers in the left (resp.\@ right) variable.
Let $T$ and $P$ be two left (resp.\@ right) Hopf monads on $\cc$ which preserve reflexive coequalizers.
Then any morphism of bimonads $T \to P$  is cross quotientable and  $P \cdiv T$ is a left (resp.\@ right) Hopf monad.
\end{prop}


\begin{proof} Let us prove the left-handed version. The morphism $f$ is cross quotientable by Remark~\ref{rem-exist-P/T}, and so $P\simeq (P\cdiv_f T) \cp T$ as bimonads.
The monoidal products of $\mo{T}{\cc}$ and $\mo{P}{\cc}$ preserve reflexive coequalizers in the left variable by Lemma~\ref{lem-tens-pres-coeq}.
Applying Lemma~\ref{lem-Q-Hopf} to the bimonads $T$ and $P\cdiv T$, we get that $P\cdiv T$ is a left Hopf monad.
\end{proof}

\begin{exa}
Let $f \co L \to H$ be a morphism of Hopf algebras over a field $\kk$, so that $H$ becomes a $L$\ti bimodule by setting $\ell \cdot h \cdot \ell'=f(\ell)hf(\ell')$.  The morphism $f$ induces a morphism of Hopf monads on
$\Vect_\kk$:
 $$f \otimes_\kk ? \co L \otimes_\kk ? \to H \otimes_\kk ?$$
  which is cross quotientable,  and $(H \otimes ?)\cdiv (L \otimes ?)$  is a \kt linear Hopf monad on the monoidal category $\lMod{L}$ given by $N\mapsto H \otimes_L N$. This construction defines an equivalence of categories
$$L \backslash \HopfAlg_\kk \to \HopfMon_\kk(\lMod{L}),$$
  where $L \backslash\HopfAlg_\kk$ is the category of Hopf \kt algebras under $L$ and $\HopfMon_\kk(\lMod{L})$ is the category of \kt linear Hopf monads on $\lMod{L}$.
\end{exa}



\section{Hopf monads associated with Hopf algebras and bosonization}\label{sect-rep}

Examples of Hopf monads may be obtained from Hopf algebras. For instance, any Hopf algebra $A$ in a braided category $\bb$ gives rise to Hopf monads $A \otimes ?$ and $? \otimes A$ on $\bb$, see Example~\ref{exa-rep-braided}. More generally, any Hopf algebra $(A,\sigma)$ in the center $\zz(\cc)$ of a monoidal $\cc$ gives rise to a Hopf monad $A\otimes_\sigma ?$ on $\cc$ (see Section~\ref{sect-rep-cent}, or \cite{BV3} for the autonomous case).
Hopf monads of this form are called representable. The main result of this section asserts that a Hopf monad on a monoidal category is representable if and only if it is augmented, that is, endowed with a Hopf monad morphism from itself to the identity (see Theorem~\ref{thm-rep}).

 More generally, given a Hopf monad $T$ on $\cc$ and a Hopf algebra $(\AT,\sigma)$ in the center $\zz(\mo{T}{\cc})$ of the category of $T$\ti modules, we construct a Hopf monad $\AT \cp_\sigma T$ on~$\cc$
of which $T$ is a retract. Conversely,
under suitable exactness conditions (involving  reflexive coequalizers), any Hopf monad $P$ of which $T$ is a retract is of the form $\AT \cp_\sigma T$ for some Hopf algebra $(\AT,\sigma)$ in $\zz(\mo{T}{\cc})$.

\subsection{Lax braidings, lax half-braidings and lax center}

A \emph{lax braiding} of a monoidal category $\cc$ is a natural transformation
$$\tau=\{\tau_{X,Y} \co X \otimes Y \to Y \otimes X\}_{X,Y \in \Ob(\cc)}$$  satisfying:
\begin{align*}
&\tau_{X, Y\otimes Z}=(\id_Y \otimes \tau_{X, Z})(\tau_{X, Y} \otimes \id_Z),\\
&\tau_{X\otimes Y,Z}=(\tau_{X,Z} \otimes \id_Y)(\id_X \otimes \tau_{Y,Z}),\\
&\tau_{X,\un}=\id_X=\tau_{\un,X}.
\end{align*}
A \emph{lax braided category} is a monoidal category endowed with a lax braiding.
A \emph{braiding} is an invertible lax braiding, and a \emph{braided category} is a monoidal category endowed with a braiding.

Let $\cc$ be a monoidal category and $M$ an object of $\cc$. A \emph{lax half-braiding} for $M$ is a natural transformation $\sigma\co M \otimes 1_\cc
\to 1_\cc \otimes M$ such that
$$
\sigma_{Y \otimes Z}=(\id_Y \otimes \sigma_Z)(\sigma_Y \otimes \id_Z) \quad \text{and} \quad
\sigma_\un=\id_M.
$$
The pair $(M,\sigma)$ is then called a \emph{lax half-braiding of $\cc$.}

The \emph{lax center of $\cc$} (see \cite{Sch,DPS}) is the lax braided category $\zzl(\cc)$ defined as follows. Objects of $\zzl(\cc)$ are left half-braidings of $\cc$. A morphism in $\zzl(\cc)$ from $(M,\sigma)$ to $(M',\sigma')$ is a morphism $f \co M \to M'$ in $\cc$ such that:
$(\id_{1_\cc} \otimes f)\sigma=\sigma'(f \otimes \id_{1_\cc})$. The monoidal product and lax braiding $\tau$ are:
\begin{equation*}
(M,\sigma) \otimes (N,\gamma)=\bigl(M \otimes N,(\sigma \otimes \id_N)(\id_M \otimes \gamma) \bigr) \quad \text{and} \quad \tau_{(M,\sigma),(N,\gamma)}=\sigma_{N}.
\end{equation*}

A \emph{half braiding} is  a lax half braiding $(M,\sigma)$ such that $\sigma$ is invertible.
The \emph{center of $\cc$} is the full monoidal subcategory  $\zz(\cc) \subset \zzl(\cc)$ whose objects are half braidings of $\cc$. It is a braided category, with braiding induced by $\tau$.

Note that if $\cc$ is autonomous, lax half braidings are in fact  half braiding, so that the lax center $\zzl(\cc)$ coincides with the center $\zz(\cc)$.

\subsection{Hopf algebras in lax braided categories}

Let $\bb$ be a lax braided category, with lax braiding $\tau$.  A \emph{bialgebra} in $\bb$ is an object $A$ of~$\bb$ endowed with an algebra structure $(m,u)$ and a coalgebra structure $(\Delta,\varepsilon)$ in~$\bb$ satisfying:
\begin{align*}
&\Delta m=(m \otimes m)(\id_A \otimes \tau_{A,A} \otimes \id_A)(\Delta \otimes \Delta),  & & \Delta u=u \otimes u, \\  &\varepsilon m=\varepsilon \otimes \varepsilon, & & \varepsilon u=\id_\un.
\end{align*}
Bialgebras in $\bb$,  together with morphisms of bialgebras (defined in the obvious way), form a category $\BiAlg(\bb)$.

Let $A$ be a bialgebra in $\bb$. An \emph{antipode} of  $A$  is a morphism $S \co A \to A$ in $\bb$ such that:
\begin{equation*}
m(S \otimes \id_A)\Delta=u \varepsilon=m(\id_A \otimes S)\Delta.
\end{equation*}
If it exists, an antipode for $A$ is unique, it satisfies:
$$Sm =m\tau_{A,A}(S \otimes S),\quad Su=u,\quad\Delta S= (S \otimes S)\tau_{A,A} \Delta,\quad \varepsilon S=\varepsilon,$$
and we have: $\tau_{A,A}=(m \otimes m) (S \otimes \Delta m \otimes S) (\Delta \otimes \Delta)$.

If $\tau_{A,A}$ is invertible, an  \emph{opantipode} of $A$ is a morphism $S' \co A \to A$ in $\bb$ such that:
\begin{equation*}
m\tau^{-1}_{A,A}(S' \otimes \id_A)\Delta=u \varepsilon=m\tau^{-1}_{A,A}(\id_A \otimes S')\Delta.
\end{equation*}
If it exists, an opantipode for $A$ is unique.

If $\tau_{A,A}$ is invertible,  the bialgebra $A$  admits an antipode and an opantipode if and only if it admits an invertible antipode, or equivalently, an invertible opantipode. When such is the case, the opantipode is the inverse of the antipode.

Let $A$ be a bialgebra on a lax braided category $\bb$, with lax half-braiding $\tau$. The \emph{fusion operator of $A$} is the morphism  
\begin{center}
$\HAl=(A \otimes m)(\Delta \otimes A) = \psfrag{A}[Bc][Bc]{\scalebox{.8}{$A$}} \rsdraw{.45}{.9}{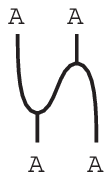}  \co A\otimes A \to A \otimes A.$
\end{center}

The \emph{opfusion operator of $A$} is the morphism
\begin{center}
$\HAr=(m \otimes A) (A \otimes \tau_{A,A})(\Delta \otimes A) = \psfrag{A}[Bc][Bc]{\scalebox{.8}{$A$}} \psfrag{X}[Bc][Bc]{\scalebox{1}{$\tau_{A,A}$}}\rsdraw{.45}{.9}{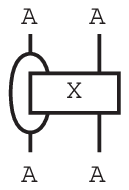}\co A\otimes A \to A \otimes A.$
\end{center}

\begin{lem} \label{lem-anti}
\begin{enumerate}
\labela
\item The bialgebra $A$ admits an antipode $S$ if and only if its fusion operator $\HAl$ is invertible.
If such is the case,  
\begin{align*}
& S=(\varepsilon \otimes A)\HAli(A \otimes u),\\
&\HAli= (A \otimes m)(A \otimes S \otimes A)(\Delta \otimes A).
\end{align*}
\item If $\tau_{A,A}$ is invertible, the bialgebra $A$ admits an opantipode $S'$ if and only if its opfusion operator $\HAr$ is invertible. If such is the
case,
\begin{align*}
&S'=(\varepsilon \otimes A)\HAri(u \otimes A),\\
&\HAri=\tau_{A,A}^{-1}(m \otimes A)(S' \otimes A \otimes A)(\tau^{-1}_{A,A} \otimes A)(A \otimes \Delta).
\end{align*}
\end{enumerate}
\end{lem}

A \emph{Hopf algebra} in a lax braided category $\bb$, with lax braiding $\tau$, is a bialgebra $A$ in $\bb$ admitting an invertible antipode and such that $\tau_{A,A}$ is invertible.
Hopf algebras in $\bb$ form a full subcategory of $\BiAlg(\bb)$ denoted by $\HopfAlg(\bb)$.

\begin{rem}\label{rem-op}
If $\bb$ is a braided category, the mirror $\overline{\bb}$ of $\bb$ is the braided category obtained when the braiding $\tau$ of $\bb$ is replaced by its mirror  $\overline{\tau}$ (defined by $\overline{\tau}_{X,Y}=\tau^{-1}_{Y,X}$). If $(A,m,u,\Delta,\varepsilon)$ is a bialgebra in a braided category $\bb$, one defines a bialgebra
$A^\opp$ in $\overline{\bb}$ by setting $A^\opp=(A,m^\opp,u,\Delta,\varepsilon)$, with $m^\opp=m \tau^{-1}_{A,A}$. We have $(A^\opp)^\opp=A$. An opantipode for $A$ is an antipode for $A^\opp$. 
The bimonads $A \otimes ?$ and $? \otimes A^\opp$ are isomorphic via $\tau_{A,-}$. See \cite[Section 1.11 and Example 2.3]{BV3} for details.
\end{rem}

\subsection{Hopf monads represented by central Hopf algebras}\label{sect-rep-cent}

 Let $\cc$ be a monoidal category. A \emph{(lax) central algebra} (resp.\@ \emph{coalgebra}, resp.\@ \emph{bialgebra}, resp.\@ \emph{Hopf algebra}) of $\cc$
is an algebra (resp.\@ coalgebra, resp.\@ bialgebra, resp.\@ Hopf algebra) in the (lax) center of $\cc$.

Any  lax central coalgebra $(A,\sigma)$ of $\cc$   gives rise to a comonoidal endofunctor of $\cc$, denoted by $A \otimes_\sigma ?$,  defined by $A \otimes ?$ as a functor and endowed with the comonoidal structure:
\begin{center}
$(A \otimes_\sigma ?)_2(X,Y)=(A \otimes \sigma_X )(\Delta \otimes X)\otimes Y=\psfrag{A}[Bc][Bc]{\scalebox{.8}{$A$}} \psfrag{U}[Bc][Bc]{\scalebox{.8}{$X$}} \psfrag{Y}[Bc][Bc]{\scalebox{.8}{$Y$}}\psfrag{X}[Bc][Bc]{\scalebox{1}{$\sigma_X$}} \rsdraw{.45}{.9}{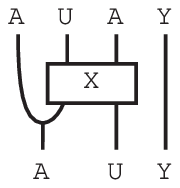}
\,, \quad
(A \otimes_\sigma ?)_0=\varepsilon=\psfrag{A}[Bc][Bc]{\scalebox{.8}{$A$}}\rsdraw{.45}{.9}{epsA},$
\end{center}
where $\Delta$ and $\varepsilon$ denote the coproduct and counit of $(A,\sigma)$.

For any  lax central bialgebra $(A,\sigma)$ of $\cc$,   the comonoidal endofunctor $A \otimes_\sigma ?$ is a bimonad on $\cc$ with  monad structure  given by:
\begin{center}
$\mu_X= m \otimes X =\psfrag{A}[Bc][Bc]{\scalebox{.8}{$A$}} \psfrag{X}[Bc][Bc]{\scalebox{.8}{$X$}}
  \rsdraw{.45}{.9}{mOTA2} \quad\text{and}\quad \eta_X=u\otimes X=\rsdraw{.45}{.9}{uOTA2},$
\end{center}
where
$m$ and $u$ are the product and unit of $A$. Denote by $\lsMod{A}{\sigma}$ the monoidal category $\mo{A \otimes_\sigma ?}{\cc}$, that is, the category of left $A$\ti module (in $\cc$) with monoidal product $(M,r) \otimes (N,s)=(M \otimes N, \omega)$, where
$$
\omega=\psfrag{A}[Bc][Bc]{\scalebox{.8}{$A$}} \psfrag{U}[Bc][Bc]{\scalebox{.8}{$M$}} \psfrag{r}[Bc][Bc]{\scalebox{1}{$r$}}  \psfrag{s}[Bc][Bc]{\scalebox{1}{$s$}} \psfrag{Y}[Bc][Bc]{\scalebox{.8}{$N$}}\psfrag{X}[Bc][Bc]{\scalebox{1}{$\sigma_M$}} \rsdraw{.45}{.9}{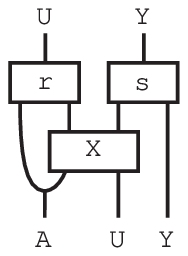} \, ,
$$
and monoidal unit $(\un,\varepsilon)$.

The bimonads of the form $A \otimes_\sigma ?$  can be characterized as follows:
\begin{lem}\label{lem-ATbimon}
Let $A$ be an object of $\cc$ and consider the endofunctor $T=A \otimes ?$ of~$\cc$. Let $\Delta \co A \to A \otimes A$ and $\varepsilon \co A \to \un$ be morphisms in $\cc$ and  $\sigma \co A\otimes ? \to ? \otimes A$ be a natural transformation such that $\sigma_\un=\id_A$. Set
$$ T_2(X,Y)=(A \otimes \sigma_X \otimes Y)(\Delta \otimes X \otimes Y) \quad \text{and} \quad T_0= \varepsilon.$$
Then the following conditions are equivalent:
\begin{enumerate}\labeli
\item $(T,T_2,T_0)$ is a comonoidal endofunctor of $\cc$;
\item $\sigma$ is a lax half braiding for $A$ and $(A,\sigma)$ is a coalgebra in $\zzl(\cc)$ with coproduct $\Delta$ and counit $\varepsilon$.
\end{enumerate}
Assume these equivalent conditions hold. Then $T=A\otimes_\sigma ?$ as comonoidal functors. Furthermore, let $m \co A \otimes A \to A$ and $u \co \un \to A$ be morphisms in $\cc$ and set:
$$\mu = m \otimes ? \co T^2 \to T \quad \text{and} \quad \eta=u \otimes ? \co 1_\cc \to T.$$
Then the following conditions are equivalent:
\begin{enumerate}\labeli
\setcounter{enumi}{2}
\item $T$ is a bimonad with product $\mu$, unit $\eta$, and comonoidal structure $(T_2,T_0)$;
\item $(A,\sigma)$ is a  lax central  bialgebra of $\cc$   with product $m$, unit $u$, coproduct $\Delta$, and counit~$\varepsilon$.
\end{enumerate}
If these equivalent conditions hold, $T=A\otimes_\sigma ?$ as bimonads.
\end{lem}
\begin{proof}
The verification, lengthy but straightforward,  is left to the reader.
\end{proof}

Let $(A,\sigma)$  be a lax central bialgebra of $\cc$, that is, a bialgebra in $\zzl(\cc)$.
The left and right fusion operators of the monad $A \otimes_\sigma ?$ are:
\begin{align*}
&H^l_{X,Y}=(A \otimes X \otimes m)(A \otimes \sigma_X \otimes A)(\Delta \otimes X \otimes A) \otimes Y= \psfrag{A}[Bc][Bc]{\scalebox{.8}{$A$}} \psfrag{U}[Bc][Bc]{\scalebox{.8}{$X$}} \psfrag{Y}[Bc][Bc]{\scalebox{.8}{$Y$}}\psfrag{X}[Bc][Bc]{\scalebox{1}{$\sigma_X$}} \rsdraw{.45}{.9}{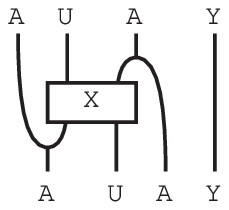}\,,\\
&H^r_{X,Y}=(m \otimes X  \otimes A)(A \otimes \sigma_{A \otimes X})(\Delta \otimes A \otimes X) \otimes Y= \psfrag{A}[Bc][Bc]{\scalebox{.8}{$A$}} \psfrag{U}[Bc][Bc]{\scalebox{.8}{$X$}} \psfrag{Y}[Bc][Bc]{\scalebox{.8}{$Y$}}\psfrag{X}[Bc][Bc]{\scalebox{1}{$\sigma_{A\otimes X}$}} \rsdraw{.45}{.9}{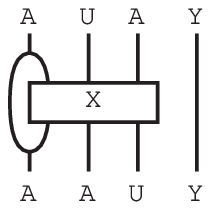}\,.
\end{align*}


\begin{prop}\label{prop-Halg-to-Hmon}
Let $(A,\sigma)$ be a lax central bialgebra in $\cc$, and let $A \otimes_\sigma ?$ be the corresponding bimonad on $\cc$. Then:
\begin{enumerate} \labela
\item The following conditions are equivalent:
    \begin{enumerate}[(i)]
    \item $A \otimes_\sigma ?$ is a left Hopf monad;
    \item $A \otimes_\sigma ?$ is a left \uHopf\ monad;
    \item $A$ admits an antipode;
    \end{enumerate}
\item The following conditions are equivalent:
    \begin{enumerate}[(i')]
    \item $A \otimes_\sigma ?$ is a right Hopf monad;
    \item $A \otimes_\sigma ?$ is a right \uHopf\ monad;
    \item $\sigma$ is invertible and $A$ admits an opantipode.
    \end{enumerate}
\end{enumerate}
In particular, the bimonad $A \otimes_\sigma ?$ is a Hopf monad if and only if $A \otimes_\sigma ?$ is a  \uHopf\ monad, if and only if   $(A,\sigma)$ is a central Hopf algebra of $\cc$, that is, a Hopf algebra in the center $\zz(\cc)$.
\end{prop}

\begin{rem}\label{rem-lmodsigaut}
Let $(A,\sigma)$ be a  central Hopf algebra of $\cc$    and $\lsMod{A}{\sigma}$ be the monoidal category of left $A$\ti modules (with monoidal product induced by $\sigma$). Then the full subcategory $\lsmod{A}{\sigma} \subset \lsMod{A}{\sigma}$ of left $A$\ti modules $(M,r)$ whose underling object $M$ has a left and a right dual is autonomous.
\end{rem}

\begin{rem}\label{rem-rep-braided}
If $\bb$ is a braided category, then its braiding $\tau$ defines a fully faithful braided functor $$\left\{\begin{array}{ccl} \bb & \to &\zz(\bb)\\ X & \mapsto &(X,\tau_{X,-})\end{array} \right.$$
which is a monoidal section of the forgetful functor $\zz(\bb) \to \bb$.
In particular if $A$ is a bialgebra in $\bb$, then $(A,\tau_{A,-})$ is a  central bialgebra of $\bb$   and we have
$$A \otimes ?= A \otimes_{\tau_{A,-}} ?$$ as bimonads on $\bb$, where $A\otimes ?$ is the bimonad constructed in Example~\ref{exa-rep-braided}. Also, if $A$ is a bialgebra in $\bb$, then $A^\opp$ is a bialgebra in the mirror $\overline{\bb}$ of $\bb$ (see  Remark~\ref{rem-op}), $(A^\opp, \overline{\tau}_{A,-})$ is a central bialgebra of $\bb$,   where $\overline{\tau}$ is the mirror braiding of $\tau$, and $$? \otimes A \simeq A^\opp\otimes_{\overline{\tau}_{A,-}} ?$$
as bimonads on $\bb$. Moreover $A$ is a Hopf algebra in $\bb$ if and only if $(A,\tau_{A,-})$ is a  central Hopf algebra of $\bb$,   if and only if $A \otimes ?$ is a Hopf monad on $\bb$, if and only if $? \otimes A$ is a Hopf monad on $\bb$.
\end{rem}

\begin{proof}[Proof of Proposition~\ref{prop-Halg-to-Hmon}]
Let $H^l$ be the left fusion operator of $T=A\otimes_\sigma ?$ and $\HAl$ be the fusion operator of $A$. We have:  $H^l_{X,Y}=H^l_{X,\un} \otimes Y$ and $H^l_{\un,\un}=\HAl$.
Thus the bimonad $T$ is a left Hopf monad if and only if $H^l_{-,\un}$ is an isomorphism, and $T$ is a left
\uHopf\ monad if and only if $\HAl$ is an isomorphism.
Hence (ii) is equivalent to (iii) since, by Lemma~\ref{lem-anti}, $\HAl$ is invertible if and only if $A$ admits an antipode. Assuming (iii) and denoting $S$ the antipode of $A$, one verifies easily that
$(A \otimes X \otimes m)(A \otimes \sigma_X \otimes A)((A \otimes S)\Delta \otimes X \otimes A)$
is inverse to $H^l_{X,\un}$. Therefore (iii) implies (i). Hence Part (a) of the proposition, since (i) implies (ii) is trivial.

Let us prove Part (b). Denote by $H^r$ the right fusion operator of $T$ and $\HAr$ the opfusion operator of $A$. Since $H^r_{X,Y}=H^r_{X,\un} \otimes Y$, the bimonad $T$ is a right Hopf monad if and only if it is a right \uHopf\ monad. Hence (i') is equivalent to~(ii').
 Moreover, we have: $H^r_{X,\un}=(A \otimes \sigma_X)(\mathbb{H'} \otimes X)$. If (iii') holds, then $\sigma$ and $\mathbb{H}'$ are invertible by Lemma~\ref{lem-anti}, and so $H^r_{-,\un}$ is an isomorphism. Hence (iii') implies (ii'). Conversely, if $H^r_{-,\un}$ is an isomorphism, then in particular $\mathbb{H'}=H^r(\un,\un)$
is invertible, and $A \otimes \sigma$ is invertible. Since $\un$ is a retract of $A$, this implies that $\sigma$ is invertible. Hence (ii') implies (iii'). This completes the proof of Part (b).
%

In particular $T$ is a Hopf monad if and only if $\sigma$ is invertible and $(A,\sigma)$ admits
an antipode and an opantipode, in other words, $(A,\sigma)$ is a Hopf algebra in $\zz(\cc)$.
\end{proof}

%

\subsection{Characterization of representable Hopf monads}\label{sect-char-rep}
Let $\cc$ be a monoidal category. A bimonad $T$ on $\cc$ is \emph{augmented} if it is endowed with an \emph{augmentation}, that is, a bimonad morphism $e\co T \to 1_\cc$.

Augmented bimonads on $\cc$ form a category
$\BiMon(\cc)/1_\cc$, whose objects are augmented bimonads on $\cc$, and morphisms between two augmented bimonads $(T,e)$ and $(T',e')$ are morphisms of bimonads $f\co T \to T'$ such that $e'f=e$.

If $(A,\sigma)$ is a  lax central bialgebra of $\cc$,   the bimonad $A\otimes_\sigma ?$ (see Section~\ref{sect-rep-cent}) is augmented, with augmentation $e= \varepsilon \otimes ? \co A \otimes_\sigma ? \to 1_\cc$, where $\varepsilon$ is the counit of $(A,\sigma)$.  Hence a functor $\BiAlg(\zzl(\cc)) \to \BiMon(\cc)/1_\cc$ which, according to Proposition~\ref{prop-Halg-to-Hmon}, induces  by restriction a functor
$$\mathfrak{R} \co \left\{ \begin{array}{ccc} \HopfAlg(\zz(\cc)) & \to & \HopfMon(\cc)/1_\cc \\ (A,\sigma) & \mapsto &(A\otimes_\sigma ?,\varepsilon \otimes ?)\end{array} \right.$$
where $\HopfMon(\cc)/1_\cc$ denotes the category of augmented Hopf monads on $\cc$.

\begin{thm}\label{thm-rep-HopfMon}
The functor $\mathfrak{R}$ is an equivalence of categories.
\end{thm}

In other words, Hopf monads representable by  central Hopf algebras    are nothing but augmented Hopf monads. Theorem~\ref{thm-rep-HopfMon} is proved in Section~\ref{sect-proofthemrep}.

\begin{rem}\label{rem-notrep}
Hopf monads are not representable in general. A counterexample is given in \cite{BV3} in terms of centralizers. Let $T$ be a centralizable Hopf monad on an autonomous category $\cc$ (see Example~\ref{exa-double}). In general the centralizer $Z_T$ of $T$ is not representable by a Hopf algebra. For example,
let $\cc=G\ti \vect$ be the category of finite-dimensional $G$\ti graded vector spaces over a field $\kk$ for some finite group $G$. The identity $1_\cc$ of $\cc$ is centralizable
and its centralizer $Z_{1_\cc}$ is representable if and only if $G$ is abelian (see~\cite[Remark 9.2]{BV3}).
\end{rem}

Hopf monads on a braided category $\bb$ which are representable by Hopf algebras in $\bb$ can also be characterized as follows:

\begin{cor}\label{cor-rep-braided}
Let $T$ be a Hopf monad on a braided category $\bb$. Then $T$ is isomorphic to the Hopf monad $A \otimes ?$ for some Hopf algebra $A$ in $\bb$
if and only if it is endowed with an augmentation $e \co T \to 1_\cc$ compatible with the braiding $\tau$ of $\bb$ in the following sense:
$$(e_X \otimes T\un)T_2(X,\un)=(e_X \otimes T\un)\tau_{T\un,TX}T_2(\un,X)$$
for all object $X$ of $\bb$.
\end{cor}
The corollary is proved in Section~\ref{sect-proofthemrep}.

\begin{rem}
Let $T$ be a centralizable Hopf monad on a braided autonomous~$\bb$ (see Remark~\ref{rem-notrep}). Then the centralizer $Z_T$ of $T$
is representable by a Hopf algebra $C_T=\int^{Y \in \bb} \ldual{T(Y)} \otimes Y$ in $\bb$ (see \cite[Theorem 8.4]{BV3}). This representability result may be recovered from Corollary~\ref{cor-rep-braided}, observing that
$$e_X=\int^{Y \in \cc} (\lev_{Y} \otimes X)(\ldual{\eta_Y} \otimes \tau_{Y,X}^{-1})\co Z_T(X) \to X$$ defines an augmentation of $Z_T$ which is compatible with the braiding $\tau$ of $\bb$.
\end{rem}

\subsection{Bosonization}\label{sect-boson} Let $\cc$ be a monoidal category.
Given a Hopf monad $(T,\mu,\eta)$ on $\cc$ and a  central Hopf algebra $(\AT,\sigma)$ of $\mo{T}{\cc}$ (that is, a Hopf algebra in the center
$\zz(\mo{T}{\cc})$ of $\mo{T}{\cc}$),    set:
$$
\AT \cp_\sigma T=(\AT \otimes_\sigma? ) \cp  T.
$$
As a cross product of Hopf monads, $\AT \cp_\sigma T$ is a Hopf monad on $\cc$ (see Proposition~\ref{prop-cp-hopf}). Set $\AT=(A,a)$, where $A=U_T(\AT)$ and $a$ is the $T$\ti action on $A$. As an endofunctor of $\cc$, $\AT \cp_\sigma T=A \otimes T$. The product $p$ and unit $v$ of $\AT \cp_\sigma T$ are:
$$p_X=\psfrag{A}[Bc][Bc]{\scalebox{.8}{$A$}} \psfrag{U}[Bc][Bc]{\scalebox{.8}{$TX$}} \psfrag{V}[Bc][Bc]{\scalebox{.8}{$T(A \otimes TX)$}} \psfrag{a}[Bc][Bc]{\scalebox{1.1}{$a$}}  \psfrag{m}[Bc][Bc]{\scalebox{1}{$\mu_X$}} \psfrag{T}[Bc][Bc]{\scalebox{.8}{$T_2(A,TX)$}}\rsdraw{.45}{.9}{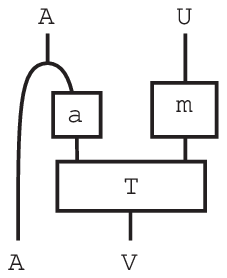} \quad \text{and} \quad  v_X=\psfrag{A}[Bc][Bc]{\scalebox{.8}{$A$}} \psfrag{X}[Bc][Bc]{\scalebox{.8}{$X$}} \psfrag{U}[Bc][Bc]{\scalebox{.8}{$TX$}} \psfrag{T}[Bc][Bc]{\scalebox{1}{$\eta_X$}} \,\rsdraw{.45}{.9}{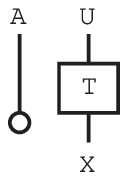}\;.$$
The comonoidal structure of $\AT \cp_\sigma T$ is given by:
$$(\AT \cp_\sigma T)_2(X,Y)=\psfrag{A}[Bc][Bc]{\scalebox{.8}{$A$}} \psfrag{U}[Bc][Bc]{\scalebox{.8}{$TX$}} \psfrag{Y}[Bc][Bc]{\scalebox{.8}{$TY$}}\psfrag{V}[Bc][Bc]{\scalebox{.8}{$T(X\otimes Y)$}}\psfrag{T}[Bc][Bc]{\scalebox{.8}{$T_2(X,Y)$}}\psfrag{X}[Bc][Bc]{\scalebox{1}{$\sigma_{F_T(X)}$}} \rsdraw{.45}{.9}{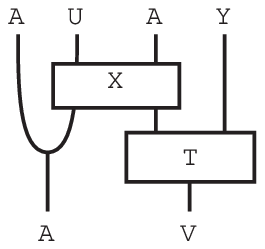} \quad \text{and} \quad  (\AT \cp_\sigma T)_0=\psfrag{A}[Bc][Bc]{\scalebox{.8}{$A$}} \psfrag{X}[Bc][Bc]{\scalebox{.8}{$T(\un)$}} \psfrag{T}[Bc][Bc]{\scalebox{.9}{$T_0$}} \rsdraw{.45}{.9}{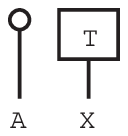}\;.$$
Denoting by $u$ and $\varepsilon$ the unit and counit of $(\AT,\sigma)$, the morphisms $$\iota=u \otimes T  \co T \to \AT \cp_\sigma T \quad \text{and}\quad \pi =\varepsilon \otimes T \co \AT \cp_\sigma T \to T$$ are Hopf monads  morphisms such that $\pi\iota=\id_T$. Hence~$T$ is a retract of $\AT \cp_\sigma T$ in the category $\HopfMon(\cc)$ of Hopf monads on $\cc$.
\begin{exa}
Let $T$ be a centralizable Hopf monad on a autonomous category~$\cc$ and $D_T$ be the double of $T$ (see Example~\ref{exa-double}). If $T$ is quasitriangular (see \cite{BV2}), then $\mo{T}{\cc}$ is braided and $T$ is a retract of $D_T$. In that case, the braided category $\mo{T}{\cc}$ admits a coend $C$, which is a Hopf algebra, and $D_T= C \cp_{\tau_{C,-}} T$ where $\tau$ is the braiding of $\mo{T}{\cc}$.
\end{exa}

Conversely, under exactness assumptions, a Hopf monad which admits $T$ as a retract is of the form $\AT \cp_\sigma T$ for some  central Hopf algebra $(\AT,\sigma)$ of $\mo{T}{\cc}$.   This results from the fact that augmented Hopf monads are representable, using the notion of cross quotient studied in Section~\ref{sect-croco}:

\begin{cor}\label{thm-rel-rep-HopfMon} Let $P$ and $T$ be Hopf monads on a monoidal category $\cc$ such that $T$ is a retract of $P$.
Assume that reflexive coequalizers exist in $\cc$ and are preserved by $P$ and the monoidal product of $\cc$.
Then
there exists a  central Hopf algebra $(\AT,\sigma)$ of $\mo{T}{\cc}$    and an isomorphism of Hopf monads $P \simeq \AT \cp_\sigma T$ such that we have a commutative diagram of Hopf monads:
$$
\xymatrix@!0 @R=3em@C=2em{
P \ar[rrrr]^-{\simeq}\ar[rrrd]&&&&\AT \cp_\sigma T  \ar[ld]\\
&T \ar[lu] \ar[urrr]|!{[ul];[rr]}\hole \ar[rr]_-{=}&&T&
}
$$
\end{cor}

\begin{proof} Denote by $f \co T \to P$ and $g \co P \to T$ the  morphisms of Hopf monads making $T$ a retract of $P$,  that is, $gf=\id_T$.
By assumption $P$ preserves reflexive coequalizers and so, since $T$ is a retract of $P$, the Hopf monad $T$ preserves reflexive coequalizers too. By Lemma~\ref{lem-tens-pres-coeq}, reflexive coequalizers exist in $\mo{T}{\cc}$ and $\mo{P}{\cc}$ and are preserved by the monoidal product. By Proposition~\ref{prop-P/T-Hopf},
the crossed quotient $P \cdiv T$ (relative to $f \co T \to P$) exists and is a Hopf monad on $\mo{T}{\cc}$. On the other hand, by functoriality of the cross quotient (see Remark~\ref{rem-funct-crossquot}), $g \co P \to T$ induces a morphism of bimonads $g \cdiv T \co P\cdiv T \to T\cdiv T\cong 1_{\mo{T}{\cc}}$. In other words the Hopf monad $P\cdiv T$ is augmented.
By Theorem~\ref{thm-rep-HopfMon}, there exists a Hopf algebra $(\AT,\sigma)$ in $\zz(\mo{T}{\cc})$ such that $P\cdiv T = \AT \otimes_\sigma ?$.
By Proposition~\ref{prop-prod-quot}, $P=(P\div T) \cp T=(\AT \otimes_\sigma ?) \cp T=\AT \otimes_\sigma T$ as Hopf monads. The commutativity of the diagram is straightforward.
\end{proof}

\begin{rem}
Let $H$ be a Hopf algebra over a field $\kk$, and $A$ a Hopf algebra in the braided category of Yetter-Drinfeld modules ${}^H_H\mathcal{YD}$. In that situation, Radford
constructed a Hopf algebra $A \# H$, known as \emph{Radford's biproduct}, or \emph{Radford-Majid bosonization}.
Radford~\cite{Rad} (see also~\cite{Maj3}) showed that if  $K$ is a Hopf algebra on a field $\kk$ and $p$ is a projection of $K$, that is, an idempotent endomorphism of the Hopf algebra $K$, then $K$ may be described as a biproduct as follows. Denote by $H$ the image of $p$, which is a Hopf subalgebra of $K$. Then there exists a
Hopf algebra $A$ in ${}^H_H \mathcal{YD}$ such that $K=A \#H$.
Corollary~\ref{thm-rel-rep-HopfMon} generalizes Radford's theorem. Indeed, in the situation of the theorem, the Hopf monad $H \otimes ?$ is a retract of the Hopf monad $K\otimes ?$ on $\Vect_\kk$. Hence, by Corollary~\ref{thm-rel-rep-HopfMon},  there exists  a Hopf algebra $(\AT,\sigma)$ in $\zz(\biMod{H}{})$ such that $K \otimes = \AT \cp_\sigma (H \otimes ?)$. Identifying the center of $\biMod{H}{}$ with the category of Yetter-Drinfeld modules, we view $(\AT,\sigma)$ as a Hopf algebra $A$ in ${}^H_H\mathcal{YD}$. Then $K \otimes ?=\AT \cp_\sigma (H \otimes ?) =A \# H \otimes ?$ as Hopf monad,  and so $K=A \# H$.
\end{rem}

\subsection{Regular augmentations}\label{sect-proofthemrep} In this section we prove Theorem~\ref{thm-rep-HopfMon}  and Corollary~\ref{cor-rep-braided}
using the notion of regular augmentation.

Let $T$ be a comonoidal endofunctor of a monoidal category $\cc$ and $e \co T \to 1_\cc$ be a comonoidal natural transformation.
Define natural transformations $u^e \co T \to T\un\, \otimes ?$
and $v^e \co T \to ? \otimes T\un$ by:
$$
u^e_X=(T\un \otimes e_X)T_2(\un,X)\quad \text{and} \quad
v^e_X=(e_X \otimes T\un)T_2(X,\un).
$$
We say that $e$ is \emph{left regular} if $u^e$ is invertible.

\begin{lem} \label{lem-halfbraiding}
Assume $e$ is left regular and set $$\sigma=v^e (u^e)^{-1} \co T\un \otimes ? \to ? \otimes T\un.$$ Then the natural transformation $\sigma$ is a lax half braiding in $\cc$ and $(T\un,\sigma)$
 is a  lax central coalgebra  of $\cc$   with coproduct $T_2(\un,\un)$ and counit $T_0$. Furthermore the natural transformation $u^e \co T \to T\un\, \otimes_\sigma ?$ is a comonoidal isomorphism.
\end{lem}

\begin{proof}
By transport of structure, the endofunctor $P=T\un \,\otimes ?$ of $\cc$ admits a unique comonoidal structure such that the natural isomorphism
$u^e\co T \to P$ is comonoidal, that is,
$$
P_0 u^e_\un=T_0 \quad \text{and} \quad P_2(X,Y)u^e_{X\otimes Y}=(u^e_X \otimes u^e_Y)T_2(X,Y).
$$
We have $e_\un=T_0$ (since $e$ is comonoidal) and so
$u^e_\un=\id_{T\un}$ and $v^e_\un=\id_{T\un}$. Hence $P_0= P_0 u^e_\un=T_0$ and $\sigma_\un=v^e_\un (u^e_\un)^{-1}=\id_{T\un}$. Moreover,
\begin{align*}
P_2(X,Y)& u^e_{X\otimes Y}=(u^e_X \otimes u^e_Y)T_2(X,Y)\\
&=(T\un \otimes e_X \otimes T\un \otimes e_Y)T_4(\un,X,\un,Y)\\
&=(T\un \otimes v^e_X \otimes Y)(T\un \otimes TX \otimes e_Y)T_3(\un,X,Y)\\
&=(T\un \otimes v^e_X \uei_X\otimes Y)(T\un \otimes u_X^e \otimes e_Y)T_3(\un,X,Y)\\
&=(T\un \otimes \sigma_X \otimes Y)(T\un \otimes T\un \otimes e_X \otimes e_Y)T_4(\un,\un,X,Y)\\
& =(T\un \otimes \sigma_X \otimes Y)(T_2(\un,\un) \otimes (e_X \otimes e_Y)T_2(X,Y))T_2(\un,X\otimes Y)\\
& =(T\un \otimes \sigma_X \otimes Y)(T_2(\un,\un) \otimes X \otimes Y )u^e_{X\otimes Y}
\quad\text{(since $e$ is comonoidal).}
\end{align*}
Therefore
$P_2(X,Y)=(T\un \otimes \sigma_X \otimes Y)(T_2(\un,\un) \otimes X\otimes Y)$
because $u^e$ is invertible. We conclude using Lemma~\ref{lem-ATbimon}.
\end{proof}

Recall that an augmentation of a bimonad $T$ on $\cc$ is a morphism of bimonads from $T$ to $1_\cc$. It is called \emph{left regular} if it is left regular as a comonoidal natural transformation.

\begin{lem} \label{prop-centbialg}
Let $(T,\mu,\eta)$ be an augmented bimonad on $\cc$. Assume its augmentation $e\co T \to 1_\cc$ is left regular. Then $\sigma=v^e (u^e)^{-1}$ is a lax half braiding for $T\un$ and $(T\un,\sigma)$ is a  lax central bialgebra of $\cc$,   with product $m=\mu_\un (u^e_{T\un})^{-1}$, unit $u=\eta_\un$, coproduct $T_2(\un,\un)$, and counit $T_0$. Moreover
 $u^e \co T \to T\un \otimes_\sigma ?$ is an isomorphism of bimonads.
\end{lem}

\begin{proof} By transport of structure, the endofunctor $P=T\un \,\otimes ?$ of $\cc$ admits a unique bimonad structure such that the natural transformation $u_e\co T \to P $ is an isomorphism of bimonads.
In view of Lemmas~\ref{lem-ATbimon} and~\ref{lem-halfbraiding}, it is enough to verify that the product
$\mu'$ and unit $\eta'$ of $P$ are given by
$\mu'=\mu_\un (u^e_{T\un})^{-1} \otimes ? \blabla{and} \eta'=\eta_\un \otimes ?$. 
Since $u^e$ is a morphism of monads, we have:
$$
\eta'_X=u^e_X \eta_X=(T\un \otimes e_X)T_2(\un,X)\eta_X=\eta_\un \otimes e_X\eta_X=\eta_\un \otimes X.
$$
Also, setting $m=\mu_\un (u^e_{T\un})^{-1}$, we have:
\begin{align*}
\mu'_Xu^e_{T\un \otimes X}T(u^e_X) &= u^e_X \mu_X  
=(T\un \otimes e_X)T_2(\un,X)\mu_X\\
&=(\mu_\un \otimes e_X \mu_X)T^2_2(\un,X)\\
&=(m u^e_{T\un} \otimes e_X T(e_X))T_2(T\un,TX)T(T_2(\un,X))\\
&=(m \otimes X)\bigl((T\un \otimes e_{T\un})T_2(\un,T\un) \otimes e_X \bigr)T_2(T\un,X)T(u^e_X)\\
&=(m \otimes X)(T\un \otimes e_{T\un \otimes X})T_2(\un,T\un \otimes X)T(u^e_X)\\
&=(m \otimes X)u^e_{T\un \otimes X}T(u^e_X),
\end{align*}
and so $\mu'_X=m \otimes X$ since $u^e$ is invertible.
\end{proof}

\begin{lem}\label{lem-hopf-reg}  Let $T$ be an augmented left \uHopf\ monad on $\cc$. Then its augmentation $e \co T \to 1_\cc$ is left regular and $(u^e)^{-1}=T(e) \Hli_{\un,-}(T\un \otimes \eta)$.
\end{lem}
\begin{proof}  Let $X$ be an object of $\cc$ and set $\theta^e_X=T(e_X) \Hli_{\un,X}(T\un \otimes \eta_X)$.
We have:
\begin{align*}
u^e_X\theta^e_X &=(T\un \otimes e_X)T_2(\un,X)T(e_X) \Hli_{\un,X}(T\un \otimes \eta_X)\\
&=(T\un \otimes e_XT(e_X))T_2(\un,TX)\Hli_{\un,X}(T\un \otimes \eta_X)\\
&=(T\un \otimes e_X \mu_X) T_2(\un,TX)\Hli_{\un,X}(T\un \otimes \eta_X)\\
&=(T\un \otimes e_X)H^l_{\un,X}\Hli_{\un,X}(T\un \otimes \eta_X)\\
&=(T\un \otimes e_X\eta_X)=\id_{T\un \otimes X}
\end{align*}
and
\begin{align*}
\theta^e_X u^e_X &=T(e_X) \Hli_{\un,X}(T\un \otimes \eta_X e_X )T_2(\un,X)\\
&=T(e_X) \Hli_{\un,X}(T\un \otimes T(e_X)\eta_{TX})T_2(\un,X)\\
&=T(e_X) T^2(e_X) \Hli_{\un,TX}(T\un \otimes \eta_{TX})T_2(\un,X)\\
&=T(e_X) T(\mu_X) \Hli_{\un,TX}(T\un \otimes \eta_{TX})T_2(\un,X)\\
&=T(e_X) \Hli_{\un,X}(T\un \otimes \mu_X\eta_{TX})T_2(\un,X) \quad\text{by Proposition~\ref{prop-fusion}}\\
&=T(e_X) \Hli_{\un,X}(T\un \otimes \mu_XT(\eta_X))T_2(\un,X)\\
&=T(e_X) \Hli_{\un,X}H^l_{\un,X}T(\eta_X)=T(e_X\eta_X)=\id_{TX}.
\end{align*}
Hence $u^e$ is invertible with inverse $\theta^e$.
\end{proof}

\begin{thm}\label{thm-rep} Let $\cc$ be a monoidal category. The functor
$$\mathfrak{R}^\mathrm{lax} \co \left\{ \begin{array}{ccc} \BiAlg(\zzl(\cc)) & \to & \BiMon(\cc)/1_\cc \\ (A,\sigma) & \mapsto &(A\otimes_\sigma ?,\varepsilon \otimes ?)\end{array} \right.$$
induces an equivalence of categories from $\BiAlg(\zzl(\cc))$ to the full subcategory of $\BiMon(\cc)/1_\cc$ of augmented bimonads $(T,e)$ such that $e$ is left regular.
\end{thm}

\begin{proof}
If $(A,\sigma)$ is a bialgebra in $\zzl(\cc)$, then $e=\varepsilon \otimes ? \co A \otimes_\sigma ? \to 1_\cc$ is a left regular bimonad morphism (since $u^e=\id_{A \otimes ?}$).
Therefore $\mathfrak{R}^\mathrm{lax}$ takes values in the full subcategory $\aaa \subset \BiMon(\cc)/1_\cc$ of augmented bimonads $(T,e)$ such that $e$ is left regular. Conversely, let $(T,e)$ be an object of $\aaa$. By Lemma~\ref{prop-centbialg}, $T\un$ is endowed with a half-braiding $\sigma$ and $(T\un,\sigma)$ is a bialgebra in $\zz(\cc)$.  This construction is functorial, that is, gives rise to a functor $\mathfrak{I}\co  \aaa \to \BiAlg(\zzl(\cc))$ defined on objects by $\mathfrak{I}(T,e)=(T\un,\sigma=v^e u^{e-1})$ and on morphisms by $\mathfrak{I}(f)=f_\un$. Moreover $\mathfrak{I}$ is quasi-inverse to $\mathfrak{R}^\mathrm{lax}$. Indeed, for $(T,e)$ in $\aaa$, $u^e$ is an isomorphism from $(T,e)$ to $\mathfrak{R}^\mathrm{lax}\mathfrak{I}(T,e)$ and, for $(A,\sigma)$ in $\BiAlg(\zzl(\cc))$, we have $\mathfrak{I}\mathfrak{R}^\mathrm{lax}(A,\sigma)=(A,\sigma)$. Hence
the Theorem.
\end{proof}

\begin{cor}\label{cor-rep}
Let $\cc$ be a monoidal category. The functor $\mathfrak{R}^\mathrm{lax}$ induces equivalences of categories between:
\begin{enumerate}
\labela
\item  Lax central left Hopf algebras of $\cc$   and augmented left Hopf monads on $\cc$.
\item  Central Hopf algebras of $\cc$   and augmented Hopf monads on $\cc$.
\end{enumerate}
Moreover an augmented left (resp.\@ right) \uHopf\ monad on $\cc$ is in fact a left (resp.\@ right) Hopf monad.
\end{cor}

\begin{proof} Let $(T,e)$ be an augmented bimonad such that $T$ is a left \uHopf\ monad. Then $e$ is left regular by Lemma~\ref{lem-hopf-reg}. By Theorem~\ref{thm-rep}, $T$ is of the form $A \otimes_\sigma ?$ for some bialgebra $(A,\sigma)$ in $\zzl(\cc)$. By Proposition~\ref{prop-Halg-to-Hmon}, $A$ admits an antipode, and $T$ is in fact a left Hopf monad.
Hence the first equivalence of categories.
Moreover, by Proposition~\ref{prop-Halg-to-Hmon}, $T$ is a Hopf monad if and only if $(A,\sigma)$ is a Hopf algebra in $\zz(\cc)$. Hence the second equivalence of categories.%
\end{proof}

\begin{proof}[Proof of Theorem~\ref{thm-rep-HopfMon}]
The theorem is just Assertion (b) of Corollary~\ref{cor-rep}.
\end{proof}

\begin{proof}[Proof of Corollary~\ref{cor-rep-braided}]
By Theorem~\ref{thm-rep-HopfMon}, the augmentation $e \co T \to 1_\bb$
defines a Hopf algebra $(A=T\un,\sigma)$ in $\zz(\bb)$ such that $T \simeq A \otimes_\sigma ?$. In view of Remark~\ref{rem-rep-braided}, the question is whether $\sigma=\tau_{A,-}$.
Recall $\sigma=v^e(u^e)^{-1}$. Therefore $\sigma_X=\tau_{A,X}$ if and only if $v^e_X=\tau_{A,X}u^e_X$, that is,  $(e_X \otimes T\un)T_2(X,\un)=(e_X \otimes T\un)\tau_{T\un,TX}T_2(\un,X)$.
\end{proof}

\section{Induced coalgebras and Hopf modules}\label{sect-Hmod}
A cocommutative coalgebra of the center of a monoidal category $\dd$ gives rise to a comonoidal comonad on $\dd$ and, under certain
exactness assumptions, to a Hopf adjunction. On the other hand, we show that the comonoidal comonad of a \uHopf\ adjunction
$(F \co \cc \to \dd, U \co \dd \to \cc)$ is represented by its \emph{induced coalgebra}, which is a cocommutative coalgebra of the categorical center of $\dd$.

As an application, we obtain a structure theorem for Hopf modules over \uHopf\ monads on monoidal categories. It generalizes Sweedler's Theorem on the structure of Hopf modules over a Hopf algebra, and is an enhanced version of \cite[Theorem~4.6]{BV2} which concerns Hopf monads on autonomous categories.

%
%
%

\subsection{From cocommutative central coalgebras to Hopf adjunctions}
Let $\dd$ be a monoidal category and $(C,\Delta,\varepsilon)$  be a coalgebra in $\dd$.
Denote by $\lComod{C}$ the category of left $C$\ti comodules in~$\dd$. The forgetful functor $V \co \lComod{C} \to \dd$ has a right adjoint, the cofree comodule functor
$$R\co \prettydef{
\dd &\to& \lComod{C}\\
X &\mapsto& (C \otimes X, \Delta \otimes X)
} .$$

The comonad (on $\dd$)  of the adjunction $(V,R)$ is $\hat{T}=(C \otimes ?,\Delta \otimes ?,\varepsilon \otimes ?)$. This adjunction is comonadic, since $\hat{T}$\ti comodules are
just left $C$\ti comodules.

On the other hand, the monad $T=RV$ (on $\lComod{C}$) of the adjunction  $(V,R)$ is defined by $T(M,\delta)=(C \otimes M, \Delta \otimes M)$ for any $C$\ti comodule $(M,\delta)$, with product $\mu_{(M,\delta)}=C \otimes \varepsilon \otimes M$ and unit $\eta_{(M,\delta)}=\delta$.
In general the adjunction $(V,R)$ is not monadic.

\begin{rem}\label{rem-monadic} The adjunction $(V,R)$ is monadic if, for instance, $\dd$ admits reflexive coequalizers  and
$C \otimes ?$ is conservative and preserves reflexive coequalizers.
\end{rem}

 Now let $(C,\sigma)$ be a lax central  coalgebra of $\dd$, that is, a coalgebra in $\zzl(\cc)$.
Then   the endofunctor $C \otimes ?$ of $\dd$ has both a comonad structure (because $C$ is a coalgebra in $\dd$) and a comonoidal structure denoted by $C\otimes_\sigma ?$ (see Lemma~\ref{lem-ATbimon}).

 A lax central coalgebra $(C,\sigma)$ is \emph{cocommutative} if its coproduct $\Delta$ satisfies $\sigma_C\Delta=\Delta$. We have:
\begin{lem}\label{lem-cocomcomonad} Let $(C,\sigma)$ be a  lax central coalgebra of $\dd$.   Then  $C \otimes_\sigma ?$  is
a comonoidal comonad if and only if $(C,\sigma)$ is cocommutative.
\end{lem}

\begin{proof}   One checks that  the coproduct $\Delta \otimes ?$ of the comonad $C \otimes ?$ is comonoidal if and only if $\sigma_C \Delta=\Delta$, and  that its counit $\varepsilon \otimes ?$ is always comonoidal.
\end{proof}

We say that a cocommutative  lax central coalgebra $(C,\sigma)$ of $\dd$ is  \emph{cotensorable}    if
for each pair $(M,\delta)$, $(N,\delta')$ of left $C$\ti comodules, the coreflexive pair
$$\xymatrix@C=4em{M \otimes N \dar{\sigma_M \delta \otimes N}{M \otimes \delta'}&  \ M \otimes C \otimes N}$$ admits an equalizer, denoted by $M\ct{\sigma}{C} N \to M \otimes N$, and the endofunctor $C \otimes ?$ preserves these equalizers.

Let $(C,\sigma)$ be a  cotensorable cocommutative lax central coalgebra of $\dd$.   Given  two left $C$\ti comodules $(M,\delta)$ and $(N,\delta')$, there exists a unique left coaction $\delta''$ of $C$ on $M \ct{\sigma}{C} N$ such that the  morphism $M\ct{\sigma}{C} N \to M \otimes N$ is a comodule morphism
$(M \ct{\sigma}{C} N, \delta'') \to (M \otimes N, \delta \otimes N)$. The assignment $(M,\delta) \times (N,\delta') \mapsto (M\ct{\sigma}{C} N,\delta'')$ defines a
functor: $$\ct{\sigma}{C} \co \lComod{C} \times \lComod{C} \to \lComod{C}.$$
Then the category $\lComod{C}$ of left $C$\ti comodules (in $\cc$) is monoidal, with monoidal product $\ct{\sigma}{C}$ and unit object $(C,\Delta)$. We denote this monoidal category by $\lsComod{\sigma}{C}$. The cofree comodule functor $R \co \dd \to \lsComod{\sigma}{C}$ is strong monoidal, so that the comonadic adjunction $(V,R)$ is comonoidal, with comonoidal comonad $C \otimes_\sigma ?$.

\begin{lem}\label{lem-cotens} Let $\dd$ be a monoidal category admitting coreflexive equalizers which are preserved by the monoidal product. Let $(C,\sigma)$ be a cocommutative central coalgebra of  $\dd$. Then $(C,\sigma)$ is cotensorable and the monoidal category $\lsComod{\sigma}{C}$ admits coreflexive equalizers which are preserved by the monoidal product $\otimes^\sigma_C$ and the forgetful functor $V \co \lsComod{\sigma}{C} \to \dd$.
\end{lem}

\begin{proof} By standard diagram chase left to the reader. \end{proof}

\begin{thm}\label{thm-coalg-hopfadj}
Let $\dd$ be a monoidal category and $(C,\sigma)$ be a  cotensorable cocommutative lax central coalgebra of $\dd$.
Then the comonoidal adjunction $$(V \co \lsComod{\sigma}{C} \to \dd, R \co \dd \to \lsComod{\sigma}{C})$$
is a left Hopf adjunction,  and its induced lax central coalgebra is $(C,\sigma)$.   Moreover, if $\sigma$ is invertible, $(V,R)$ is a Hopf adjunction.
\end{thm}

\begin{proof}
Let $d$ be an object of $\dd$ and $(M,\delta)$ be a left $C$\ti comodule.
Then the morphism $\sigma_M\delta \otimes d\co M \otimes d \to M \otimes C \otimes d$ is an equalizer of the pair
$$\xymatrix@C=6em{M \otimes C \otimes d \dar{\sigma_M\delta \otimes C \otimes d}{M \otimes \Delta \otimes d}& M \otimes C \otimes C \otimes d}.$$
Hence an isomorphism $M \ct{\sigma}{C} R(d) \iso M \otimes d$ which is the left fusion operator of the comonoidal adjunction $(V,R)$.
Similarly if $\sigma$ is invertible, the morphism
$(\sigma^{-1}_d \otimes M)(d \otimes \delta)\co d \otimes M \to C \otimes d \otimes M$ is an equalizer of the pair
$$\xymatrix@C=8em{C \otimes d \otimes M \dar{(\sigma_{C \otimes d}\otimes M)(\Delta \otimes d \otimes M)}{C \otimes d \otimes \delta}& C \otimes d \otimes C \otimes M}.$$
Hence an isomorphism $R(d) \ct{\sigma}{C} M \iso d \otimes M$ which is the right fusion operator of the comonoidal adjunction $(V,R)$.
\end{proof}

\subsection{Induced coalgebra and comonad of a comonoidal adjunction}
Let $\cc$, $\dd$ be monoidal categories and $(F \co \cc \to \dd, U \co \dd \to \cc)$ be a comonoidal adjunction,
with adjunction unit $\eta \co 1_\cc \to UF$ and counit $\varepsilon \co FU \to 1_\dd$.

 Being comonoidal, $F$ sends the trivial coalgebra $\un$ in~$\cc$ to a coalgebra $\hat{C}=F(\un)$ in $\dd$, with coproduct
$\Delta=F_2(\un,\un)$ and counit $\epsilon=F_0$, called the
\emph{induced coalgebra of the comonoidal adjunction}. 

The endofunctor $\hat{T}=FU$ of $\dd$ is a comonoidal comonad, with coproduct $F(\eta_U) \co \hat{T}\to \hat{T}^2$ and counit $\varepsilon$ (see Section~\ref{sect-bimonadadj}).


In this situation we have three comonads on the category $\dd$, namely:
\begin{itemize}
\item $? \otimes \hat{C}$ (with coproduct $? \otimes \Delta$ and counit $? \otimes \epsilon$);
\item $\hat{C} \otimes \,?$ (with coproduct $\Delta \otimes \,?$ and counit $\epsilon \otimes \,?$);
\item the (comonoidal) comonad $\hat{T}=FU$ of the adjunction $(F,U)$.
\end{itemize}
 How are they related?

\begin{lem}\label{lem-comon-coalg}
The  Hopf operators $\FO^l$ and $\FO^r$ define
morphisms of comonads:
$$
\FO^l_{\un,-} \co \hat{T} \to
\hat{C} \otimes ?\quad\text{and}\quad
\FO^r_{-,\un} \co \hat{T} \to ? \otimes \hat{C}.$$
\end{lem}

\begin{proof} The commutativity of the following diagrams:
\begin{equation*}
\xymatrix@C=4.5pc@R=1pc{
FUd \ar[ddd]_{F_2(\un,Ud)}\ar[r]^{F\eta_{Ud}}\ar[rdd]_{F_2(\un,Ud)}
&FUFUd\ar[d]^{F_2(\un,UFUd)}\\
&F\un\otimes FUFUd\ar[d]^{F\un \otimes \varepsilon_{FUd}}\\
&F\un\otimes FUd\ar[d]^{F\un \otimes F_2(\un,Ud)}\\
F\un\otimes FUd\ar[d]_{F\un \otimes \varepsilon_d}\ar[r]^-{F_2(\un,\un) \otimes FUd}
&F\un\otimes F\un \otimes FUd\ar[d]^{F\un \otimes F\un \otimes \varepsilon_d}\\
F\un\otimes d\ar[r]_-{F_2(\un,\un) \otimes d}
&F\un\otimes F\un \otimes d
}
\quad
\xymatrix{
FUd \ar[rd]^{\varepsilon_d}\ar[d]_{F_2(\un,Ud)}&\\
F\un \otimes FUd \ar[d]_{F\un \otimes \varepsilon_d}& d\\
F\un \otimes d \ar[ru]_{F_0 \otimes d}&}
\end{equation*}
which results from the fact that the adjunction $(F,U)$ is comonoidal, means that $\FO^l_{\un,-}$ is a morphism of comonads. The proof for $\FO^r_{-,\un}$ is similar.
\end{proof}

\subsection{From Hopf adjunctions to cocommutative central coalgebras}

In the case of a left \uHopf\ adjunction,  the induced coalgebra is endowed with a canonical lax half braiding, making it a cocommutative lax central coalgebra which represents the induced comonoidal comonad:

\begin{thm}\label{prop-comon-coalg}
Let $(F \co \cc \to \dd,U\co \dd \to \cc)$ be a left \uHopf\ adjunction, with induced coalgebra $\hat{C}$. Then:
\begin{enumerate}
\labela
\item The natural transformation $\hat{\sigma}=\FO^r_{-,\un} \FOli_{\un,-} \co \hat{C} \otimes ? \to ? \otimes \hat{C}$ is a lax half-braiding of $\dd$ such that,
 for any object $c$ of $\cc$, the diagram:
$$ \xymatrix@R=1.4em@C=.8em{ & Fc \ar[ld]_-{F_2(\un,c)}\ar[rd]^-{F_2(c,\un)}& \\
\hat{C} \otimes Fc  \ar[rr]^-{\hat{\sigma}_{Fc}}
 & & Fc \otimes \hat{C}
}$$
is commutative.
\item $(\hat{C},\hat{\sigma})$ is a cocommutative   lax central coalgebra of $\dd$   and $\FO^l_{\un,-} \co \hat{T} \to \hat{C} \otimes_{\hat{\sigma}} ?$ is an isomorphism of comonoidal comonads.
\end{enumerate}

\end{thm}

\begin{proof}
Let $(F,U)$ be a left \uHopf\ adjunction, so that ${\FO^l}_{\un,-}$ is invertible.  Since $U$ is strong monoidal, we identify $\hat{C}=F(\un)$ and $\hat{T}(\un)=FU(\un)$ as coalgebras in $\dd$.   We apply Lemma~\ref{lem-halfbraiding} to the comonoidal endofunctor $\hat{T}$ of $\dd$ and the comonoidal morphism
$\varepsilon \co FU \to 1_\cc$. The natural transformations $u^\varepsilon \co FU \to FU(\un) \otimes ?$ and $v^\varepsilon \co FU \to ? \otimes FU(\un)$ of the lemma are nothing but $\FO^l_{\un,-}$ and $\FO^r_{-,\un}$ respectively. Therefore $u^e=\FO^l_{\un,-}$ being invertible, we conclude that  $\hat{\sigma}=v^e (u^e)^{-1}$ is a lax braiding on $\dd$ and $(\hat{C},\hat{\sigma})$ is a coalgebra in $\zz(\dd)$ such that $u^e$ is a comonoidal isomorphism.
 Now, for any object $c$ of $\cc$, we have $$\FO^l_{\un,Fc} F(\eta_{c})=F_2(\un,c) \quad\text{and}\quad \FO^r_{Fc,\un} F(\eta_c)=F_2(c,\un),$$
from which the equality $\hat{\sigma}_{Fc}F_2(\un,c)=F_2(c,\un)$ follows directly.  Hence Part (a).

Applying  this equality  to $d=\un$ gives the cocommutativity of the coalgebra $(\hat{C},\hat{\sigma})$, so that $\hat{C} \otimes_{\hat{\sigma}} ?$ is a comonoidal comonad by Lemma~\ref{lem-cocomcomonad}. Thus $\FO^l_{\un,-}=u^e$ is an isomorphism of comonoidal comonads, hence Part (b).
\end{proof}

As a consequence, the comonoidal comonad of a \uHopf\ adjunction is canonically represented by a  cocommutative central coalgebra of  the upper category.   More precisely:

\begin{cor}\label{cor-Hopf-Adj-Rep}
Let $(F\co \cc \to \dd,U\co \dd \to \cc)$ be a \uHopf\ adjunction, with induced coalgebra $\hat{C}$. Then $\hat{\sigma}=\FO^r_{-,\un} \FOli_{\un,-}$ is a half-braiding for $\hat{C}$.  Moreover $(\hat{C}, \hat{\sigma})$ is a  cocommutative central coalgebra in  $\dd$ called the \emph{induced central coalgebra of the \uHopf\ adjunction $(F,U)$.}
\end{cor}
\begin{proof}
Since $(F,U)$ is a \uHopf\ adjunction, the \uHopf\ operators $\FO^r_{-,\un}$ and  $\FO^l_{\un,-}$ are invertible. Thus $\hat{\sigma}$ is invertible and the corollary follows then directly from Theorem~\ref{prop-comon-coalg}.
\end{proof}

\begin{exa}
Let $\cc$ be a monoidal category and $(A,\sigma)$ be a Hopf algebra in $\zz(\cc)$, with product $m$, coproduct $\Delta$, and counit $\varepsilon$. Consider the Hopf monad $T=A \otimes_\sigma ?$ on $\cc$ (see Proposition~\ref{prop-Halg-to-Hmon}). Recall $\lsMod{A}{\sigma}$ denotes the monoidal category $\mo{A \otimes_\sigma ?}{\cc}$ of left $A$\ti modules (in $\cc$), with monoidal product induced by $\sigma$ (see Section~\ref{sect-rep-cent}). The induced coalgebra $\hat{C}$ of $A \otimes_\sigma ?$ is the left $A$\ti module $\hat{C}=(A,m)$, with coproduct $\Delta$ and counit $\varepsilon$. Its associated half-braiding is given by
$$
\hat{\sigma}_{(M,r)}= \psfrag{A}[Bc][Bc]{\scalebox{.8}{$A$}} \psfrag{U}[Bc][Bc]{\scalebox{.8}{$M$}} \psfrag{H}[Bc][Bc]{\scalebox{1}{$\sigma_M$}}\psfrag{r}[Bc][Bc]{\scalebox{1}{$r$}}\psfrag{X}[Bc][Bc]{\scalebox{1}{$\sigma_A$}} \rsdraw{.45}{.9}{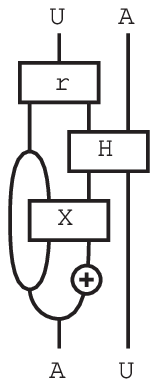}
$$
for any left $A$\ti module $(M,r)$.
Then $(\hat{C}, \hat{\sigma})$ is a cocommutative coalgebra in the center $\zz(\lsMod{A}{\sigma})$ of $\lsMod{A}{\sigma}$.
\end{exa}

 \begin{prop}\label{prop-HM-strong}
Let $(F\co \cc \to \dd,U\co \dd \to \cc)$ be a comonadic \uHopf\ adjunction,  with induced central coalgebra $(\hat{C},\hat{\sigma})$.
Assume that for all $X, Y$ objects of $\cc$, the morphism $F_2(X,Y) \co F(X \otimes Y) \to F(X) \otimes F(Y)$ is an equalizer of the coreflexive pair
$$\xymatrix@C=5em{FX \otimes FY \dar{F_2(X,\un) \otimes FY}{FX \otimes F_2(\un,Y)}&FX \otimes F\un \otimes FY},$$
and these equalizers are preserved by the endofunctor $F(\un) \otimes ?$.
Then the cocommutative central coalgebra $(\hat{C},\hat{\sigma})$ is cotensorable and the comparison functor
$$K \co \cc \to \lsComod{\hat{\sigma}}{\hat{C}}$$
 is a strong monoidal equivalence. In particular $(F,U)$ is a Hopf adjunction.
\end{prop}

\begin{proof} The cotensorability assumption means that for each pair $(M,\delta)$, $(N,\delta')$ of left $\hat{C}$\ti comodules, the coreflexive pair
$$\xymatrix@C=4em{M \otimes N \dar{\sigma_M \delta \otimes N}{M \otimes \delta'}&  \ M \otimes C \otimes N}$$ admits an equalizer, and the endofunctor $\hat{C} \otimes ?$ preserves these equalizers.
Now recall that the comparison functor $K$ is defined by $K(X)=(FX,F_2(\un,X))$ for $X$ in~$\cc$.
If $X$ is an object of $\cc$ then by Theorem~\ref{prop-comon-coalg}, Part (a), we have $\hat{\sigma}_{FX}F_2(\un,X)= F_2(X,\un)$.
Since $K$ is an equivalence, we conclude that $(\hat{C},\hat{\sigma})$ is cotensorable. Moreover, we have $K(X\otimes Y)=K(X) \otimes^{\hat{\sigma}}_{\hat{C}} K(Y)$ so that  $K$ is a strong monoidal equivalence. By Theorem~\ref{thm-coalg-hopfadj}, $(F,U)$ is a Hopf adjunction.
\end{proof}

\subsection{Descent}\label{sect-descent}

Let  $(T,\mu,\eta)$ be a monad on a category $\cc$.  Its adjunction $(F_T,U_T)$
has unit $\eta\co 1_\cc \to U_T F_T = T$ and has counit denoted by
$\varepsilon\co F_T U_T \to 1_{\mo{T}{\cc}}$.
Let $\hat{T}$ be the comonad of the adjunction $(F_T,U_T)$, that is, $\hat{T}=F_T U_T$ on $\mo{T}{\cc}$, with coproduct $\delta=F_T (\eta_{U_T})$ and counit $\varepsilon$.
Denote by $\HM(T)$ the category $\como{\hat{T}}{(\mo{T}{\cc})}$ of $\hat{T}$\ti comodules in the category of $T$\ti modules in $\cc$.
Objects of $\HM(T)$ are triples $(B,r,\rho)$, where $B$ is an object of $\cc$, $r\co TB \to B$, and
$\rho\co B \to TB$ are morphisms in $\cc$, such that $(B,r)$ is a $T$\ti module, that is,
$$
rT(r)=r\mu_B \quad \text{and} \quad  r\eta_B=\id_B,
$$
and $(B,\rho)$ is a $\hat{T}$\ti comodule whose coaction is $T$\ti linear, that is,
$$
T(\rho)\rho=\delta_B\rho, \quad  r\rho=\id_B, \quad \text{and} \quad \rho r=\mu_B T(\rho).
$$
Morphisms in $\HM(T)$ from $(M,r,\rho)$ to $(N,s,\varrho)$ are morphisms $f\co M \to N$ in $\cc$ which are morphisms of $T$\ti modules and $\hat{T}$\ti comodules:
$$
fr=sT(f) \quad \text{and} \quad T(f)\rho=\varrho f.
$$
The comparison functor of the comonad $\hat{T}$ is the functor
$$
\chi \co \left\{ \begin{array}{ccl} \cc & \to &\HM(T)\\
X & \mapsto &  (TX,\mu_X, T\eta_X)
\end{array} \right..
$$
The question whether $\chi$ is an equivalence is a descent problem.

The \emph{coinvariant part} of an object $\mathbb{B}=(B,r,\rho)$  of $\HM(T)$ is the equalizer of the  coreflexive pair
$$\xymatrix{B \dar{\eta_B}{\rho}&TB}.$$
If the coinvariant part of $B$ exists, it is denoted
by  $i_\mathbb{B}\co \mathbb{B}_T \to B$.
We say that $T$ \emph{admits coinvariant parts} if any object of $\HM(T)$ admits
a coinvariant part. We say that $T$ \emph{preserves coinvariant parts} if, for any object $\mathbb{B}$ of $\HM(T)$ which admits a coinvariant part  $i_\mathbb{B}\co \mathbb{B}_T \to B$,
the morphism $T(i_\mathbb{B})$ is an equalizer of the pair $(T\eta_B, T\rho)$.

The following characterization of monads $T$ for which $\chi$ is an equivalence is a reformulation of \cite[Theorem 1]{FM1}.

\begin{thm}\label{thm-comonadicity}
Let $T$ be a monad on a category $\cc$. The following assertions are equivalent:
\begin{enumerate}
\labeli
\item The functor $\chi \co \cc \to \HM(T)$ is an equivalence of categories;
\item $T$ is conservative, admits coinvariant parts, and preserves coinvariant parts.
\end{enumerate}
If such is the case, the functor `coinvariant part'   $\mathbb{B} \mapsto \mathbb{B}_T$ is quasi-inverse to~$\chi$.
\end{thm}

\subsection{Hopf modules for \uHopf\ monads}
Let $T$ be a bimonad on a monoidal category $\cc$.
The \emph{induced coalgebra of $T$}, denoted by $\hat{C}$, is the induced coalgebra of the comonoidal adjunction $(F_T,U_T)$. Explicitly
$\hat{C}=(T(\un), \mu_\un)$, with coproduct $T_2(\un,\un)$ and counit $T_0$. Note that $U_T(\hat{C})=T(\un)$ is a coalgebra in $\cc$.

A \emph{left Hopf $T$-module} (as defined in~\cite{BV2}) is a
left $\hat{C}$-comodule in $\mo{T}{\cc}$, that is, a triple $(M, r, \rho)$ such that $(M,r)$ is a $T$\ti module,
$(M,\rho)$ is a left $T(\un)$\ti comodule, and
\begin{equation*}
\rho r= (\mu_\un \otimes r) T_2(T\un,M) T(\rho).
\end{equation*}
A \emph{morphism of Hopf $T$-modules} between two left Hopf $T$-modules $(M, r, \rho)$ and $(N, s, \varrho)$ is a morphism of
$\hat{C}$-comodules in $\mo{T}{\cc}$, that is, a morphism $f \co M \to N$ in~$\cc$ such that
\begin{equation*}
f r=s T(f) \text{\quad and \quad} (\id_{T(\un)}\otimes f) \rho= \varrho f.
\end{equation*}
We denote by $\HM^l(T)$ the category of left Hopf $T$\ti modules.

The \emph{coinvariant part of a left Hopf module $\mathbb{M}=(M,r,\rho)$} is the equalizer of the coreflexive pair
$$\xymatrix@C=4em{M\dar{\eta_\un \otimes M}{\rho}&T(\un) \otimes M}.$$
If it exists, it
is denoted by $\mathbb{M}_T$. We say that $T$ preserves coinvariant parts of left Hopf modules if, whenever a left Hopf module $\mathbb{M}=(M,r,\rho)$ admits a coinvariant part $i_T \co \mathbb{M}_T \to M$, then $T(i_T)$ is an equalizer of $(T(\eta_\un \otimes M),T\rho)$.

\begin{thm}\label{Hopfmodcor}
Let $T$ be a left  \uHopf\ monad on a monoidal category~$\cc$.
The following assertions are equivalent:
\begin{enumerate}
\labeli
\item The functor $$\mathfrak{h}^l\co \left\{ \begin{array}{ccl} \cc & \to & \HM^l(T) \\ X & \mapsto &\bigl(TX,\mu_X,T_2(\un,X) \bigr) \end{array} \right. $$ is an equivalence of categories;
\item $T$ is conservative, left Hopf $T$\ti modules admit coinvariant parts, and~$T$ preserves them.
\end{enumerate}
If these hold, the functor `coinvariant part'  $\mathbb{M} \mapsto \mathbb{M}_T$ is quasi-inverse to $\mathfrak{h}^l$.
\end{thm}

\begin{rem}
Similarly, we define the category $\HM^r(T)$ of \emph{right Hopf $T$-modules}. Since  $\HM^r(T)=\HM^l(T^\cop)$, Theorem~\ref{Hopfmodcor} holds also for right  \uHopf\ monads and right  Hopf modules (see Remark~\ref{rem-copHopf}).
\end{rem}

\begin{exa}
Let $(A,\sigma)$ be a central Hopf algebra in a monoidal category $\cc$, that is, a Hopf algebra in the center $\zz(\cc)$ of $\cc$. Consider the left Hopf monad $T=A\otimes_\sigma ?$
on~$\cc$, see Proposition~\ref{prop-Halg-to-Hmon}. A \emph{left Hopf module over $A$} is left Hopf $T$\ti module, that is, a triple $(M,r\co A\otimes M \to M, \rho\co M \to A \otimes M)$ such that $(M,r)$ is a left $A$-module,
$(M,\rho)$ is a left $A$-comodule, and $\rho r= (m \otimes r) (\id_A \otimes \sigma_A \otimes \id_M) (\Delta \otimes
\rho)$, where $m$ is the product of $A$ and $\Delta$ is coproduct of $A$. Assume now
that $\cc$ splits idempotents. Then  the morphism $r(S \otimes \id_M)\rho$ (where $S$ denotes the antipode of $A$) is an idempotent of $A \otimes M$ and its image is the coinvariant part of $M$. One verifies that $T$ is conservative and preserves coinvariant parts. Applying  Theorem~\ref{Hopfmodcor},
we obtain the fundamental theorem of Hopf modules for central Hopf algebras.
In view of Remark~\ref{rem-rep-braided}, we recover the fundamental decomposition theorem of Hopf modules for Hopf algebras in a braided category (see \cite{BKLT}) which, for the category of vector spaces over a field, is just Sweedler's classical theorem. For a detailed treatment of the case of Hopf algebras over a field,  we refer to \cite[Examples 6.2 and 6.3]{BN}.
\end{exa}

\begin{proof}[Proof of Theorem~\ref{Hopfmodcor}]
Let $\hat{T}$ be the comonad of the adjunction $(F_T,U_T)$ and $\hat{C}$ be the induced coalgebra of $T$. Since $T$ is a left \uHopf\ monad, $\FO^l_{\un,-}\co \hat{T} \to \hat{C}\otimes ?$ is an isomorphism of comonads by Lemma~\ref{lem-comon-coalg}. It induces an isomorphism of categories
$$\kappa^l_T \co \HM(T)=(\mo{T}{\cc})_{\hat{T}} \to (\mo{T}{\cc})_{\hat{C}\otimes ?}=\HM^l(T)$$ such that $\kappa^l_T  \chi=\mathfrak{h}^l$.
We conclude using Theorem~\ref{thm-comonadicity}.
\end{proof}

\subsection{Summary}\newcommand{\Ag}{\mathfrak{a}}
\newcommand{\Mg}{\mathfrak{m}}
\newcommand{\Og}{\mathfrak{o}}
\newcommand{\Cg}{\mathfrak{c}}
In this section we summarize the relationships between Hopf monads, Hopf adjunctions, and cocommutative central coalgebras.

We have constructed several correspondences between these objects:
\begin{itemize}
\item
A Hopf adjunction $(F\co \cc \to \dd, U \co \dd  \to \cc)$ gives rise to  a Hopf monad $\Mg(F,U)=UF$ on $\cc$ by Proposition~\ref{thm-adj-frob-hopf}, and to a cocommutative central coalgebra
$\Cg(F,U)=(\hat{C},\hat{\sigma})$ in $\dd$ by Corollary~\ref{cor-Hopf-Adj-Rep};

\item
A Hopf monad $T$ on a monoidal category $\cc$ defines a Hopf adjunction $$\Ag(T)=(F_T \co \cc \to \mo{T}{\cc},U_T\co \mo{T}{\cc} \to \cc)$$ by Theorem~\ref{thm-hopf-frob};

\item
 A \emph{cotensorable} cocommutative central coalgebra $(C,\sigma)$ on a monoidal category $\dd$ yields a \emph{Hopf} adjunction  $$\Og(C,\sigma)=(U \co \lsComod{\sigma}{C} \to \dd, R \co \dd \to \lsComod{\sigma}{C})$$ by Theorem~\ref{thm-coalg-hopfadj}.
\end{itemize}

Hence the following triangle:
$$\xymatrix@R=4.5em@C=1.5em{
&  **{[F]+}\txt{Hopf adjunctions}\ar@/^1pc/[dl]^-{\txt{$\Mg$}} \ar@/^1pc/[dr]^-{\txt{$\Cg$}}&\\
**{[F]+}\txt{Hopf monads} \ar@/_1pc/[rr]_-{\txt{$\Cg\Ag$}} \ar@/^1pc/[ur]^-{\txt{$\Ag$}} && **{[F]+}\txt{cocommutative\\ central coalgebras}\ar@/_1pc/@{-->}[ll]_{\txt{$\Mg\Og$}} \ar@/^1pc/@{-->}[ul]^-{\txt{$\Og$}}}$$

We have:
\begin{itemize}
\item  $\Mg \Ag(T)=T$;
\item $\Ag\Mg(F,U)\simeq (F,U)$ if and only if the adjunction $(F,U)$ is monadic;
\item  $\Cg\Og(C,\sigma)=(C,\sigma)$;
\item  assuming $\Cg(F,U)$ is cotensorable, we have $\Og\Cg(F,U) \simeq (F,U)$ if the comonoidal adjunction $(F,U)$  satisfies the conditions of Proposition~\ref{prop-HM-strong}.
\end{itemize}

With suitable exactness assumption, we have in fact equivalences:

 \begin{thm} The following data are equivalent via the assignments $\Ag$ and $\Cg$:
\begin{enumerate}[(A)]
\item A Hopf monad $T$ on a monoidal category $\cc$ such that:
    \begin{itemize}
    \item $\cc$ admits reflexive coequalizers and coreflexive equalizers, and its monoidal product preserves coreflexive equalizers;
    \item $T$ is conservative and preserves reflexive coequalizers and coreflexive equalizers;
    \end{itemize}
\item A Hopf adjunction $(F \co \cc \to \dd, U \co \dd \to \cc)$ such that:
  \begin{itemize}
  \item $\cc$ and $\dd$ admit reflexive coequalizers and coreflexive equalizers, and their monoidal products preserve coreflexive equalizers;
  \item $F$ and $U$ are conservative, $U$ preserves reflexive coequalizers and $F$ preserves coreflexive equalizers.
 \end{itemize}
  \item A cocommutative central coalgebra $(C,\sigma)$ in a monoidal category $\dd$ such that:
 \begin{itemize}
   \item $\dd$ admits reflexive coequalizers and coreflexive equalizers, and its monoidal product preserves coreflexive equalizers (in particular the central coalgebra $(C,\sigma)$ is co\-tensorable);
       \item the endofunctor $C \otimes ?$ of $\dd$ is  conservative and preserves reflexive coequalizers.
    \end{itemize}
\end{enumerate}
 Moreover,  a Hopf adjunction satisfying the conditions of (B) is a monadic and comonadic Hopf adjunction.
\end{thm}

\begin{proof} Firstly, we show the equivalence of (A) and (B).
Let $T$ be a Hopf monad on a monoidal category $\cc$ satisfying the conditions of (A). Then $U_T$, being the forgetful functor of a monad, preserves and creates limits and in particular equalizers. As a result, the monoidal category $\mo{T}{\cc}$ admits coreflexive equalizers and $U_T$ preserves them. From this one deduces that, since
the monoidal product of $\cc$ preserves coreflexive equalizers, so does that of $\mo{T}{\cc}$.
Moreover, since $T$ preserves reflexive coequalizers, $U_T$ creates and preserves them. Consequently:  $\mo{T}{\cc}$ admits
reflexive coequalizers, and $F_T$ preserves reflexive coequalizers.
The forgetful functor $U_T$ is conservative, and since by assumption $T=U_T F_T$ is conservative, so is $F_T$. Thus $\Ag(T)=(F_T,U_T)$ is a Hopf adjunction satisfying the conditions of (B), and we have $\Mg\Ag(T)=T$.
Conversely, let $(F,U)$ be a Hopf adjunction satisfying the conditions of (B).
By adjunction $F$ preserves colimits and $U$ preserves limits. The Hopf monad $T=\Mg(F,U)=UF$ is conservative and preserves reflexive coequalizers and coreflexive equalizers, so that it satisfies the conditions of (A). Moreover by Beck's monadicity theorem, the adjunction $(F,U)$ is monadic, so $\Ag\Mg(F,U)\simeq (F,U)$, hence the equivalence of (A) and (B).

Secondly, we show the equivalence of (B) and (C). Let $(C,\sigma)$ be a cocommutative central comonad in a monoidal category $\dd$ satisfying the conditions of (C). Then $(C,\sigma)$ is cotensorable, and the adjunction $$\Og(C,\sigma)=(V \co \lsComod{\sigma}{C} \to \dd, R \co \dd \to \lsComod{\sigma}{C})$$ is a Hopf adjunction. It is comonadic, with comonoidal comonad $\hat{T}=C \otimes_\sigma ?$. It is also monadic, see Remark~\ref{rem-monadic}.
Moreover the cotensor product $\otimes^\sigma_C$ preserves coreflexive equalizer by Lemma~\ref{lem-cotens}.
Thus the adjunction $(V,R)$ satisfies the conditions of (B). We have
$\Cg\Og(C,\sigma)=(C,\sigma)$.

Let us prove conversely that if $(F,U)$ is a Hopf adjunction satisfying the conditions of (B), then its induced central coalgebra
$(\hat{C},\hat{\sigma})=\Cg(F,U)$ satisfies the conditions of (C) and $\Og \Cg(F,U) \simeq (F,U)$ as Hopf adjunctions.
We will need  the following lemmas.

\begin{lem}\label{lem-OgCg1}
 Let $\cc$ be a category admitting coreflexive equalizers and let $T$ be a conservative monad on $\cc$ preserving coreflexive equalizers. Then for each object $X$ of $\cc$,
$\eta_X$ is an equalizer of the pair $(T(\eta_X),\eta_{TX})$.
\end{lem}
\begin{proof}
Let $X$  be an object of $\cc$. Observe that $T(\eta_X)$ is an equalizer of the coreflexive pair $(T^2(\eta_X), T(\eta_{TX}))$.
Since $T$ is conservative and $\cc$ admits coreflexive equalizers preserved by $T$, $\eta_X$ is an equalizer of the coreflexive pair  $(T(\eta_X),\eta_{TX})$.
\end{proof}
\begin{lem}\label{lem-OgCg2}
Let $\cc$ be a monoidal category  whose monoidal product preserves coreflexive equalizers in the left variable. Let $T$ be a left Hopf monad on $\cc$ which preserves coreflexive equalizers. Assume furthermore that for each object $X$ of~$\cc$, $\eta_X$ is an equalizer of the pair $(T(\eta_X),\eta_{TX})$.
Then for all objects $X, Y$ of $\cc$, $T_2(X,Y) \co T(X \otimes Y) \to TX \otimes TY$ is an equalizer of the coreflexive pair
$$\xymatrix@C=5em{T(X) \otimes T(Y) \dar{T_2(X,\un) \otimes T(Y)}{T(X) \otimes T_2(\un,Y)}&T(X) \otimes T(\un) \otimes T(Y)}.$$
\end{lem}
\begin{proof}  The following diagram:
$$\xymatrix@C5em{
T(X\otimes Y)\ar[r]^{T(X\otimes \eta_Y)}\ar[d]_{=} & T(X \otimes TY)\dar{T(X\otimes \eta_{TY})}{T(X \otimes T(\eta_Y))}\ar[d]_{H^l_{X,Y}} & T(X \otimes T^2Y)\ar[d]^{(TX \otimes H^l_{\un,Y})H^l_{X,TY}} \\
T(X \otimes Y)\ar[r]_{T_2(X,Y)} & TX \otimes TY\dar{T_2(X,\un) \otimes TY}{TX \otimes T_2(\un,Y)} & TX \otimes T\un \otimes TY,
 }$$
 is commutative (in the sense that the left square and the two
 right squares obtained by taking respectively the top and bottom arrow of each pair, are commutative); this results easily from Proposition~\ref{prop-fusion}.
 The top row is an equalizer because the endofunctor $T(X \otimes ?)$ preserves coreflexive equalizers.
Since $H^l$ is invertible, we conclude that the bottom row is exact, hence the lemma.
\end{proof}
Now let $T=UF$ be the Hopf monad of $(F,U)$. Then $T$ satisfies the hypotheses of Lemmas~\ref{lem-OgCg1} and~\ref{lem-OgCg2}, so that for all objects $X, Y$ of $\cc$, $T_2(X,Y) \co
T(X\otimes Y) \to T(X) \otimes T(Y)$ is an equalizer of the pair
$$\xymatrix@C=5em{T(X) \otimes T(Y) \dar{T_2(X,\un) \otimes T(Y)}{T(X) \otimes T_2(\un,Y)}&T(X) \otimes T(\un) \otimes T(Y)}.$$
Moreover, the adjunction $(F,U)$ is monadic. In particular the functor $U$ creates and preserves equalizers; thus
$F_2(X\otimes Y)$ is an equalizer of the pair
$$\xymatrix@C=5em{F(X) \otimes F(Y) \dar{F_2(X,\un) \otimes F(Y)}{F(X) \otimes F_2(\un,Y)}&F(X) \otimes F(\un) \otimes F(Y)}.$$
We may therefore apply Proposition~\ref{prop-HM-strong} to the adjunction $(F,U)$, and we conclude that the comparison functor $\cc \to \lsComod{\hat{\sigma}}{\hat{C}}$ is a strong monoidal equivalence and $\Cg(F,U)$ satisfies the conditions of (C), hence $\Og \Cg(F,U) \simeq (F,U)$ as Hopf adjunctions.
\end{proof}

\section{Hopf algebroids and finite abelian tensor categories}\label{sect-hopf-oid}

In this section, we study  bialgebroids which, according to Szlach\'anyi~\cite{Szl03}, are linear bimonads on categories of bimodules admitting a right adjoint. A bialgebroid corresponds with a Hopf monad if and only if it is a Hopf algebroid in the sense of Schauenburg~\cite{Sch}. We also use Hopf monads to prove that any finite tensor category is naturally equivalent (as a tensor category) to the category of finite-dimensional modules over some finite dimensional Hopf algebroid.

\subsection{Bialgebroids and bimonads}

Let $\kk$ be a commutative ring and $R$ be a \kt algebra. Denote by $\biMod{R}{R}$ the category of $R$\ti bimodules. It is a monoidal category, with monoidal product $\otimes_R$ and unit object $R$.  We identify $\biMod{R}{R}$ with the category $\lMod{R^e}$ of left $R^e$\ti modules, where $R^e=R \otimes_\kk R^\opp$. Hence a monoidal product~$\biotimes$ on  $\lMod{R^e}$ (corresponding to $\otimes_R$ on $\biMod{R}{R}$), with unit $R$ (whose $R^e$\ti action is $(r \otimes r') \cdot x=rxr'$).

If $f\co B \to A$ is \kt algebra morphism, we denote by $\lmo{f}{A}$ the left $B$\ti module $A$ with left action $b\cdot a=f(b)a$, and by $A_f$ the right $B$\ti module $A$ with right action $a\cdot b=af(b)$.

A \emph{left bialgebroid with base $R$} (also called Takeuchi $\times_R$-bialgebra) consists of data  $(A,s,t,\Delta,\varepsilon)$ where:
\begin{itemize}
\item $A$ is a \kt algebra;
\item $s\co R \to A$ and $t\co R^\opp \to A$ are \kt algebra morphisms whose images in $A$ commute. Hence a \kt algebra morphism $$ e \co \left\{ \begin{array}{ccc} R^e & \to & A \\ r \otimes r' & \mapsto &s(r)t(r') \end{array} \right.,$$ which gives rise to a $R^e$\ti bimodule $\bi{e}{A}{e}$.

\item $(\lmo{e}{A},\Delta,\varepsilon)$ is a coalgebra in the monoidal category $(\lMod{R^e}, \biotimes,R)$.
\end{itemize}
In this situation the Takeuchi product $A\times_R A  \subset \lmo{e}{A} \biotimes \lmo{e}{A}$, defined by
$$A\times_R A =\{{\textstyle\sum a_i \otimes b_i \in \lmo{e}{A} \biotimes \lmo{e}{A}  \mid \forall r \in R, \sum a_it(r)  \otimes b_i=\sum a_i \otimes b_i s(r)}\}$$
is a \kt algebra, with product defined by  $(a\otimes b)(a' \otimes b')=aa'\otimes b b'$,  and
one requires:
\begin{itemize}
\item $\Delta(A) \subset A\times_R A$;
\item $\Delta \co A \to A\times_R A$ is a \kt algebra morphism;
\item $\varepsilon(a \,s(\varepsilon(a')))=\varepsilon(aa')=\varepsilon(a \,t(\varepsilon(a')))$;
\item $\varepsilon(1_A)=1_R$.
\end{itemize}

The notion of left bialgebroid has a nice interpretation in terms of bimonads. A bialgebroid $A$ with base $R$ gives rise to an endofunctor of $\lMod{R^e}\simeq \biMod{R}{R}$:
$$
T_A\co \left\{ \begin{array}{ccl} \lMod{R^e} & \to & \lMod{R^e} \\ M & \mapsto &T_A(M)= \bi{e}{A}{e} \otimes_{R^e} M \end{array} \right.
$$
The axioms of a left bialgebroid are such that $T_A$ is a \kt linear  bimonad admitting a right adjoint.    These properties characterize left bialgebroids:
\begin{thm}[\cite{Szl03}]\label{thm-szl} Let $\kk$ be a ring and $R$ a $\kk$\ti algebra.
Via the correspondence $A \mapsto T_A$, the following data are equivalent:
\begin{enumerate}[(A)]
\item A left bialgebroid $A$ with base $R$;
\item A $\kk$\ti linear bimonad $T$ on the monoidal category $\biMod{R}{R} \simeq \lMod{R^e}$ admitting a right adjoint.
\end{enumerate}
\end{thm}

\subsection{Hopf algebroids}
We define a \emph{left}, resp.\@ \emph{right}, \emph{(pre\ti)Hopf algebroid} to be a bialgebroid $A$ whose associated bimonad $T_A$ is a left, resp.\@ right, (pre\ti)Hopf monad. A \emph{(pre\ti)Hopf algebroid} is a left and right (pre\ti)Hopf algebroid.

Let  $A$ be a bialgebroid and $T_A$ be its associated bimonad on  $\biMod{R}{R} \simeq \lMod{R^e}$.
The fusion operators $H^l$ and $H^r$ of $T_A$ are:
$$H^l_{M,N} \co \left\{ \begin{array}{ccc} \bi{e}{A}{e} \otimes_{R^e}
\bigl( M  \biotimes (\bi{e}{A}{e} \otimes_{R^e} N) \bigr)
 &\to &
(\bi{e}{A}{e} \otimes_{R^e} M) \biotimes
(\bi{e}{A}{e} \otimes_{R^e} N)\\
a \otimes m \otimes b \otimes n &\mapsto & a_{(1)} \otimes m \otimes a_{(2)}b \otimes n\end{array} \right.
$$
and
$$H^r_{M,N} \co \left\{ \begin{array}{ccc} \bi{e}{A}{e} \otimes_{R^e}
\bigl((\bi{e}{A}{e} \otimes_{R^e} M ) \biotimes  N \bigr)
 &\to &
(\bi{e}{A}{e} \otimes_{R^e} M) \biotimes
(\bi{e}{A}{e} \otimes_{R^e} N)\\
a \otimes b \otimes m  \otimes n &\mapsto & a_{(1)}b \otimes m \otimes a_{(2)} \otimes n\end{array} \right. .
$$
Using the fact that $R^e$ and $R$ are respectively a projective generator and
the unit object of $\lMod{R^e}$,
we obtain the following characterization of Hopf bialgebroids and \uHopf\ algebroids.

\begin{prop}
Let $A$ be a bialgebroid with base $R$. Then:
\begin{enumerate}
\labela
\item The bialgebroid $A$ is a left Hopf algebroid if and only if  the $R^e$\ti linear map
$$
H^l_{R^e,R^e}\co \left\{ \begin{array}{ccc} \bi{e}{A}{t} \otimes_{R^\opp} {\lmo{t}{A}} & \to & \lmo{e}{A} \biotimes \lmo{e}{A}  \\ a \otimes b  & \mapsto &a_{(1)} \otimes a_{(2)} b\end{array} \right.
$$
is bijective.
\item The bialgebroid $A$ is a right Hopf algebroid if and only if  the $R^e$\ti linear map
$$
H^r_{R^e,R^e}\co \left\{ \begin{array}{ccc}  \bi{e}{A}{s} \otimes_R  \lmo{s}{A} & \to &  \lmo{e}{A} \biotimes \lmo{e}{A} \\ a \otimes b  & \mapsto & a_{(1)} b\otimes a_{(2)}\end{array} \right.
$$
is bijective.
\item  The bialgebroid $A$ is a left \uHopf\ algebroid if and only if the $R^e$\ti linear map
$$
H^l_{R,R^e} \co \left\{ \begin{array}{ccc} \bi{e}{A}{e} \otimes_{R^e} \lmo{e}{A} & \to & \lmo{e}{\Bar{A}} \biotimes \lmo{e}{A} \\ a \otimes a'  & \mapsto &a_{(1)} \otimes a_{(2)} a'\end{array} \right.
$$
is bijective, where $$\lmo{e}{\Bar{A}}=\bi{e}{A}{e}\otimes_{R^e} R=A/\{as(r)=at(r)\mid a \in A, r \in R\}.$$
\item  The bialgebroid $A$ is a right \uHopf\ algebroid if and only if the $R^e$\ti linear map
and
$$
H^r_{R^e,R} \co \left\{ \begin{array}{ccc} \bi{e}{A}{e} \otimes_{R^e} \lmo{e}{A}   & \to & \lmo{e}{A} \biotimes \lmo{e}{\Bar{A}} \\ a \otimes a'  & \mapsto & a_{(1)} a'\otimes a_{(2)} \end{array} \right.
$$
is bijective.
\end{enumerate}
\end{prop}

\begin{rem}
The notion of $\times_R$-Hopf algebra introduced by Schauenburg in \cite{Sch} coincides with our notion of left Hopf algebroid.
\end{rem}

\begin{rem}
The category $\lMod{R^e}$ is monoidal closed with internal Homs:
$$
\lhom{M,N}=\Hom_{R^\opp}(\lmo{R^e}{M},\lmo{R^e}{N}) \quad\text{and}\quad \rhom{M,N}=\Hom_{R}(\lmo{R^e}{M},\lmo{R^e}{N}).
$$
By Theorem~\ref{thm-hopfmon-closed}, a left bialgebroid $A$ with base $R$ is a Hopf algebroid if and only if it admits a left antipode
$$
s^l_{M,N}\co \bi{e}{A}{e} \biotimes \Hom_{R^\opp}(\bi{e}{A}{e} \otimes_{R^e} M,N) \to \Hom_{R^\opp}(M,\bi{e}{A}{e} \otimes_{R^e} N)
$$
and a right antipode
$$
s^r_{M,N}\co \bi{e}{A}{e} \biotimes \Hom_R(\bi{e}{A}{e} \otimes_{R^e} M,N)\to \Hom_R(M,\bi{e}{A}{e} \otimes_{R^e} N).
$$
\end{rem}

\begin{rem}
Let $A$ be a \uHopf\ algebroid with base $R$. Since $\lMod{R^e}$ is abelian, the bimonad $T_A$ admits coinvariant parts. If $A_e$ is a faithfully flat right $R^e$\ti module, then $T_A$ is conservative and preserves coinvariant parts.  Thus, the Hopf module decomposition theorem (see Theorem~\ref{Hopfmodcor}) applies to (pre\ti)Hopf algebroids which are faithfully flat on the right over the base ring.
\end{rem}

\subsection{Existence of fibre functors for finite tensor categories}
A \emph{tensor category over $\kk$} is an autonomous category endowed with a structure of $\kk$-linear abelian category such
that the monoidal product $\otimes$ is bilinear and $\End(\un)=\kk$.

We say that a \kt linear abelian category  $\aaa$ is \emph{finite} if it is \kt linearly equivalent to the category $\bimod{R}{}$ of finite-dimensional left modules over some finite-dimensional \kt algebra $R$.
Note that if $\aaa$ is a finite, then so is $\aaa^\opp$, since the functor
$$
\left\{ \begin{array}{ccc} (\bimod{R}{})^\opp  & \to &  \bimod{R^\opp}{} \\
N & \mapsto & \Hom(N,\kk) \end{array} \right.
$$
is a \kt linear equivalence.

\begin{thm}\label{thm-tannaka}
Let $\cc$ be a finite tensor category over a field $\kk$. Then $\cc$ is equivalent, as a tensor category, to the category of modules over a finite-dimensional left Hopf algebroid over $\kk$.
\end{thm}

We  first  state and prove an analogue of this theorem in terms of Hopf monads in a somewhat more general setting. 
Let $\cc$ be a monoidal category. Recall that the category $\EndNat(\cc)$ of endofunctors of $\cc$ is strict monoidal with composition for monoidal product and $1_\cc$ for monoidal unit.
The functor
$$
\Omega \co \left\{ \begin{array}{ccc} \cc &\to & \EndNat(\cc) \\
X &\mapsto & X \otimes ? \end{array} \right.
$$
is strong monoidal.

\begin{thm}\label{lem-adj-mon} Let $\ee$ be a full monoidal subcategory of $\EndNat(\cc)$ containing $\Omega(\cc)$. Denote by
$\omega \co \cc \to \ee$
the corestriction of $\Omega$ to $\ee$.
Then
\begin{enumerate} \labela
\item If $\omega$ has a left adjoint $\ff$, the adjunction $(\ff,\omega)$ is monadic, its monad $T=\omega \ff$ is a bimonad on $\ee$, and the comparison functor $\cc \to \ee^T$ is a monoidal equivalence.
\item If $\cc$ is right autonomous, then $\omega$ has a left adjoint if and only if  the coend $$\ff(e)=\int^{X \in \cc} e(X) \otimes \rdual{X}$$ exists for all $e \in \Ob(\ee)$. In that case,  the assignment  $e \mapsto \ff(e)$ defines a functor which is a left adjoint of $\omega$,  and the bimonad $T=\omega \ff$ is a right Hopf monad.
\item If $\cc$ is autonomous and $\omega$ has a left adjoint $\ff$, then the bimonad $T=\omega \ff$ is a Hopf monad.
\end{enumerate}
\end{thm}

\begin{proof} Assume $\omega$ has a left adjoint $\ff$. Then the adjunction $(\ff,\omega)$ is a comonoidal adjunction, so that the comparison functor $K \co \cc \to \mo{T}{\ee}$ is strong monoidal. Besides, $\omega$ has a left quasi-inverse $e \mapsto e(\un)$, and so satisfies conditions (a) and (b) of Beck's monadicity Theorem~\ref{thm-beck}, so that the adjunction $(\ff,\omega)$ is monadic and $K$ is a monoidal equivalence. Hence Part (a).

Assume $\cc$ is right autonomous.
For $e \in \Ob(\ee)$ and $X \in \Ob(\cc)$, we have a natural bijection between natural transformations $e \to X \otimes ?$ and
dinatural transformations $\{e(Y)\otimes \rdual{Y} \to X\}_{Y \in \Ob(\cc)}$. Therefore $\omega$ has a right adjoint if and only if the coends $\ff(e)$ exist for any object $e$ of $\ee$. Assume that such is the case. Then  the assignment $e \mapsto  \ff(e)$  gives a left adjoint of $\omega$.  For $X \in \Ob(\cc)$ and $e \in \Ob(\ee)$, we have:
$$\ff(\omega(X) \circ e)= \int^{Y \in \cc} X \otimes e(Y) \otimes \rdual{Y} \simeq X \otimes \int^{Y\in \cc} e(Y) \otimes \rdual{Y} =X \otimes \ff e$$ because $X \otimes ?$ has a right adjoint $\rdual{X} \otimes ?$ and so preserves colimits. One checks that this isomorphism is the right \frob\ operator $\FO^r_{e,X} \co \ff(\omega(X) \circ e) \to X \otimes \ff e $
of the adjunction $(\ff,\omega)$. Thus $T$ is a right Hopf monad by Theorem~\ref{thm-hopf-frob}.  Hence Part~(b).

Finally  assume that $\cc$ is also left autonomous.  Let $X \in \Ob(\cc)$ and $e \in \Ob(\ee)$.  Since the functor $\ldual{X} \otimes ?$ is left adjoint to $X \otimes ?$ and the functor $? \otimes X$ preserves colimits (because it has a right adjoint $? \otimes \ldual{X}$), we have:
\begin{align*}
\ff(e \circ \omega(X)) &= \int^{Y \in \cc}e(X \otimes Y) \otimes \rdual{Y} \simeq \int^{Y \in \cc} e(Y) \otimes \rdual{(\ldual{X} \otimes Y)}\\
&\simeq \int^{Y \in \cc} e(Y) \otimes \rdual{Y} \otimes X \simeq \ff e \otimes X.
\end{align*}
One checks that the composition of these isomorphisms is the left \frob\ operator $\FO^l_{X,e} \co \ff(e \circ \omega(X)) \to \ff e  \otimes X$ of the adjunction $(\ff,\omega)$.  Therefore $T$ is also a left Hopf monad  by Theorem~\ref{thm-hopf-frob}. Hence Part (c).
\end{proof}

\begin{proof}[Proof of Theorem~\ref{thm-tannaka}]
We apply Theorem~\ref{lem-adj-mon} to a finite tensor category $\cc$ over a field $\kk$.
If $\aaa$ is a  $\kk$\ti linear abelian category, we denote by $\EndNat^{ra}_\kk(\aaa)$ the full monoidal
subcategory of $\EndNat(\aaa)$ of $\kk$\ti linear endofunctors which admit a right
adjoint.

Set $\ee=\EndNat^{ra}_\kk(\cc)$. For $X \in \Ob(\cc)$, the endofunctor
$X \otimes ?$ is $\kk$\ti linear and has a right adjoint, namely $\rdual{X} \otimes ?$, so we have $\Omega(\cc) \subset \ee$.
Denoting by $\omega \co \cc \to \ee$ the corestriction of $\Omega$ to $\ee$, we have a commutative triangle of strong monoidal $\kk$\ti linear functors:
$$\xymatrix@R=1em@C=1.8em{
\cc \ar[rr]^-{\Omega}\ar[rd]_-{\omega}& & \EndNat(\cc)\\
& \ee\ar[ru]_{\text{ inclusion}} &
}$$

By assumption, there exists a finite dimensional $\kk$\ti algebra $R$ and a $\kk$\ti linear equivalence
$\Upsilon \co\cc \to \bimod{R}{}$, with  quasi-inverse $\Upsilon^*$ of $\Upsilon$, hence a $\kk$\ti linear strong monoidal equivalence:
$$
\left\{ \begin{array}{ccc} \ee=\EndNat^{ra}_\kk(\cc) &\to &\EndNat^{ra}_\kk(\bimod{R}{})\\
 E &\mapsto &\Upsilon\circ  E \circ \Upsilon^*\end{array} \right.
$$
Composing this with the well-known strong monoidal  $\kk$\ti linear equivalence
$$
\left\{ \begin{array}{ccc}\EndNat^{ra}_\kk(\bimod{R}{}) &\to &\bimod{R}{R}\\
 e &\mapsto &{e(_RR)}_R \end{array} \right.
$$
we obtain
a  $\kk$\ti linear strong monoidal equivalence
$$\Theta\co \ee \to \bimod{R}{R} \simeq \bimod{R^e}{}\,.$$
In particular $\ee$ is a finite $\kk$\ti linear abelian category.
The category $\EndNat(\cc)$ is abelian as a category of functors to an abelian
category, $\Omega$ is exact (the tensor product of~$\cc$ being exact in each variable), and the inclusion
$\ee \hookrightarrow \End(\cc)$ is fully faithful, so $\omega$ is exact.
It is a well-known fact that a right (resp.\@ left) exact $\kk$\ti linear functor between finite $\kk$\ti linear abelian categories admits a left (resp.\@ right) adjoint.
 Thus  $\omega$ has a left adjoint $\ff$, as well as a right adjoint $\rr$. By Theorem~\ref{lem-adj-mon}, we conclude that the comonoidal adjunction $(\ff,\omega)$ is monadic and its monad $T= \omega\ff$ is a Hopf monad. Moreover $T$ is $\kk$\ti linear and has a right adjoint $\omega\rr$.

Now we transport $T$ along the  $\kk$\ti linear monoidal equivalence $\Theta \co \ee \to \bimod{R}{R}$.
Pick a quasi-inverse  $\Theta^*$ of $\Theta$. The adjunction $(\ff\Theta^*,\Theta\omega)$
is a monadic \frob\ adjunction. Its monad  $T'$  is a $\kk$\ti linear Hopf monad on $\bimod{R}{R}$ with a right adjoint $\Theta \omega \rr\Theta^*$.
By Theorem~\ref{thm-szl}, $T$ is of the form $T_A$ for some bialgebroid $A$  with  base~$R$, which is by definition a Hopf algebroid. Monadicity ensures that
the comparison functor
$\cc \to \mo{T}{(\bimod{R}{R})}=\bimod{A}{}$ is a $\kk$\ti linear monoidal equivalence of categories.
This concludes the proof of Theorem~\ref{thm-tannaka}.
\end{proof}

\end{document}